\input amstex
\input epsf
\documentstyle{amsppt}

\def\al{\alpha}
\def\th{\theta}
\def\be{\beta}
\def\ep{\varepsilon}
\def\de{{\delta}}
\def\ka{\kappa}
\def\la{\lambda}
\def\La{\Lambda}
\def\ga{\gamma}
\def\Ga{\Gamma}
\def\si{\sigma}
\def\Si{\Sigma}
\def\de{{\delta}}

\def\th{\theta}

\def\cA{\Cal A}
\def\cB{\Cal B}
\def\cC{\Cal C}
\def\cD{\Cal D}
\def\cE{\Cal E}
\def\cG{\Cal G}
\def\cF{\Cal F}

\def\cI{\Cal I}
\def\cJ{\Cal J}
\def\cK{\Cal K}
\def\cL{\Cal L}
\def\cM{\Cal M}
\def\cN{\Cal N}
\def\cO{\Cal O}
\def\cR{\Cal R}

\def\cL{\Cal L}
\def\cP{\Cal P}
\def\cQ{\Cal Q}
\def\cS{\Cal S}

\def\limt{\lim_{t\to\infty}}
\def\limn{\lim_{n\to\infty}}

\def\wt{\widetilde}
\def\wh{\widehat}
\redefine\l{\ell}

\def\qqquad{\quad \quad \quad}
\def\sobre#1#2{\lower 1ex \hbox{ $#1 \atop #2 $ } }
\def\bajo#1#2{\raise 1ex \hbox{ $#1 \atop #2 $ } }

\def\qedsymbol{\hfill $\blacksquare$}
\def\today{\rightline{\ifcase\month\or
  January\or February\or March\or April\or May\or June\or
  July\or August\or September\or October\or November\or December\fi
  \space\number\day,\space\number\year}}
\overfullrule=0pt

\global\newcount\numsec\global\newcount\numfor
\global\newcount\numfig
\global\newcount\numcon
\global\newcount\numkon
\gdef\profonditastruttura{\dp\strutbox}
\def\senondefinito#1{\expandafter\ifx\csname#1\endcsname\relax}
\def\SIA #1,#2,#3 {\senondefinito{#1#2}
\expandafter\xdef\csname #1#2\endcsname{#3}\else
\write16{???? ma #1,#2 e' gia' stato definito !!!!} \fi}

\def\etichetta(#1){(\veroparagrafo.\veraformula)
\SIA e,#1,(\veroparagrafo.\veraformula)
 \global\advance\numfor by 1
 \write16{ EQ \equ(#1) #1  }}

\def\letichetta(#1){\veroparagrafo.\veraformula
\SIA e,#1,{\veroparagrafo.\veraformula}
\global\advance\numfor by 1
 \write16{ Sta \equ(#1) == #1}}

\def\tetichetta(#1){\veroparagrafo.\veraformula 
\SIA e,#1,{(\veroparagrafo.\veraformula)}
\global\advance\numfor by 1
 \write16{ tag \equ(#1) == #1}}

\def \FU(#1)#2{\SIA fu,#1,#2 }

\def\etichettaa(#1){(A\veroparagrafo.\veraformula)
 \SIA e,#1,(A.\veroparagrafo.\veraformula)
\global\advance\numfor by 1
 \write16{ EQ \equ(#1) == #1  }}

\def\getichetta(#1){Fig. \verafigura
 \SIA e,#1,{\verafigura}
 \global\advance\numfig by 1
 \write16{ Fig. \equ(#1) ha simbolo  #1  }}

\def\cetichetta(#1){C_\veracon
 \SIA c,#1,{\veracon}
 \global\advance\numcon by 1
 \write16{ Const \cequ(#1) =  #1 }}

\def\ketichetta(#1){K_{\verakon}
 \SIA k,#1,{\verakon}
 \global\advance\numkon by 1
 \write16{ Const \kequ(#1) =  #1 }}

\newdimen\gwidth
\def\BOZZA{
\def\alato(##1){
 {\vtop to \profonditastruttura{\baselineskip
 \profonditastruttura\vss
 \rlap{\kern-\hsize\kern-1.2truecm{$\scriptstyle##1$}}}}}
\def\galato(##1){ \gwidth=\hsize \divide\gwidth by 2
 {\vtop to \profonditastruttura{\baselineskip
 \profonditastruttura\vss
 \rlap{\kern-\gwidth\kern-1.2truecm{$\scriptstyle##1$}}}}}
}
\def\alato(#1){}
\def\galato(#1){}
\def\calato(#1){}
\def\kalato(#1){}
\def\veroparagrafo{\number\numsec}
\def\veraformula{\number\numfor}
\def\verafigura{\number\numfig}
\def\veracon{\number\numcon}
\def\verakon{\number\numkon}
\def\geq(#1){\getichetta(#1)\galato(#1)}
\def\Eq(#1){\eqno{\etichetta(#1)\alato(#1)}}
\def\Ceq(#1){\cetichetta(#1)\calato(#1)} 
\def\Keq(#1){\ketichetta(#1)\kalato(#1)} 
\def\teq(#1){\tag{\tetichetta(#1)\hskip-1.6truemm\alato(#1)}}  
\def\eq(#1){\etichetta(#1)\alato(#1)}
\def\Eqa(#1){\eqno{\etichettaa(#1)\alato(#1)}}
\def\eqa(#1){\etichettaa(#1)\alato(#1)}
\def\eqv(#1){\senondefinito{fu#1}$\clubsuit$#1\else\csname fu#1\endcsname\fi}
\def\equ(#1){\senondefinito{e#1}\eqv(#1)\else\csname e#1\endcsname\fi}
\def\cequ(#1){C_\senondefinito{c#1}\eqv(#1)\else\csname c#1\endcsname\fi}
\def\kequ(#1){K_\senondefinito{k#1}\eqv(#1)\else\csname k#1\endcsname\fi}
\def\Lemma(#1){Lemma \letichetta(#1)\hskip-1.6truemm}
\def\Theorem(#1){{Theorem \letichetta(#1)}\hskip-1.6truemm}
\def\Proposition(#1){{Proposition \letichetta(#1)}\hskip-1.6truemm}
\def\Corollary(#1){{Corollary \letichetta(#1)}\hskip-1.6truemm.}
\def\Remark(#1){{\noindent{\bf Remark \letichetta(#1)\hskip-1.6truemm.}}}

\def\include#1{
\openin13=#1.aux \ifeof13 \relax \else
\input #1.aux \closein13 \fi}
\openin14=\jobname.aux \ifeof14 \relax \else
\input \jobname.aux \closein14 \fi

\magnification=\magstep1
\vcorrection{.5truecm}
\vsize 19cm

\TagsOnRight
\NoBlackBoxes
\document
\topmatter
\title
A shape theorem for the spread of an infection
\endtitle
\author
Harry Kesten and Vladas Sidoravicius
\endauthor
\leftheadtext{Harry Kesten and Vladas Sidoravicius}
\rightheadtext{Shape theorem for spread of an infection}

\abstract
In \cite{KSb} we studied the following model for the spread of a rumor
or infection: There is a
``gas'' of so-called $A$-particles,
each of which performs a continuous time simple
random walk  on $\Bbb Z^d$, with jumprate $D_A$.
We assume that ``just before the start''  the number of $A$-particles at $x$,
$N_A(x,0-)$,
has a mean $\mu_A$ Poisson distribution and that the $N_A(x,0-), \, x
\in \Bbb Z^d$, are independent.
In addition, there are
$B$-particles which perform continuous time simple random walks with
jumprate $D_B$. We start with a finite number of $B$-particles
in the system at time 0. The positions of these initial $B$-particles
are arbitrary, but they are non-random. The $B$-particles move
independently of each other. The only interaction is that when a
$B$-particle and an $A$-particle coincide, the latter instantaneously
turns into a $B$-particle. \cite {KSb} gave some basic estimates for
the growth of the set $\wt B(t):= \{x \in \Bbb Z^d:$ a $B$-particle visits
$x$ during $[0,t]$\}. In this article we show that if $D_A=D_B$, then
$B(t) = \wt B(t) + [-\frac 12, \frac 12]^d$ grows linearly in time with
an asymptotic shape, i.e., there
exists a non-random set $B_0$ such that $(1/t)B(t) \to B_0$, in a sense
which will be made precise.
\endabstract

\address
Harry Kesten,
Department of Mathematics,
Malott Hall,
Cornell University,
Ithaca NY 14853, USA
\endaddress
\email
kesten\@math.cornell.edu
\endemail

\address
Vladas Sidoravicius,
IMPA,
Estr. Dona Castorina 110,
Rio de Janeiro,
Brasil,
\endaddress
\email
vladas\@impa.br
\endemail

\endkeywords

\endtopmatter

\subhead
1. Introduction
\endsubhead
\numsec=1
\numfor=1

We study the model described in the abstract. One interpretation of
this model is that the $B$-particles represent
individuals who are infected, and the $A$-particles represent
susceptible individuals; see \cite {KSb} for another interpretation.
$\wt B(t)$ represents the collection of sites which have been
visited by a $B$-particle during $[0,t]$, and $B(t)$ is a slightly
fattened up version of $\wt B(t)$, obtained by adding a unit cube
around each point of $\wt B(t)$. This fattened up version is
introduced merely to simplify the statement of our main result. It is simpler
to speak of the shape of the set $(1/t)B(t)$ as a subset of $\Bbb
R^d$, than of the discrete set $(1/t)\wt B(t)$.

The aim of this paper is to describe how the infection spreads
throughout space as time goes on. In \cite {KSb} we proved a first
result in this direction in the case $D_A=D_B$. We proved that under
this condition there exist constants $0 < C_2 \le C_1 <
\infty$ such that almost surely
$$
\cC(C_2t)  \subset B(t) \subset \cC(2C_1t)
\text{ for all large }t,
\teq(1.1)
$$
where
$$
\cC(r) := [-r,r]^d.   
\teq(1.2)
$$
\equ(1.1) gives upper and lower bounds which are linear in time,
for $B(t)$, the region which
has been visited by the infection during $[0,t]$. However, the upper
and lower bounds in \equ(1.1) are not the same. The principal result
of this paper is a so-called shape theorem which gives the first
order asymptotic behavior of the region $B(t)$. It shows that $(1/t)B(t)$
converges to a fixed set $B_0$. Thus, not only is the growth linear
in time, but $B(t)$ looks asymptotically like (a scaled version of)
$B_0$. This of course sharpens \equ(1.1) by `bringing the upper and
lower bound together'. However, the result \equ(1.1) is a crucial
tool for proving the shape theorem. We do not know of a shortcut which
proves the shape theorem without much of the development of \cite
{KSb} for \equ(1.1). The precise form of the shape theorem here is as follows:
\proclaim{Theorem 1}
Consider the model described in the abstract. If $D_A = D_B$,
then there exists a non-random, compact, convex set $B_0$ such
that for all $\ep > 0$ almost surely
$$
(1-\ep)B_0 \subset \frac 1t B(t) \subset (1+\ep)B_0 \text{ for all
large }t.
\teq(1.3)
$$
The origin is an interior point of $B_0$, and $B_0$ is invariant under
reflections in coordinate hyperplanes and under permutations of the
coordinates.
\endproclaim

\noindent
{\bf Remark 1.} It follows immediately from Theorem 1 and
Proposition B below that the particle distribution at a large time $t$
looks as follows: The numbers of particles, irrespective of type,
that is $N_A(x,t) + N_B(x,t), x \in \Bbb Z^d$, is a collection of
i.i.d. mean $\mu_A$ Poisson variables plus a finite number of particles
which started at time zero at fixed locations (these are the particles
added as $B$-particles at the start). For every $\ep >0$ there are
almost surely no $A$ particles in $(1-\ep)tB_0$ and no $B$-particles
outside $(1+\ep)t B_0$ for all large $t$.

\bigskip
Shape theorems have a fairly long history and have become the first
goal of many investigations of stochastic growth models. To the best
of our knowledge Eden (see \cite {E}) was the first one to ask for a
shape theorem for his celebrated `Eden model'. The problem turned out
to be a stubborn one. The first real progress was due to Richardson, who
proved in \cite {Ri} a shape theorem not only for the Eden model, but
also for a more general class of models, now called Richardson
models. In these models one typically thinks of the sites of $\Bbb
Z^d$ as cells which can be of two types (for instance $B$ and $A$
or infected and susceptible).
Cells can change their type to the type
of one of their neighbors according to appropriate rules. One starts
with all cells off the origin type $A$ and cell of type $B$ at the origin and
tries to prove a shape theorem for the set of cells of type $B$ at a large
time. An important example of such a model is `first-passage
percolation', which was introduced in \cite {HW} (this includes the
Eden model, up to a time change). A quite good shape theorem for
first-passage percolation is known (see \cite {Ki}, \cite {CD}, \cite
{Ke}). In more recent first-passage percolation papers even sharper
information has been obtained which gives estimates on the rate at
which $(1/t)B(t)$ converges to its limit $B_0$ (see \cite {Ho} for a
survey of such results).

Shape theroems for quite a few variations of Richardson's model and
first-passage percolation have been proven (see for instance
\cite{BG} and \cite{GM}), but as far as we know these are all for
models in which the cells do not move over time, with one
exception. This exception is the so-called frog model which follows
the rules given in our abstract, but which has $D_A = 0$, i.e., the
susceptibles or type $A$ cells stand still (see \cite {AMP} and \cite
{RS} for this model). The present paper may be the first one which
allows both tyes of particles to move.

In nearly all cases shape theorems are proven by means of Kingman's
subadditive ergodic theorem (see \cite {Ki}). This is also what is
used for the frog model. For this model one can show that the family
of random variables $\{T_{x,y}\}$ is subadditive, were $T_{x,y}$ is a
version of the
first time a particle at $y$ is infected, if one starts with one
infected particle at $x$ and one susceptible at each other site.
More precisely, the $T_{x,y}$ can all be defined on one probability
space such that $T_{x,z} \le T_{x,y} + T_{y,z}$ for all $x,y,z \in
\Bbb Z^d$ and such that their joint distribution is invariant under
translations.
Unfortunately this subadditivity property is no longer valid if one
allows both types of particles to move. Nevertheless, subadditivity
methods are still heavily used in the proof of Theorem 1. However, we now
use subadditivity only for certain `half-space' processes which
approximate the true process. Moreover, these half-space processes have only
approximate superconvolutive properties (in the terminology of \cite
{Ha}). There is no obvious family of random variables with
properties like those of the $T_{x,y}$. One only has some relation between
the distribution functions of the $H(t,u)$ for a fixed unit vector
$u$, where $H(t,u)$ is basically the maximum of $\langle x,u \rangle$
over all $x$ which have been reached by a $B$-particle by time $t$
($\langle x,u \rangle$ is the inner product of $x$ and $u$; for technical
reasons $H(t,u)$ will be calculated in a process in which the starting
conditions are slightly different from our original process).
These properties are strong enough
to show that for each unit vector $u$ there exists a constant $\la(u)$
such that almost surely
$$
\limn \frac 1t H(t,u) = \la(u),
\teq(1.4)
$$
Thus the $B$-particles reach in time $t$ half-spaces in a fixed
direction $u$ at distances which grow linearly in $t$.
Except in dimension 1, it then still requires a considerable amount of
technical work to go from this result about the linear growth
of the distances of reached half-spaces to the full asymptotic shape result.
We will give more heuristics before some of our lemmas.

\medskip
\noindent
{\bf Remark 2.} Our proof in \cite {KSb} shows that the right hand
inclusion in  \equ(1.1) remains valid
for arbitrary jumprates of the $A$ and the $B$-particles.
However, it is still not known whether the left hand inclusion holds in
general. The lower bound for $B(t)$ is known only when $D_A=D_B$,
or when $D_A=0$, that  is, when the $A$ and $B$-particles move
according to the same random walk (see \cite{KSb}),
or in the frog model, when the $A$-particles stand
still (see \cite {AMP},\cite{RS}).

\medskip
Here is some general notation which will be used throughout the paper.
$\|x\|$ without
subscript denotes the $\l^\infty$-norm of a vector $x = (x(1), \dots,
x(d)) \in \Bbb R^d$, i.e.,
$$
\|x\| = \max_{1 \le i \le d} |x(i)|.
$$
We will also use the Euclidean norm of $x$; this will be denoted by
the usual $\|x\|_2$.
$\langle x,u \rangle$ denotes the (Euclidean) inner product of two
vectors $x, u \in \Bbb R^d$, and $\bold 0$ denotes the origin (in
$\Bbb Z^d$ or $\Bbb R^d$).
For an event $\cE, \;\cE^c$ denotes its complement.

$K_1, K_2, \dots$ will denote various strictly positive, finite
constants whose precise value is
of no importance to us. The same symbol $K_i$ may have different
values in different formulae. Further, $C_i$ denotes a strictly positive
constant whose value remains the same throughout this paper.
a.s. is an abbreviation of almost surely.

\bigskip
\noindent
{\bf Acknowledgement.} The research for this paper was started during
 a stay by H. Kesten at the Mittag-Leffler Inst. in 2001-2002.
H. Kesten thanks the
 Swedish Research Council for awarding
 him a Tage Erlander Professorship for 2002.
  Further support for HK
 came from the NSF under Grant DMS 9970943 and from Eurandom. HK
 thanks Eurandom for appointing him as Eurandom Professor in the fall
 of 2002. He also thanks the Mittag-Leffler Inst. and Eurandom
for providing him with excellent facilities and
 for their hospitality during his visits.

V.Sidoravicius thanks Cornell University and the
Mittag-Leffler Institute for their hospitality and travel support.
His research was supported by FAPERJ Grant E-26/151.905/2001,
CNPq (Pronex).

\subhead
2. Results from [KSb]
\endsubhead
\numsec=2
\numfor=1

{\it Throughout the rest of this paper we assume that
$$
D_A=D_B
\teq(2.0)
$$
and we
abbreviate their common value to} $D$.
We begin this section with some further facts about the setup. We
concentrate on the special case $D_A=D_B$.
More details can be found in Section 2
of \cite {KSb}  which deals with the construction of our particle system.
$\{S_t\}_{t \ge 0}$ will be a continuous time simple random walk on
$\Bbb Z^d$ with jumprate $D$ and starting at $\bold 0$.
To each initial particle $\rho$ is
assigned a path $\{\pi_A(t,\rho)\}_{t \ge
0}$ which is distributed like $\{S_t\}_{t \ge 0}$.
The paths $\pi_A(\cdot, \rho)$ for different $\rho$'s are independent
and they are all independent of the initial $N_A(x,0-), x \in \Bbb
Z^d$. The position of $\rho$ at time $t$ equals $\pi(0, \rho) +
\pi_A(t,\rho)$, and this can be assigned to $\rho$ without knowing the
paths of any of the other particles. The
type of $\rho$ at time $s$ is denoted by $\eta(s,\rho)$. This equals
$A$ for $0 \le s < \th(\rho)$ and equals $B$ for $s \ge
\th(\rho)$, where $\th(\rho)$, the so-called switching time of $\rho$,
is the first time at which
$\rho$ coincides with an initial $B$-particle. Note that this is
simpler than in the construction of \cite {KSb} for the general case
which may have $D_A \ne D_B$. In that case we had simple random walks
$\{S^\eta\}_{t \ge 0}$ with jumprate $D_\eta$ for $\eta \in \{A,B\}$,
and there were two paths associated with each initial particle $\rho:
\pi_\eta(\cdot, \rho), \eta \in \{A,B\}$, with $\{\pi_\eta(t,\rho)\}$
having the same distribution as $\{S^\eta_t\}$. If $\rho$ had initial
position $z$, its position was then equal to $z+ \pi_A(0,\rho)$ until
$\rho$ first coincided with a $B$-particle at time $\th(\rho)$; for $t
\ge \th(\rho)$ the position of $\rho$ was $z+\pi_A(\th(\rho),\rho) +
[\pi_B(t,\rho)- \pi_B(\th(\rho), \rho)]$. This depends on $\th(\rho)$
and therefore on the movement of all the other particles. In the
present case we can take $\pi_B = \pi_A$, which has the great
advantage that the path of $\rho$ does not depend on the paths of the
other particles. This is the reason why the case $D_A=D_B$ is special.
We proved in \cite {KSb} that on a certain state space $\Si_0$, the
collection of positions and types of all particles at time $t$,
with $t$ running from 0 to $\infty$,  is well
defined and forms a strong Markov process with respect to the
$\si$-fields $\cF_t = \cap_{h > 0} \cF_{t+h}^0,\; t \ge 0$,
where $\cF_t^0$ is the $\si$-field generated by the positions and types
of all particles during $[0,t]$. The elements of these $\si$-fields
are subsets of $\Si^{[0, \infty)}$, where $\Si = \prod_{k \ge
1}\big((\Bbb Z^d \cup \partial_k) \times \{A,B\}\big)$. $\Si^{[0,\infty)}$
is the pathspace
for the positions and types of all particles.
More explicit definitions are given in \cite
{KSb} but are probably not needed for this paper.
It was also shown in \cite {KSb}
that if one chooses the number of initial $A$-particles at $z$, with
$z$ varying over $\Bbb Z^d$, as
i.i.d. mean $\mu_A$ Poisson variables, then the process starts off
in $\Si_0$ and stays in $\Si_0$ forever, almost surely.

We write $N_\eta(z,t)$ for the number of particles of type
$\eta$ at the space-time point $(z,t), \, z \in \Bbb Z^d,
\eta \in \{A,B\}$, while
$N_A(z,0-)$ denotes the number of $A$-particles
to be put at $z$ ``just before'' the system starts evolving.
 Note that our model always has only particles of
one type at each given site, because an $A$-particle which meets a
$B$-particle changes instantaneously to a $B$-particle. Thus,
if $N_A(z,0-) =N$ for some site $z$ and we
add $M > 0$ $B$-particles at $z$ at time 0, then we have to say
that $N_A(z,0) = 0, N_B(z,0) = N+M$.

We shall rely heavily on basic upper and lower bounds for the growth
of $B(t)$ which come from Theorems 1 and 2 in \cite {KSb}.
\proclaim{Theorem A} If $D_A = D_B$, then
there exist constants $0 < C_2\le C_1 <
\infty$ such that for every fixed K
$$
P\big\{\cC(C_2t) \subset B(t) \subset \cC(2C_1t)\big\} \ge 1 -
\frac 1{t^K}
\teq(2.1)
$$
for all sufficiently large $t$.
\endproclaim
We also have some information about the presence of $A$-particles in
the regions which have already been visited by $B$-particles. The
following is Proposition 3 of \cite {KSb}.
\proclaim{Proposition B} If $D_A = D_B$, then
for all $K$ there exists a constant $C_3 = C_3(K)$ such that
$$
\align
&P\{\text{there is a vertex $z$ and an  $A$-particle at the space-time
point $(z,t)$ while}\\
&\phantom{M}\text{there also was a $B$-particle at $z$ at some time }
\le t - C_3[t\log t]^{1/2}\} \\
&\le \frac 1{t^K}\text{ for all sufficiently large }t.
\teq(2.2)
\endalign
$$
Consequently, for large $t$
$$
\align
&P\{\text{at time $t$ there is a site in $\cC\big(C_2t/2\big)$ which}\\
&\phantom{MMMMMMMM}\text{is occupied
by an $A$-particle}\}
\le \frac 2{t^K}.
\teq(2.3)
\endalign
$$
\endproclaim

Finally we reproduce here Lemma 15 of \cite {KSb} which gives
an important monotonicity property. We repeat that in the present
setup, with the $N_A(x,0-)$ i.i.d. Poisson variables, our process
a.s. has  values in $\Si_0$ at all times (see Proposition 5 of \cite
{KSb}).
\proclaim{Lemma C}
Assume $D_A = D_B$ and let $\si^{(2)} \in
\Si_0$. Assume further that $\si^{(1)}$ lies below $\si^{(2)}$ in the following
sense:
$$
\text{for any site $z \in \Bbb Z^d$, all particles present in
$\si^{(1)}$ at $z$ are also present in $\si^{(2)}$ at $z$},
\teq(7.1a)
$$
and
$$
\align
&\text{at any site $z$ at which the particles in $\si^{(2)}$ have type $A$},\\
&\text{the particles also have type $A$ in $\si^{(1)}$}.
\teq(7.2a)
\endalign
$$
Let $\pi_A(\cdot,\rho)$
be the random walk paths associated to the
 various particles and assume that the Markov processes
 $\{Y_t^{(1)}\}$ and $\{Y_t^{(2)}\}$ are constructed by means of
 the same set of paths $\pi_A(\cdot,\rho)$
and starting with state $\si^{(1)}$ and $\si^{(2)}$, respectively
 {\rom (}as defined in Section 2 of \cite {KSb}, but with
$\pi_A(\cdot,\rho)=\pi_B(\cdot,\rho)$ for all $s,\rho$; see (2.6),
(2.7) there{\rom )}. Then, almost surely,
$\{Y_t^{(1)}\}$ and
 $\{Y^{(2)}_t\}$ satisfy
 \equ(7.1a) and \equ(7.2a) for all $t$
with $\si^{(i)}$ replaced by $Y^{(i)}_t,\;
 i=1,2$. In particular, $\si^{(1)} \in \Si_0$.
\endproclaim
In particular, this monotonicity property says that if $\si^{(1)}$
is obtained from $\si^{(2)}$ by removal of some particles and/or changing
some $B$-particles to $A$-particles, then the process starting from
$\si^{(1)}$ has no more $B$-particles at each space-time point
than the process starting from $\si^{(2)}$. We note that this
monotonicity property holds only under our basic assumption that $D_A=D_B$.

\subhead
3. A subadditivity relation
\endsubhead
\numsec=3
\numfor=1

In this section we shall prove a basic subadditivity relation and
deduce from it that the $B$-particles spread in each fixed direction
over a distance which grows asymptotically linearly with time. This
statement is ambiguous because we haven't made precise what ``spread
in a fixed direction'' means. Here this will be
measured by
$$
\max\{\langle x,u \rangle: x \in \wt B(t)\},
\teq(3.1)
$$
where $u$ is a given unit vector (in the Euclidean norm) in $\Bbb
R^d$ (see the abstract for $\wt B$).
In addition we
will not prove subadditivity (which is an almost sure relation),
but only superconvolutivity, in the terminology of
\cite {Ha} (which is a relation between distribution functions).
The tool of superconvolutivity in other models
with no obvious subadditivity in the strict sense goes back to
\cite {R}, and was also used in \cite {BG} and \cite{W}.

We define the closed half-space
$$
\cS(u,c)= \{x \in \Bbb R^d: \langle x, u \rangle  \ge c\}.
$$
Given a $u \in S^{d-1}$ and $r \ge
0$ we
consider the {\it half-space process corresponding to }$(u,-r)$ (also called
$(u,-r)$ half-space-process).
We define this to be the
process whose initial state is of the form
$$
\align
N_A(x, 0-) = 0 \text{ if } x \notin \cS(u, -r)& \text{ and the }
N_A(x,0-),\; x \in \cS(u,-r),\\
&\text{ are i.i.d., mean $\mu_A$ Poisson variables.}
\endalign
$$
In addition the particles at $x_{0,-r}$ are turned into $B$-particles at
time 0, where $x_{0,-r}$ is the site in $\cS(u,-r)$ nearest to the origin (in
$\l^\infty$-norm) with $N_A(x_{0,-r},0-) > 0$;
 if there are several possible
choices for $x_{0,-r}$, the tie is broken according to some deterministic
rule chosen in advance.
There will be many other occasions were ties
may occur. These will be broken in the same way as here, but we shall
not mention ties or the breaking of them anymore.
Note that no extra $B$-particles are
introduced at time 0, but that only the type of the particles at $x_{0,-r}$
is changed. Thus,
$$
N_A(x,0) + N_B(x,0) = N_A(x,0-) \text{ for
all $x$}.
\teq(3.0abc)
$$
From time 0 on the particles move and change
type as described in the abstract.
Note that only the initial state is restricted
to $\cS(u,-r)$. Once the particles start to move they are free to leave
$\cS(u,-r)$. The $(u,-r)$ half-space process will often be denoted by
$\cP^h(u,-r)$.

We further define the $(u,-r)$ {\it half-space process
starting at }$(x,t)$.
 This process is defined for times $t' \ge t$ only. We define it
as follows: at time $t$ let $x_{0,-r}(t)$ be the nearest site to $x$ which is
occupied in the $(u,-r)$ half-space process.
We then reset the types of the particles at $x_{0,-r}(t)$ to $B$ and the
types of all other particles present in the $(u,-r)$ half-space process
at time $t$ to $A$.  The
particles then move along the same path in the $(u,-r)$ half-space
process starting at $(x,t)$ as in $\cP^h(u,-r)$
(which starts at $(\bold 0,0)$). However, the types of the
particles in the $(u,-r)$ half-space
process starting at $(x,t)$ are determined on the basis of the reset
types at time $t$.  Thus the half-space process starting
at $(x,t)$ has at any time only particles which were in $\cS(u,-r)$ at
time 0. Moreover, at any site $y$ and time
$t'\ge t$, $\cP^h(u,-r)$ and the $(u,-r)$ half-space process started at $(x,t)$
contain exactly the same particles.
We see from this that the {\it paths} of the particles in the $(u,-r)$
half-space processes starting at $(x,t)$ and at $(\bold 0,0)$
are coupled so that they coincide from time $t$ on, but the types of a
particle in these two processes may differ. Lemma C
shows that if there is a $B$-particle in $\cP^h(u,-r)$
at $x$ at time $t$, then
in this coupling any $B$-particle in the
$(u,-r)$ half-space process starting at $(x,t)$ also has to have type
$B$ in $\cP^h(u,-r)$.

The coupling between the two half-space processes clearly relies heavily on
the assumption $D_A =D_B$, so that we can assign the same path to a
particle in the two processes, even though the types of the particle
in the two processes may be different.

It is somewhat unnatural to start the $(u,-r)$ half-space process with
$B$-particles at $x_{0,-r}$ in case $r <0$, so that the origin does not lie
in the half-space $\cS(u,-r)$. We shall avoid that situation. We can,
however, use the $(u,-r)$ half-space process starting at $(x,t)$. This
is well defined for all $r$. We merely need to find
the site nearest to $x$ which has at time $t$ a particle which started
in $\cS(u,-r)$ at time 0. We can then reset the type of the particles
at this site to $B$ at time $t$.
We shall consider the $(u,-r)$ half-space process starting at $(x,t)$
mostly in cases where we already know that $x$ itself is occupied at time $t$
in the $(u,-r)$ half-space process.

Finally we shall occasionally talk about the {\it full-space process}
and the {\it full-space process starting at} $(x,t)$. These are
defined just as the half-space processes, but with $r = \infty$. In
particular, the full-space process starts with $B$-particles
only at the nearest occupied site to the origin and \equ(3.0abc)
applies. The full-space
process starting at $(x,t)$ has $B$-particles at time $t$ only at
the nearest occupied site to $x$. The type of all particles at other
sites are reset to $A$ at time $t$.
By stationarity in time, the full-space process started at $(x,t)$ has
the same distribution at the space-time point $(x+y,t+s)$
as the full-space process (started at
$(\bold 0,0)$) at the point $(y,s)$. Again we shall use the same
random walk paths $\pi_A$ for all the full
state processes and the half-space processes, so that these processes
are automatically coupled. We shall denote the full-space process by $\cP^f$.

We point out that if $0 \le r_1 \le r_2$, and if $\|x_{0,-r}\|
\le r_1/\sqrt d$,
then $x_{0,-r} \in \cS(u,-r_1) \subset \cS(u,-r_2)$. In this case, both
$\cP^h(u,-r_1)$ and $\cP^h(u,-r_2)$ start with changing the type to
$B$ at the site $x_{0,-r}$ only. By Lemma C, at any time
$$
\text{any $B$-particle in $\cP^h(u,-r_1)$ is also a $B$-particle in
$\cP^h(u,-r_2)$}.
\teq(3.0cde)
$$
 This comment also applies if $\cP^h(u,-r_2)$ is
replaced by $\cP^f$ (which is the case $r_2 = \infty$).

It seems worthwhile to discuss more explicitly the relation of the
full-space process to our process as described in  the
abstract. The latter has some $B$-particles introduced at time 0 at
one or more sites, in
addition to the Poisson numbers of particles, $N_A(x,0-), x \in \Bbb
Z^d$. If exactly one $B$-particle is added at time 0, and this
particle is placed at $\bold 0$, then we shall call the resulting
process the {\it original process}.

 Suppose we want to estimate $P\{\cA(x_0)\}$ in the full-space
process,  where
$$
x_0 := \text{the nearest occupied site to the origin at
time 0 in }\cP^f,
\teq(3.0xyz)
$$
 $\cA$ is some event and $\cA(x)$ is the translation by $x$ of
this event  (which takes $N_A(\bold 0,s)$ to $N_A(x,s)$). Then, for
$C$ a subset of $\Bbb Z^d$,
$$
\align
&P\{x_0 \in C,\cA(x_0) \text{ in }\cP^f\} = \sum_{x
\in C} P\{x_0 = x, \cA(x)\}\\
&\le \sum_{x \in C} P\{\text{$x$ is occupied at time 0}, \cA(x)
\text{ in }\cP^f\}\\
&= \sum_{x \in C} \sum_{k=1}^\infty e^{-\mu_A} \frac {[\mu_A]^k}{k!}
P\{\cA|\text{there are $k$ $B$-particles at $\bold 0$ at time 0}\}.
\teq(3.0)
\endalign
$$
(The probability in the last sum is the same in $\cP^f$ as in the
original process.)
On the other hand, in the original process we have
$$
\align
&P\{\cA \text{ in original process}\} \\
&= \sum_{k=1}^\infty
e^{-\mu_A} \frac{[\mu_A]^{k-1}}{(k-1)!}P\{\cA|\text{there are $k$
$B$-particles at $\bold 0$ at time 0}\}.
\teq(3.0a)
\endalign
$$
Comparison of the right hand sides in \equ(3.0) and \equ(3.0a)
yields the crude bound
$$
\align
&P\{x_0 \in C,\cA(x_0) \text{ in the full-space process}\}\\
&\qquad\le \text{(cardinality of $C$)}\mu_A P\{\cA \text{ in original
process}\}.
\teq(3.0b)
\endalign
$$
We shall repeatedly use a somewhat more general version of this
inequality (see for instance (3.25), (3.77), (3.78), (5.33)).
Suppose $s \ge 0$ is fixed and $X$ is a
random vertex in $\Bbb Z^d$, and suppose further that
$$
P\{\text{in }\cP^f, \cA(X) \text{ but $(X,s)$ is not occupied}\}=0.
\teq(3.0cba)
$$
(Note that this is satisfied if $(X,s)$ is occupied almost surely in $\cP^f$.)
Let $C \subset \Bbb Z^d$ as before.
Now, given that there are $k \ge 1$ particles at the (non-random)
space-time point $(x,s)$,
the full-space process starting at $(x,s)$ is simply a translation by
$(x,s)$ in space-time of the original process, conditioned to start
with $k-1$ points at the origin and one $B$-particle added at the origin.
Therefore,
essentially for the same reasons as for \equ(3.0b),
$$
\align
&P\{X \in C, \cA(X) \text{ in the full-space process starting at $(X,s)$}\}\\
&\qquad\le \text{(cardinality of $C$)}\mu_A P\{\cA \text{ in original
process}\}.
\teq(3.0c)
\endalign
$$

For a rather trivial comparison in the other direction we note that if
$P\{\cA \text{ in } \cP^f\} =0$
for the full-space process, then we certainly have
for each $k \ge 1$ that
$$
\align
&0=P\{\cA \text{ in } \cP^f, x_0 = \bold 0, k \text{ particles at }x_0\}\\
&=P\{\cA \text{ in } \cP^f, k \text{ particles at } \bold 0\}\\
&= e^{-\mu_A} \frac{[\mu_A]^k}{k!}P\{\cA|\text{there are $k$
$B$-particles at $\bold 0$ at time 0}\}.
\teq(3.0d)
\endalign
$$
This implies, via \equ(3.0a) that also $P\{\cA \text{ in } \cP^f\} = 0$.

It is somewhat more complicated to compare $\cP^f$ with the process
described in the abstract if more than one $B$-particle is introduced
at time 0. Rather than develop general results in this direction we
merely show in our first lemma that it suffices to prove \equ(1.3)
for the full-space process.
\proclaim{Lemma 1} If \equ(1.3) holds in $\cP^f$, then it also holds in
the original process of the abstract with any fixed finite number of
$B$-particles added at time 0.
\endproclaim
\demo{Proof} The preceding discussion shows that if \equ(1.3) has
probability 1 in $\cP^f$, then it has probability 1 in the
original process (with one particle added at the origin at time 0). By
translation invariance \equ(1.3) will then have probability 1 in the
process of the abstract with one particle added at any fixed site at time 0.

Lemma C implies that one can couple two processes as in the abstract, with
collections of $B$-particles $A^{(1)} \subset A^{(2)}$
added at time 0, respectively, in such a way that the process
corresponding to $A^{(1)}$ always has no more $B$-particles than the
one corresponding to $A^{(2)}$. Therefore, if the left hand inclusion
in \equ(1.3) holds when
only one $B$-particle is added at time 0, then it certainly
holds if more than one $B$-particle are added.

It follows that we only have to prove the right hand inclusion in
\equ(1.3) for the process from the abstract with more than one particle added,
if we already know it when exactly one particle is added.
Assume first that we run this last process with one $B$-particle
$\rho_0$ added at $z_0$. We now
have to refer the reader to the genealogical paths
introduced in the proof of Proposition 5 of \cite {KSb}. The right
hand inclusion in \equ(1.3) then says that for all $\ep > 0$
$$
\align
&P\{\text{there exist genealogical paths from $z_0$ to some point}\\
&\phantom{MMMM}\text{outside
$(1+\ep)tB_0$ for arbitrarly large $t$}\} = 0.
\teq(3.0de)
\endalign
$$
From the construction of the genealogical paths in Proposition 5 of
\cite {KSb} and the fact
that a.s. there are only finitely many $B$-particles at finite times
(see (2.18) in \cite {KSb}) it is not hard to deduce that
$$
\align
&\{\wt B(t) \not \subset (1+\ep)tB_0 \text{ at time $t$
if one adds a $B$-particle
$\rho_i$}\\
&\phantom{MMMM}\text{at $z_i,\; 1 \le i \le k$, at time 0}\}\\
&=\{\text{there is a genealogical path from some $z_i,\; 1 \le i
\le k$,}\\
&\phantom{MMMM}\text{to the complement of $(1+\ep)tB_0$ at time
$t$ if one}\\
&\phantom{MMMM}\text{adds a $B$-particle
$\rho_i$ at $z_i,\; 1 \le i \le k$, at time 0}\}\\
&\subset \bigcup_{i=1}^k
\{\text{there is a genealogical path from $z_i$
to the complement of}\\
&\phantom{MMMM}\text{$(1+\ep)tB_0$ at time
$t$ if one adds a $B$-particle $\rho_i$ at $z_i$ at time 0}\}
\teq(3.0ef)
\endalign
$$
(the $z_i$ do not have to be distinct here).
It follows that
$$
\align
&P\{\wt B(t) \not \subset (1+\ep)tB_0 \text{ for arbitrarily large times $t$
if one}\\
&\phantom{MMMM}\text{adds a $B$-particle
$\rho_i$ at $z_i,\; 1 \le i \le k$, at time 0}\}\\
&\le \sum_{i=1}^k
P\{\text{there are genealogical paths from $z_i$
to the complement of $(1+\ep)tB_0$}\\
&\phantom{MMMM}\text{at arbitrarily large times
$t$ if one adds a $B$-particle $\rho_i$ at $z_i$ at time 0}\}\\
&= 0 \text{ (by \equ(3.0de))}.
\endalign
$$
Thus the right hand inclusion in \equ(1.3) holds a.s., even if one
adds $k$ $B$-particles at time 0.
\hfill $\blacksquare$
\enddemo

We recall that
$$
\cP^h(u,-r) \text{ is short for the $(u,-r)$ half-space process},
$$
$$
\cP^f \text{ is short for the full-space process},
$$
and we further introduce
$$
\cB^h(y,s;u,-r):=\{\text{there is a $B$-particle at $(y,s)$ in }
\cP^h(u,-r)\},
\teq(3.0f)
$$
$$
h(t,u,-r) = \max \{\langle x,u \rangle:\cB^h (x,t;u,-r) \text{ occurs}\}.
\teq(3.1a)
$$
$P^{or}$ will denote the probability measure for the original process
(with one $B$-particle added at the origin at time 0); $E^{or}$ is
expectation with respect to $P^{or}$.
(The superscripts $h, f$ and $or$ are added to various symbols
which refer to a
half-space process, the full-space process, or the original process,
respectively). We use $P$ without superscript if it is clear from the
context with which process we are dealing or when we are discussing the
probability of an event which is described entirely in terms of
the $N_A(x,0-)$ and the paths $\pi_A$.

The following technical lemma will be useful.
It tells us that, with high probability, $\cP^h(u,-r)$ moves out in the
direction of $u$ at least at the speed $C_4$, provided $r$ is large
enough (see (3.15)).
Its proof would be nicer
if we made use of the fact that even the $(u,0)$ halfspace-process has,
with a probability at least $1-t^{-K},\; B$-particles at time $t$
at sites $x$ with $\langle x,u \rangle \ge
Ct$, for some constant $C >0$. However, it takes some work to prove
this fact and we decided to do without it.
The lemma itself is proven by recursively
constructing a sequence
of space-time points which move out in the direction of $u$ along an
exponentialy growing sequence, so that there is
only an exponentially  small (in $k$) probability that the $k$-th
point is not occupied in the $(u,-r)$ half-space process.
\proclaim{Lemma 2}
Let $C_1, C_2$ be as in Theorem A and let
$$
C_4 = \frac{2\sqrt d C_1C_2}{32 \sqrt d C_1 + C_2}.
\teq(3.9ab)
$$
For all constants $K \ge 0$, there exists a constant $r_0 = r_0(K) \ge 0$
such that for $r \ge r_0$
$$
P\Big\{h(t,u,-r) \le C_4t \text{ for some }t \ge t_1 := \frac 1{4\sqrt d
C_1}\Big
[1+ \frac{C_2}{32\sqrt d C_1}\Big]r\Big\}
\le r^{-K}.
\teq(3.10)
$$
\endproclaim
\demo{Proof} \newline
{\bf Step 1.} For $k
\ge 1$ define the times
$$
t_k = \frac 1{4 \sqrt d C_1}\Big[1+\frac {C_2}{32 \sqrt d C_1}\Big]^k r,
$$
and the real numbers
$$
d_k = \frac {C_2}{32 \sqrt d C_1}\Big[1+\frac {C_2}{32 \sqrt d C_1}\Big]^k r.
$$
Also define for each $k \ge 1$ the event
$$
\align
\cD_k :=\big\{&\cB^h(x_k,t_k;u,-r) \text{ occurs for some $x_k$
which}\\
&\text{satisfies
$\langle x_k,u \rangle \ge d_k$ and } \|x_k\| \le 2C_1 t_k\big\}.
\teq(3.12)
\endalign
$$
In this step we shall reduce the lemma to an estimate for the
probability that $\cD_k$
fails for some $k \ge 1$. Indeed, assume that $\cD_k$ occurs for all
$k \ge 1$. By definition,
there is then a $B$-particle at
$(x_k,t_k)$ in the $(u,-r)$ half-space process (starting at $(\bold
0,0)$), so that
$$
h(t_k,u,-r) \ge \langle x_k,u \rangle \ge d_k
= \frac{C_2}{32 \sqrt d C_1}
\Big[1+\frac {C_2}{32 \sqrt d C_1}\Big]^k r , \quad k \ge 1.
\teq(3.12b)
$$
Recall that $\cF_t$ is defined in the beginning of Section 2.
In addition to \equ(3.12b), we have
on the event $\{\langle x_k,u \rangle \ge d_k\}$, for $k \ge 1$,
$$
\align
&P\{h(t,u,-r) \le \frac 12 d_k \text{ for some }t \in [t_k, t_{k+1})|
\cF_{t_k}\}\\
&\le P\{\text{each $B$-particle in $\cP^h(u,-r)$ at $(x_k,t_k)$
moves during }\\
&\phantom {MMMMMMMMMM}\text{$[t_k,t_{k+1}]$ to some site $x$
with $\langle x,u \rangle \le \frac 12 d_k$}\}\\
&\le P\{\min_{q \le t_{k+1}-t_k} \langle S_q,u \rangle
\le  -\frac12 d_k = -C_4 t_{k+1} \}\\
& \le K_1\exp[-K_2t_{k+1}]
\teq(3.12a)
\endalign
$$
for some constants $K_1, K_2$ depending on $d, D_A$ only;
see (2.42) in \cite {KSa} for the last inequality.
It follows that the left hand side of \equ(3.10) is bounded by
$$
P\{\text{$\cD_k$ fails for some }k \ge 1\}
+ \sum_{k=1}^\infty K_1\exp[-K_2t_k].
\teq(3.14)
$$
\newline
{\bf Step 2.} We shall now derive a recursive bound for
$\cap_{1 \le j \le k}\cD_j$.
Assume that  $\cap_{1 \le j \le k-1}\cD_j$
occurs for some $k \ge 2$.
Consider now the full-space process starting at
$(x_{k-1},t_{k-1})$. Define
the following events for this process:
$$
\align
\cE_{k,1} &:= \{\text{at time }t_k \text{ all occupied
sites in }\\
&\phantom{MM}x_{k-1} + \cC \big((C_2/2)(t_k-t_{k-1})\big)
\text{ contain in fact a $B$-particle}\},\\
\cE_{k,2}& :=\{\text{at time $t_k$ there is an occupied site in }\\
&\phantom{MM}x_{k-1} + (C_2/4)(t_k-t_{k-1})u + \cC\big([\log t_k]^2\big)\},\\
\cE_{k,3} &:= \{\text{all particles in $x_{k-1}+
\cC\big(2C_1(t_k-t_{k-1})\big)$}\\
&\phantom{MM}\text{at time $t_{k-1}$ started at time 0 in }\cS(u,-r)\},\\
\cE_{k,4} &:= \{\text{there is no $B$-particle outside $x_{k-1} +
\cC\big(C_1(t_k -t_{k-1})\big)$}\\
&\phantom{MM}\text{during }[t_{k-1}, t_k]\},\\
\cE_{k,5} &:= \{\text{no particle which is outside $x_{k-1} +
\cC\big(2C_1(t_k - t_{k-1})\big)$}\\
&\phantom{MM} \text{at time $t_{k-1}$ enters
$x_{k-1} +\cC\big(C_1(t_k - t_{k-1})\big)$ during }[t_{k-1},t_k]\},\\
\endalign
$$

We claim that on
$$
\cD_{k-1} \cap \bigcap_{1 \le i \le 5}\cE_{k,i}
\teq(3.19aa)
$$
also $\cD_k$ occurs, provided $r \ge $ some suitable $r_1$,
independent of $k$, and $k \ge 2$. We merely
need to make sure that $\sqrt d [\log t_k]^2 \le (C_2/8)(t_k - t_{k-1})$
whenever $r \ge r_1$.
To prove our claim when $k \ge 2$,
observe first that the occurrence of
$\cE_{k,1} \cap \cE_{k,2}$ guarantees that at time $t_k$
there is a $B$-particle at some $x_k$ in
$x_{k-1} + (C_2/4)(t_k-t_{k-1})u + \cC\big([\log t_k]^2)
\subset x_{k-1} + \cC\big((C_2/2)(t_k-t_{k-1})\big)$. Such a particle
automatically satisfies
$$
\langle x_k,u \rangle \ge  \langle x_{k-1},u \rangle +
\frac {C_2}4 (t_k-t_{k-1}) - \sqrt d [\log t_k]^2 \ge d_{k-1} +
\frac {C_2}8 (t_k-t_{k-1}) =  d_k.
\teq(3.20)
$$
It also satisfies $\|x_k\| \le 2C_1t_k$, because
$\|x_{k-1}\| \le 2C_1t_{k-1}$, and on $\cE_{k,2},\; \|x_k\| \le
\|x_{k-1}\| +(C_2/4)(t_k-t_{k-1}) +[\log t_k]^2$, while $C_2 \le C_1$.
This particle at $(x_k,t_k)$ is a
$B$-particle in the full-space process starting at
$(x_{k-1},t_{k-1})$. We are going to show that, in fact, it is also a
$B$-particle in the $(u,-r)$ half-space process starting at
$(x_{k-1},t_{k-1})$. This will prove our claim,
because the monotonicity property of Lemma C
implies that any $B$-particle
in the $(u,-r)$ half-space process starting at $(x_{k-1},t_{k-1})$ is also a
$B$-particle in the $(u,-r)$ half-space process (starting at $(\bold 0,0)$),
provided that there is a $B$-particle at $(x_{k-1},t_{k-1})$ in the $(u,-r)$
half-space process. (Note that this proviso is satisfied
 by the induction assumption that $\cD_{k-1}$ occured.)
We first observe that the particle at $(x_k,t_k)$ must at
time $t_{k-1}$ have been in $x_{k-1} +\cC\big(2C_1(t_k-t_{k-1})\big)$,
because $x_k \in x_{k-1} + \cC\big((C_2/2)(t_k-t_{k-1})\big) \subset
x_{k-1} + \cC\big((C_1/2)(t_k-t_{k-1})\big)$ and
$\cE_{k,5}$ occurs.
By virtue of $\cE_{k,3}$ this particle then belongs to $\cP^h(u,-r)$
as well as to the $(u,-r)$ half-space
process
starting at $(x_{k-1},t_{k-1})$. We still have to show that this
particle  also has
type $B$ in the $(u,-r)$ half-space process starting at $(x_{k-1},t_{k-1})$.
To this end we note that the particles starting
outside $x_{k-1} +
\cC\big(2C_1(t_k - t_{k-1})\big)$ at time $t_{k-1}$ do not influence
the type of any particle at time $t_k$
in the full-space
process starting at $(x_{k-1}, t_{k-1})$. Indeed, in this process
the particles outside $x_{k-1} +
\cC\big(2C_1(t_k - t_{k-1})\big)$ start as $A$-particles,
and since $\cE_{k,4}\cap\cE_{k,5}$
occurs, these particles  do not meet any $B$-particle at or before time
$t_k$.
Thus, whether the particle at $(x_k,t_k)$ is also a $B$-particle in the
$(u,-r)$ half-space process starting at $(x_{k-1}, t_{k-1})$ depends
only on the paths of the particles which were in
$x_{k-1} +
\cC\big(2C_1(t_k - t_{k-1})\big)$ at time $t_{k-1}$
(compare the lines following (2.37) in \cite{KSb}). All these
particles were particles in $\cP^h(u,-r)$ at time
$t_{k-1}$ (on $\cE_{k,3}$), and hence also are in this half-space
process at time $t_k$. Thus the type of the particle at $(x_k,t_k)$
is the same in the full-space process starting at $(x_{k-1}, t_{k-1})$
and in the $(u,-r)$ half-space process  starting at $(x_{k-1}, t_{k-1})$.
This justifies our claim that $\cD_k$ occurs for $k \ge 2$. We leave
it to the reader to make some simple modifications in the above
argument to show that $\cD_1$ occurs on
$$
\cD_0 \cap \bigcap_{1 \le i \le 5}\cE_{1,i},
$$
where
$$
t_0 = 0 \text{ and } \cD_0 = \{\|x_0\| \le K_3\log
r\},
$$
provided $r_1$ is chosen large enough; $x_0$ is defined in
\equ(3.0xyz) and $K_3$ is chosen right after
(3.25) and depends on $K, d$ and $\mu_A$ only.

We have now shown that on the event \equ(3.19aa) also $\cD_k$
occurs. If this is the case and also $\cap_{1\le i\le 5}\cE_{k+1,i}$
occurs, then $\cD_{k+1}$ occurs etc. Consequently, for $r \ge r_1$,
$$
\align
&P\{\text{$\cD_0$ occurs, but some $\cD_k$ fails}\}\\
&\le\sum_{i=1}^5 P\{\text{for some $x_0$ with } \|x_0\| \le K_3\log r,
\text{ $x_0$ is occupied but $ \cE_{1,i}$ fails}\}\\
&+
\sum_{k=2}^\infty \sum_{i=1}^5 P\{\text{for some
$x_{k-1}$ with }\|x_{k-1}\|\le 2C_1t_{k-1}\text{ and }
\langle x_{k-1},u\rangle \ge d_{k-1}, \\
&\phantom{MMMMMMMM}\text{$\cB^h(x_{k-1}, t_{k-1}; u,-r)$ occurs,
but $\cE_{k,i}$ fails}\}.
\teq(3.20a)
\endalign
$$
{\bf Step 3.} In this step we shall give most of the estimates
for the terms in the
right hand side here for $k \ge 2$. The basic inequalities remain
valid for $k=1$ by trivial modifications which we again leave to the reader.
For the various estimates we have to take
all $t_k$ large. This will automatically be the case if $r$ is large;
we shall not explicitly mention this in the estimates below.

We start with the estimate for the failure of $\cE_{k,1}$.
\comment
 note that if
$x$ is occupied by a $B$-particle at time $t$ in $\cP^f$, then,
by the monotonicity property of Lemma C, any point in
space-time which is occupied by an $A$-particle in the full-space
process is also occupied by an
$A$-particle in the full-space process started at $(x,t)$.
Therefore, if
\endcomment
If $\cE_{k,1}$
fails, for a given $(x_{k-1}, t_{k-1})$, then there must be some
$y \in x_{k-1} + \cC\big((C_2/2)(t_k-t_{k-1})$ such that $y$ is
occupied by an $A$-particle at time $t_k$ in the full-space process
started at $(x_{k-1}, t_{k-1})$. Recall that if we shift
the full-space process starting at $(x,t)$ by
$(-x,-t)$ in space-time, then we obtain the full state process
starting at
$(\bold 0,0)$.
Moreover, if we condition on the event that $x$ is
occupied at time $t$, then, after the shift by $(-x,-t)$ the
$N_A(y,0),\; y \ne \bold 0$, are i.i.d. mean $\mu_A$ Poisson random variables.
Therefore, by summing over the possible values for $x_{k-1}$,
$$
\align
&P\{\text{for some
$x_{k-1}$ with }\|x_{k-1}\|\le 2C_1t_{k-1}\text{ and }
\langle x_{k-1},u\rangle \ge d_{k-1}, \\
&\qqquad \text{ $\cB^h(x_{k-1},t_{k-1};u,-r)$ occurs,
but $\cE_{k,1}$ fails}\}\\
&\le \sum_{\|x \| \le 2C_1t_{k-1}} P\{\text{$\cB^h(x,t_{k-1};u,-r)$ occurs
and in the full-space process started}\\
&\qqquad \text{ at $(x,
t_{k-1})$
there is an
$A$-particle in $x+\cC\big((C_2/2)(t_k-t_{k-1})\big)$ at time }t_k\}\\
&\le \sum_{\|x \| \le 2C_1t_{k-1}} P\{\text{$\bold 0$ is occupied at
time 0 and in $\cP^f$ there is an
$A$-particle}\\
&\qqquad \text{ in $\cC\big((C_2/2)(t_k-t_{k-1})\big)$ at time }t_k
-t_{k-1}\}.
\endalign
$$
To the right hand side here we can apply
\equ(3.0b) (with $C = \{\bold 0\}$). This shows that the right hand
side is at most
$$
K_4[t_{k-1}]^d \mu_A P^{or}\{\text{at time }t_k-t_{k-1},
\text{ there is an $A$-particle in }\cC \big((C_2/2)(t_k-t_{k-1})\big)\}.
$$
The probability in the right hand side here is calculated for the
original process with one particle added at $\bold 0$ at time 0. By
\equ(2.3)  (with K replaced by $K+d+2$) this probability is at most
$2[t_k-t_{k-1}]^{-K-d-2}$. Therefore,
$$
\align
&P\{\text{for some
$x_{k-1}$ with }\|x_{k-1}\|\le 2C_1t_{k-1}\text{ and }
\langle x_{k-1},u\rangle \ge d_{k-1}, \\
&\qqquad \text{ $\cB^h(x_{k-1},t_{k-1};u,-r)$ occurs,
but $\cE_{k,1}$ fails}\}\\
&\le 2K_4[t_{k-1}]^d \mu_A [t_k-t_{k-1}]^{-K-d-2}
\le K_5 t_k^{-K-2}.
\endalign
$$

It turns out that in the estimates for $\cE_{k,2},\cE_{k,3}$ and $\cE_{k,5}$
we can ignore the type of the particle at $(x_{k-1},t_{k-1})$; we just
need that this space-time point is occupied.
For $\cE_{k,2}$ we again shift by $(-x_{k-1},-t_k)$,
sum over the possible values of $x_{k-1}$ and apply \equ(3.0b). This gives
$$
\align
&P\{\text{for some } x_{k-1} \text{ with }\|x_{k-1}\|\le 2C_1t_{k-1},
\text{ $(x_{k-1},t_{k-1})$ is occupied but $\cE_{k,2}$ fails}\}\\
&\le K_4[t_{k-1}]^d \mu_A P\{N_A(y,0-)= 0 \text{ for all }y \in
(C_2/4)(t_k-t_{k-1}) u +\cC([\log t_k]^2)\}\\
&\le t_k^{-K-2},
\endalign
$$
for large $r$, because the $N_A(y,0-)$ are independent.

Next, for
$\cE_{k,3}$ we use that on $\cD_{k-1}$, the distance between
$x_{k-1} + \cC\big(2C_1(t_k-t_{k-1})\big)$
and the complement of $\cS(u,-r)$ is at least
$$
\align
\langle x_{k-1},u \rangle +r - 2\sqrt d C_1 (t_k-t_{k-1})
&\ge d_{k-1}+r - 2\sqrt d C_1 (t_k-t_{k-1})\\
& = \frac 12 d_{k-1}+ r.
\endalign
$$
Thus, if we take the restriction $\langle x_{k-1},u
\rangle \ge d_{k-1}$ into account we find that
$$
\align
&P\{\text{for some
$x_{k-1}$ with }\|x_{k-1}\|\le 2C_1t_{k-1}\text{ and }
\langle x_{k-1},u\rangle \ge d_{k-1}, \\
&\phantom{MMMMM}\text{$x_{k-1}$ is occupied at time $t_{k-1}$
but $\cE_{k,3}$ fails}\}\\
&\le \sum \Sb \|x\| \le 2C_1t_{k-1}\\
\langle x,u \rangle \ge d_{k-1}
\endSb
\;\sum_{y \notin \cS(u,-r)} \;\sum_{z \in x+
\cC\big(2C_1(t_k-t_{k-1})\big)} P\{S_{t_{k-1}} = z-y\}\\
&\le K_6 t^d_{k-1} [t_k-t_{k-1}]^d
P\{\|S_{t_{k-1}}\| \ge \frac 12 d_{k-1} + r\}\\
&\le K_7t_k^{2d} \exp \Big[- K_8\frac {(d_{k-1}+r)^2}{t_{k-1} +
d_{k-1} +r}\Big] \text{ (by (2.42) in \cite {KSa})}\\
& \le t_k^{-K-2}.
\teq(3.18abc)
\endalign
$$

The estimate for $\cE_{k,4}$ comes from Theorem A, or rather Theorem 1
 in \cite {KSb}, which is the basis for the right hand inclusion in
 Theorem A. Indeed, we have, again by summing over the possible values
 of $x_{k-1}$ and using \equ(3.0c),
$$
\align
&P\{\text{for some }x_{k-1}\text{ with }\|x_{k-1}\| \le 2C_1t_{k-1},
\text{ $(x_{k-1},t_{k-1})$ is occupied but $\cE_{k,4}$ fails}\}\\
&\le K_4 [t_{k-1}]^d \mu_A P^{or}\{\text{there are $B$-particles outside }
\cC\big(C_1(t_k-t_{k-1})\big)\\
&\phantom{MMMMMMMMMMMMMMMMMMM}\text{ at some time } \le t_k-t_{k-1}\}.
\teq(3.18)
\endalign
$$
To estimate the probability on the right we argue as in the proof of
Theorem 3 of \cite{KSb}. If a particle has
type $B$ at some time $s \le t_k-t_{k-1}$ and is outside
the cube $\cC\big(C_1(t_k-t_{k-1})\big)$ at that time, then by
symmetry of the random walk $\{S\}$, the particle has a conditional
probability given $\cF_s$ at
least 1/2 of being outside $\cC\big(C_1(t_k-t_{k-1})\big)$ at time
$t_k-t_{k-1}$. Therefore (with $E^{or}$ denoting expectation with
respect to $P^{or}$),
$$
\align
&E^{or}\{\text{number of particles outside $\cC\big(C_1(t_k-t_{k-1})\big)$
at some time}\\
&\phantom{MMMMMMMMMMMMMMMMMMMM}\text{during } [0, t_k-t_{k-1}]\}\\
&\le 2E^{or}\{\text{number of particles outside $\cC\big(C_1(t_k-t_{k-1})\big)$
at time }t_k - t_{k-1}\}.
\endalign
$$
The expectation in the  right hand side here is exponentially small
in $(t_k - t_{k-1})$ by Theorem 1 of \cite {KSb} and is an upper bound for the
probability in the right hand of \equ(3.18). Thus
the left hand side of \equ(3.18) is at most $t_k^{-K-2}$ again.

The probability involving $\cE_{k,5}^c$ is also
$O(t_k^{-K-2})$. This can be
shown by large deviation estimates for random walks, analogously
to the terms involving $\cE_{k,3}^c$.

This provides the necessary bounds of the summands in
\equ(3.20a). Finally, we have
$$
\align
P\{\cD_0 \text{ fails} \} &\le P\{N_A(x,0-) = 0 \text{ for all
$x$ with }\|x\| <K_3 \log r\}\\
& = \exp\big[-\mu_A K_9[K_3 \log r]^d\big].
\teq(3.13)
\endalign
$$
Thus we can take $K_3= K_3(K,d, \mu_A)$
so large that the left hand side of \equ(3.13) is at most
$r^{-K-1}$ for all $r\ge 3$.
\equ(3.14), \equ(3.13) and \equ(3.20a) together
now show that the left hand side of \equ(3.10) is for large $r$ at
most
$$
K_{10} r^{-K-1}+ \sum_{k=1}^\infty K_1\exp[-K_2t_k]
+\sum_{k=1}^\infty K_{11}t_k^{-K-2} \le K_{12}r^{-K-1}.
\tag "$\blacksquare$"
$$
\enddemo

For any vector $v \in \Bbb R^d$, we define
$$
v^\perp = v^\perp(u):=v - \langle v,u \rangle u.
$$
We further introduce the following (semi-infinite) cylinders
with axis in the direction of $u$, for $\al, \be \in \Bbb R$ and $\ga
\in \Bbb R^d$ a vector orthogonal to $u$ (see Figure 1):
$$
\Ga(\al, \be,\ga) = \Ga(\al, \be, \ga, u) := \{x \in \Bbb R^d: \langle x,u
\rangle \ge \al,\| x^\perp- \ga\| \le \be\}
$$
and the events
$$
\align
\cG(\al, \be,\ga, \cP, t) &= \cG(\al, \be,\ga,  \cP, t,u)\\
:&=\{\text{in the process $\cP$ there is a $B$-particle
in $\Ga(\al,\be, \ga )$ at time $t$}\}.
\endalign
$$
The last definition will be used with $\cP$ taken equal to
some half-space, full-space or the original process.

\topinsert
\vspace{2truecm}
\epsfverbosetrue
\epsfxsize=300pt
\epsfysize=250pt
\centerline {\epsfbox{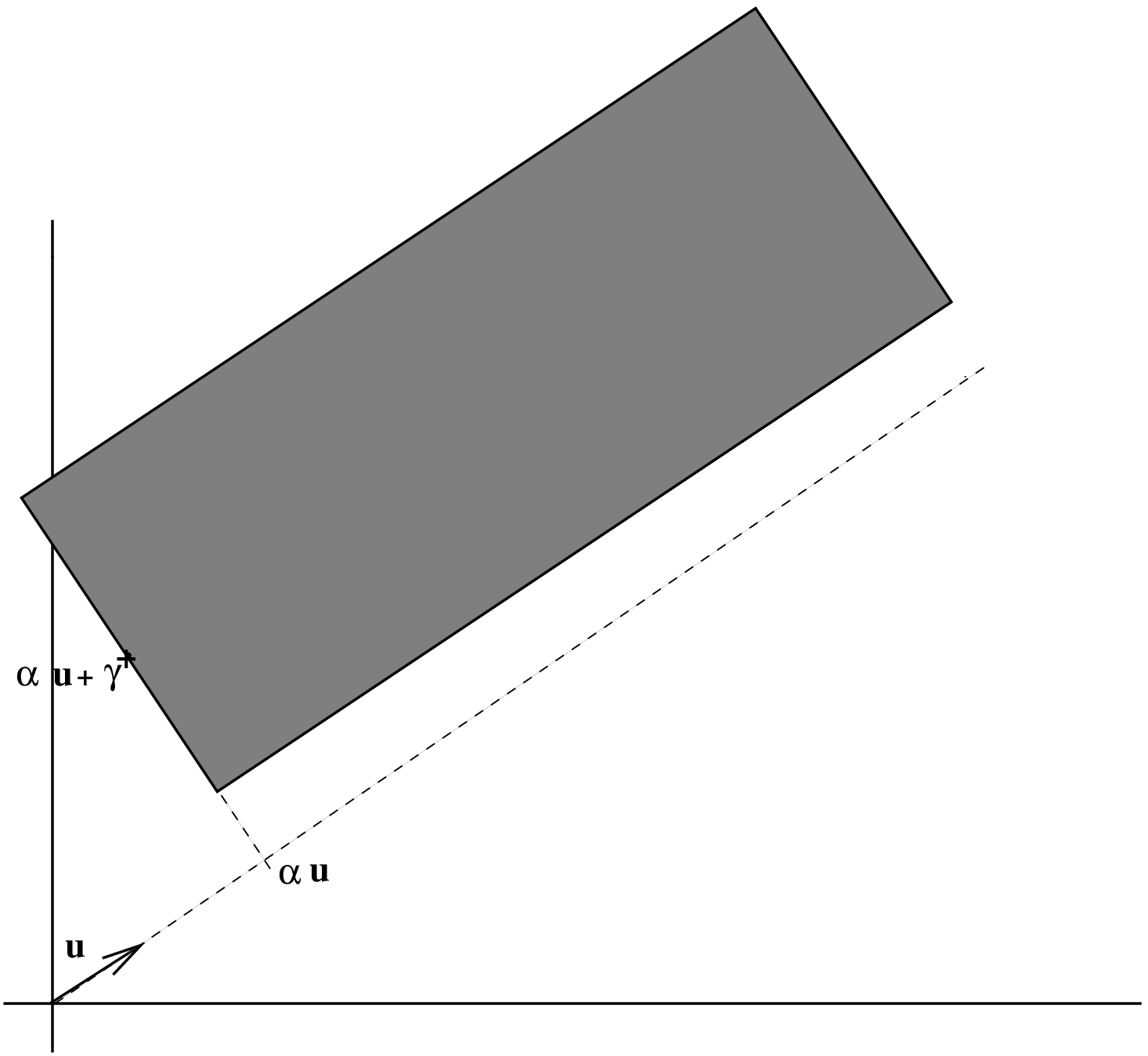}}

\botcaption
{Figure 1} The shaded region represents a cylinder $\Ga(\al,
\be,\ga)$; it extends to infinity on the upper right.
\endcaption
\endinsert
We remind the reader that $\cP^h(u,-r)$, the $(u,-r)$
half-space process, only uses particles which at time 0 are in the
half-space $\cS(u,-r) = \{x:\langle x,u \rangle \ge -r\}$.
We shall work a great deal with the process $\cP^h\big(u,
-C_5\ka(s)\big)$, where
$$
\ka(s) = [(s+1) \log (s+1)]^{1/2}
$$
and $C_5$ is a constant to be determined below (see the
line preceding (3.63)).
We make several more definitions:
$$
\align
h^*(s,u) :&= h(s,u, -C_5\ka(s))=\max \{\langle x,u \rangle:
\cB^h\big(x,s;u, -C_5\ka(s)\big) \text{ occurs}\}\\
&= \max \{\langle x,u \rangle:
x \text{ is occupied by a $B$-particle in }\\
&\phantom{MMMMMMMMMMMMMMMM}\cP^h\big(u, -C_5\ka(s)\big)
\text{ at time $s$}\}.
\teq(3.13aa)
\endalign
$$
We order the sites of $\Bbb Z^d$ in some deterministic way, say
lexicographically, and take
$$
\align
\l^*(s,u)  :=& \text{ the first $x$ in this order for which
$\cB^h\big(x,s;u, -C_5\ka(s)\big)$ occurs}\\
&\text{ and which satisfies $\langle x,u\rangle = h^*(s,u)$}.
\teq(3.13bb)
\endalign
$$
Thus, $h^*(s,u)$ is the furthest displacement in the direction of $u$ among
the $B$-particles in the process $\cP^h(u, -C_5\ka(s))$
at time $s$ and $\l^*(s,u)$ is the first site occupied by a
$B$-particle in this process at time $s$ on which this maximal
displacement is
reached. We shall write $m^*(s,u)$ for $[l^*(s,u)]^\perp$ so that we have the
orthogonal decomposition
$$
\l^*(s,u) = h^*(s,u) u + m^*(s,u).
\teq(3.13cc)
$$

The following proposition contains our principal ``subadditivity''
property. If we take $\be = \infty$, that is, if we only look at its
statement about displacements in the direction of $u$, then the
proposition says that (up to error terms) the
maximal displacement in the direction $u$ at time $s+t+C_6\ka(t)$
in the process
$\cP^h\big(u, -C_5\ka(s+t+C_6\ka(t))\big)$
is stochastically larger than the sum of two independent
such displacements, which are distributed like the maximal
displacement in $\cP^h\big(u, -C_5\ka(s)\big)$ at time $s$ and the
maximal displacement in
$\cP^h\big(u, -C_5\ka(t)\big)$ at time $t$, respectively (see
Corollary 5 for more details). The basic idea of the proof is that if
$\l^*$ is a point where $\cP^h\big(u, -C_5\ka(s)\big)$ achieves its
maximum displacement in the direction $u$ at time $s$, then we can start a new
half-space process ``near'' $\l^*$ which is nearly a copy of $\cP^h\big(u,
- C_5\ka(t)\big)$ and which is nearly independent of the first process
$\cP^h\big(u, -C_5\ka(s)\big)$. If we run the second process for $t$
units of time the sum of the displacements in the direction of $u$ in the first
and second process is nearly a lower bound for the displacement of the
original process at time $s+t$.
\proclaim{Proposition 3} Let $u \in S^{d-1},
\al, \be \ge 0$ and $\ga_s, \ga_t \in
\Bbb R^d$ orthogonal to $u$.
For any $K >0$ there exist constants $0 < C_5-C_8,s_0 < \infty$,
which depend on $K$, but are independent of  $u \in S^{d-1}$ and of
$\al, \be, \ga_s, \ga_t$, such that for
$$
s_0 \le s \le t \text{ and } t\log t \le C_7s^2
\teq(3.2cba)
$$
it holds
$$
\align
&P\big\{\cG\big(\al, \be,\ga_s + \ga_t,
\cP^h\big(u, -C_5\ka(s+t+C_6\ka(t))\big),s+t +C_6\ka(t)\big)\big\} \\
&\ge \int_{0 \le h < \infty} \int_{m \in \Bbb R^d}
P\{h^*(s,u) \in dh, m^*(s,u) - \ga_s \in dm\}\\
&\phantom{MMMMMMM}
\times P\big\{\cG\big(\al- h, \be-d,\ga_t - m,
\cP^h\big(u,-C_5\ka(t)\big),t\big)\big\}\\
&\phantom{MMMMMMMMMMMMMMMMMMMMMMMM}-C_8s^{-K-1}.
\teq(3.2)
\endalign
$$
\endproclaim
\demo{Proof}
The constants $C_i$ and $s_0$ will be fixed later on. $K_i$ will be
used to denote other auxiliary constants.
For the time being we only
do manipulations which do not depend on the specific value of the
$C_i, K_i$. We break the proof up into three steps.
The first two steps reduce the proof of \equ(3.2) to estimating the
probabilities of a number of small events. These estimates will be
carried out in step 3.
\newline
{\bf Step 1.}
Run $\cP^h(u,-C_5\ka(s))$ till time $s$.
Let $h^*(s,u) =  h \in \Bbb R,\, \l^*(s,u) = y \in \Bbb Z^d$. Set
$\overline y :=\lfloor y+4C_5\ka(t)u \rfloor$ (the meaning of
this last notation is that we take the largest integer for each
coordinate separately). Next we run the
$(u,\langle y ,u \rangle + 2C_5\ka(t))$ half-space
process starting at
the space-time point $\big(\overline y, s+ C_6\ka(t)\big)$
for $t$ units of time. This latter half-space process will be shown to
be almost an independent copy of the translate by
$\big(\overline y, s+C_6\ka(t)\big)$ of $\cP^h\big(u, -C_5\ka(t)\big)$.
Define $z(s,t)$ to be the nearest site in $\Bbb Z^d$ to $\overline y$
which is occupied at time $s+C_6\ka(t)$ by a particle which started at
time 0 in the halfspace
$\cS\big(u,\langle y,u \rangle + 2C_5\ka(t)\big)$.
\comment
To run the $(u,\langle y ,u \rangle + 2C_5\ka(t))$ half-space
process starting from $\overline y, s+ C_6\ka(t)\big)$
we use only particles which at time 0 were in
$\cS\big(u, \langle y,u \rangle + 2C_5\ka(t)\big)$.
In addition, we have to  reset at time $s+ C_6\ka(t)$ the types of all
particles in the $\big(u,\langle y,u \rangle + 2C_5\ka(t)\big)$
half-space process which are not at
$z(s,t)$ to $A$,
and further give the particles at $z(s,t)$ type $B$.
\endcomment
It will be useful to define for general $v \in \Bbb Z^d$
$$
\align
z_v = &\text{ the nearest site on $\Bbb Z^d$ to $\overline v :=
\lfloor v + 4C_5 \ka(t)u \rfloor$ which is}\\
&\text{ occupied at time $s+C_6\ka(t)$ by
a particle which started}\\
&\text{ at time 0 in }
\cS\big(u,\langle v,u \rangle + 2C_5\ka(t)\big).
\teq(3.29aa)
\endalign
$$
Thus, $z_v$ has the same relation to $v$ as $z(s,t)$ has to $y$. In
particular, $z_{y}=z(s,t)$. We can now define, still for any
$v \in \Bbb Z^d$, the sets
$$
\align
A_1(s,t) = \{x:\; & x
\text{ is occupied by one or more $B$-particles at time}\\
&\text{$s+ t+C_6\ka(t)$ in
the process }\cP^h\big(u, - C_5\ka(s+t+C_6\ka(t))\big)\},
\endalign
$$
$$
\align
A_2(s,t,v) = \{x:\;&  x
\text{ is occupied by one or more $B$-particles at time}\\
&\text{$s+t+ C_6\ka(t)$ in
the }\big(u, \langle v,u \rangle + 2C_5\ka(t)\big)\text{ half-space}\\
&\text{process starting at }(\overline v, s +  C_6\ka(t)\},
\teq(3.29bab)
\endalign
$$
and
$$
\align
A_3(t)= \{x:\; & x
\text{ is occupied by one or more $B$-particles at time $t$ in }\\
&\text{an independent copy of
the process }\cP^h\big(u, - C_5\ka(t)\big)\}.
\endalign
$$
We stress that $A_3$ is defined by means of a new copy of all
initial data and paths. It is independent of the processes we worked
with so far.

Our aim is to prove the following
two statements, and to show that they imply \equ(3.2).
The first statement is that
outside an event
of probability at most
$s^{-K-1}$ it is the case that
$$
A_1(s,t) \supset A_2(s,t,y).
\teq(3.3)
$$
The second statment is that
$$
\align
A_1(s,t)-\overline y &\text{ is at least as large as $A_3(t)$, outside }\\
&\text{an event of probability at most $s^{-K-1}$}
\teq(3.4)
\endalign
$$
(still $y = \l^*(s,u)$ in these relations).
The relation \equ(3.4) is stated somewhat imprecisely, but a precise
version will be given below (see (3.51)).
In this first step we shall reduce the proofs of \equ(3.3) and
\equ(3.4) to a number of probability estimates.

To begin with the inclusion \equ(3.3), we claim that this
holds on the intersection of the event
$$
\align
\{\langle y,u \rangle \ge 0\} \cap
\{&z(s,t)
\text{ is occupied by a $B$-particle at time}\\
&\text{$s+ C_6\ka(t)$ in }\cP^h\big(u,-C_5\ka(s)\big)\}
\teq(3.30)
\endalign
$$
with the event (see \equ(3.0xyz) for $x_0$)
$$
\{\|x_0\| \le K_3 \log s\}.
\teq(3.4ab)
$$
This follows from two applications of the monotonicity
property in Lemma C. Indeed, under \equ(3.4ab) (and $s
\ge s_1$ for a large enough $s_1$),
both the $(u,-C_5\ka(s))$ and the $\big(u, -C_5\ka(s+t+C_6\ka(t))\big)$
half-space process begin with B-particles at $x_0$.
One application of the
monotonicity property therefore gives us that (under \equ(3.4ab))
$\cP^h\big(u, -C_5\ka(t+s+C_6\ka(t))\big)$ has more $B$-particles than
$\cP^h\big(u, - C_5\ka(s)\big)$ at each space-time point, and therefore
$$
\align
A_1(s,t) \supset \{x:\; & x
\text{ is occupied by one or more $B$-particles at time}\\
&\text{$s+ t+C_6\ka(t)$ in
the process }\cP^h\big(u, - C_5\ka(s)\big)\}.
\teq(3.4bcd)
\endalign
$$
For the second application we recall that
$z(s,t)$ is
occupied at time $s+C_6\ka(t)$ by a particle which started in $\cS\big(u,
\langle y,u \rangle +2C_5\ka(t)\big)$, and in fact
is the closest occupied site to $\overline y$ with this property.
To run the $\big(u,\langle y,u \rangle +2C_5\ka(t)\big)$
half-space process starting at
$\big(\overline y,s+C_6\ka(t)\big)$ and to
find $A_2(s,t,y)$ we first
remove all particles which at time 0 were in the halfspace $\{x: \langle x,u
\rangle < \langle y,u \rangle + 2C_5\ka(t)\}$. After that we
reset to $A$ the types of all
particles not at $z(s,t)$ at time $s+C_6\ka(t)$ and give all particles
at $z(s,t)$ type $B$.
Note that in the first step all particles which do not belong to
$\cP^h\big(u, -C_5\ka(s)\big)$ are removed, since
$$
-C_5\ka(t) \le 2C_5\ka(t) \le
 \langle y,u \rangle + 2C_5\ka(t) \text{ (on \equ(3.30))} .
$$
Thus, (at time $s+C_6\ka(t)$)
after both steps, all remaining particles are also in
$\cP^h\big(u, -C_5\ka(s)\big)$, and
the particles which have type $B$, i.e., only the particles at $z(s,t)$,
also have type $B$ in $\cP^h\big(u, -C_5\ka(s)\big)$ (still on the
event \equ(3.30)).
By virtue of the monotonicity property of Lemma C,
at time $s+t+C_6\ka(t)$, any $B$-particle present in the $\big(u,
\langle y,u \rangle +2C_5\ka(t)\big)$ half-space process starting at
$\big(\overline y,s+C_6\ka(t)\big)$ is also a
$B$-particle in $\cP^h\big(u,-C_5\ka(s)\big)$. Therefore,
on the event \equ(3.30),
$$
\align
A_2(s,t,y) \subset \{x:\; & x
\text{ is occupied by one or more $B$-particles at }\\
&\text{time $s+ t+C_6\ka(t)$ in
the process }\cP^h\big(u, - C_5\ka(s)\big)\}.
\teq(3.4cde)
\endalign
$$
Combining \equ(3.4bcd) and \equ(3.4cde) gives \equ(3.3) on the
intersection of the events \equ(3.30) and \equ(3.4ab).
We postpone the
proof that this intersection indeed has probability at least $1-
s^{-K-1}$ to step 3.

To prepare for the desired precise form of \equ(3.4) we shall prove
that there exist constants $K_1$ and $s_2$ such that for
$t \ge s \ge s_2$, $\La$ any non-random subset of $\Bbb Z^d$,
and any fixed $v \in \Bbb Z^d$,
$$
P\{A_2(s,t,v)\text{ intersects }\La\} \ge P\big\{\big(\overline v +A_3(t)\big)
\text{ intersects }\La\big\} -
K_1t^{-K-d-1}.
\teq(3.33)
$$
To prove this inequality we remind the reader
that $A_2(s,t,v)$ is the collection of sites where
$B$-particles are present at time
$s+t + C_6\ka(t)$, if one starts at time
$s+C_6\ka(t)$ in the state obtained
by removing the particles which started outside
$\cS\big(u,\langle v,u \rangle +2C_5\ka(t)\big)$ at time 0, and
by resetting all particles not at
$z_v$ to type $A$, while setting the type of the particles at $z_v$ to
$B$. To find the distribution of
$A_2(s,t,v)$ we must first describe the state at time $s+C_6\ka(t)$
(after the removal of particles and resetting of types) in some more detail.
First let us look how many particles there are at the various sites,
irrespective of their type.
We began at time 0 with $N_A(w,0-)$ particles at $w$, for $w \in
\cS\big(u, \langle v,u \rangle +2C_5\ka(t)\big)$ and with 0 particles
at any $w$ outside
$\cS\big(u, \langle v,u \rangle +2C_5\ka(t)\big)$.
The $N_A(w,0-)$ are i.i.d. mean $\mu_A$
Poisson random variables. We let these particles perform their random
walks till time $s+C_6\ka(t)$. Let us write
$\wh N\big(w,s+C_6\ka(t)\big)$ for the number of particles (of either
type) at $w$ at this time. By properties of the Poisson distribution, all the
$\wh N\big(\overline v +w,s+C_6\ka(t)\big), \; w \in \Bbb Z^d$,
are still independent Poisson variables, but
$$
\align
E\wh N&\big(\overline v+w,s+C_6\ka(t)\big)\\& =
\sum_{w' \in \cS\big(u,\langle v,u \rangle +2C_5\ka(t)\big)}
\mu_A P\{S_{s + C_6\ka(t)} = \overline v+w-w'\}=:\nu(v,w,s,t).
\teq(3.34)
\endalign
$$
$z_v$ is now the nearest lattice point to $\overline v$ which is
occupied by some particle at time $s+C_6\ka(t)$. We then reset all
particles not at $z_v$ to type $A$, and the ones at $z_v$ to type
$B$. If we shift everything by $\big(-\overline v, -s
-C_6\ka(t)\big)$ (i.e., move $(w,r)$ to $\big(w-\overline v,
r-s-C_6\ka(t)\big)$), then, at $(w,0)$ we have $M(w) :=
\wh N\big(\overline v+w,s+C_6\ka(t)\big)$ particles, all of which will
be reset to type $A$, except those at $w_0 :=$ the nearest lattice
site to the
origin with $M(w)> 0$. In fact, $w_0 = z_v - \overline v$. The $M(w)$
are independent Poisson variables, and $M(w)$ has mean
$\nu(v,w,s,t)$.
It follows from the definition of $A_2(s,t,v)$ and of the $\big(u,
\langle v,u\rangle + 2C_5\ka(t)\big)$ half-space process started at
$\big(\overline v, s+C_6\ka(t)\big)$
that $A_2(s,t,v) - \overline v$ has the distribution of
$$
\{x: \text{ there is a $B$-particle at $x$ at time $t$ in
this shifted system}\}.
\teq(3.35a)
$$
For the $w \in \cS\big(u, -C_5\ka(t)\big)$, the means
$\nu(v,w,s,t)$ are close to $\mu_A$. In fact, it follows from
\equ(3.34) that for $w \in \cS\big(u, -C_5\ka(t)\big)$
and $t \ge t_0 \lor s$, for some $t_0$ (independent of $v,w,u$),
$$
\align
\mu_A \ge \nu(v,w,s,t) &\ge  \mu_A\Big
[1- \sum_{\wt w: \langle \wt w,u \rangle
\ge  C_5\ka(t)-d} P\{S_{s+C_6\ka(t)}= \wt w\}\Big]\\
&\ge \mu_A \big[1-P\{\langle S_{s+ C_6\ka(t)},u\rangle  \ge
C_5\ka(t)-d\}\big]\\
&\ge \mu_A \big[1- K_2\exp[-K_3C_5^2 \log t]\big].
\teq(3.35abc)
\endalign
$$
for some constants $K_2,K_3$ that depend on $d,D$ only (see (2.42) in
\cite {KSa} and recall that we assume \equ(3.2cba)).  From now on
we take $C_5$ so large that for large $t$
$$
\mu_A[1-t^{-K-2d-1}] \le \nu(v,w,s,t) \le \mu_A \text{ for all }w \in
\cS\big(u, -C_5\ka(t)\big).
\teq(3.36a)
$$
It suffices for this that $K_3 C_5^2 \ge K+2d+2$. We may have to
raise $C_5$ in the proof of (3.62) and (3.64) in step
3, but that can only improve the present estimates.
With such a choice of $C_5$
the distribution of the particle numbers $\{M(w): w \in
\cS\big(u, -C_5\ka(t)\big) \cap \cC(3C_1t)\}$ differs in total
variation from the distribution of an i.i.d. collection
of mean $\mu_A$ Poisson variables on
$\cS\big(u, -C_5\ka(t)\big) \cap  \cC(3C_1t)$ by at most
$$
\sum_{w \in \cS\big(u, -C_5\ka(t)\big) \cap\; \cC(3C_1t)}
\mu_At^{-K-2d-1} \le K_4 t^{-K-d-1}
\teq(3.36)
$$
for some constant $K_4 = K_4(\mu_A,d)$.

Now consider an auxiliary process which
starts at time 0 with $N_A(w,0-)$ particles only at the vertices $w
\in  \cS\big(u, -C_5\ka(t)\big) \cap \cC(3C_1t)$, and with no particles
outside this set. Let $w_0$ be the
nearest vertex in $\cS\big(u, -C_5\ka(t)\big)$
to the origin with $N_A(w_0,0-) > 0$. Take the type
of all particles not at $w_0$ equal to $A$ and the type of the
particles at $w_0$ equal to B.
 If $w_0$ lies outside
$\cS\big(u, -C_5\ka(t)\big) \cap \cC(3C_1t)$, then this auxiliary
process has never any $B$-particles. On the other hand, if $w_0 \in
\cS\big(u, -C_5\ka(t)\big) \cap \cC(3C_1t)$, then the auxiliary
process is obtained from
$\cP^h\big(u, -C_5\ka(t)\big)$ by removing at time 0 all
particles in   $\cS\big(u, -C_5\ka(t)\big) \setminus \cC(3C_1t)$.
Finally, let
$$
A_4(t) = \{x: \text{ there is a $B$-particle at $x$ at time $t$ in
this auxiliary system}\}.
$$
From our considerations above (in particular \equ(3.35a), \equ(3.36))
we have that
$$
P\{A_2(s,t,v)\text{ intersects }\La\}
\ge P\{\overline v + A_4(t) \text{ intersects }\La\} - K_4t^{-K-d-1}.
\teq(3.37)
$$
Indeed, were it not for the fact that $N_A(w,0-)$ is a Poisson
variable of mean $\mu_A$ instead of $\nu(v,w,s,t)$, the auxiliary system
would be obtained from the system in which $A_2(s,t,v)$ is computed
by translation by $\big(-\overline v,-s-C_6\ka(t)\big)$ and by removing
the particles outside $\cS\big(u, -C_5\ka(t)\big) \cap \cC(3C_1t)$. The
term $-K_4t^{-K-d-1}$ corrects for increasing the mean from $\nu(v,w,s,t)$
to $\mu_A$.

To come to \equ(3.33) we still want to prove the
inequality
$$
\align
P\{\overline v +A_4(t) \text{ intersects }\La\} \ge P\{\overline v + A_3(t)
\text{ intersects }\La\} -K_5t^{-K-d-1}.
\teq(3.38)
\endalign
$$
This follows from the fact that if, in $\cP^h\big(u, -C_5\ka(t)\big)$
all $B$-particles stay inside $\cC(2C_1t)$ during
$[0,t]$, and no particle which starts outside $\cC(3C_1t)$ at time 0
enters $\cC(2C_1t)$ during $[0,t]$, then the particles which start
outside $\cC(3C_1t)$ do not
interact with any particle, and do not cause the creation of any
$B$-particles during $[0,t]$ (compare the argument for (2.36) in \cite
{KSb}).
In these circumstances
$\cP^h\big(u, -C_5\ka(t)\big)$ has no more $B$-particles
at time $t$ than the auxiliary process, which is obtained by removing
the particles which start outside $\cC(3C_1t)$ at time 0, as described
above. Therefore
$$
\align
&\big|P\{\overline v +A_4(t) \text{ intersects }\La\}
- P\{\overline v + A_3(t)
\text{ intersects }\La\big|\\
&\le P\{w_0 \notin \cS\big(u, -C_5\ka(t)\big)
\cap  \cC(3C_1t)\}\\
&+ P\{\text{some $B$-particles in
$\cP^h\big(u, -C_5\ka(t)\big)$ leave
$\cC(2C_1t)$ during }[0,t]\}\\
&+P\{\text{in $\cP^h\big(u,-C_5\ka(t)\big)$
some particles  which start outside $\cC(3C_1t)$} \\
&\phantom{iMMMMMMMMMMMMMMMMM}\text{enter $\cC(2C_1t)$ during $[0,t]$}\}.
\teq(3.39)
\endalign
$$
The first term in the right hand side here is trivially $o(t^{-K-d-1})$
(compare \equ(3.13)).

To estimate the second term in
the right hand side of \equ(3.39) we first observe that
if $x_0$, the nearest occupied point to the origin at time 0, lies in
the half-space $\cS\big(u, -C_5\ka(t)\big)$, then at all times the
$B$-particles in $\cP^h\big(u, -C_5\ka(t)\big)$ are also $B$-particles
in $\cP^f$ (see \equ(3.0cde)). Therefore a decomposition  with
respect to the value of $x_0$ and an application of \equ(3.0b) show
that the second term in  the
right hand side of \equ(3.39) is (for $t \ge$ some $t_2$) bounded by
$$
\align
&P\{\|x_0\| \ge K_4 \log t\}
+ \sum_{\|x\| \le K_4 \log t}\mu_A P^{or}\{\text{there are $B$-particles}\\
&\phantom{MMMMMMMMMMMMMMM}\text{outside
$\cC(C_1t)$ during $[0,t]$}\}\\
&\le K_6t^{-K-d-1}
\teq(3.39aa)
\endalign
$$
(see \equ(3.0cde), \equ(3.13) and the estimates for \equ(3.18)).

The third term in the right hand side of \equ(3.39) is at most
$$
\align
\sum_{w \notin \cC(3C_1t)}& E\{N_A(w,0-)\}
P\{\sup_{r \le t}\|S_r\| \ge \|w\|- 2C_1t\}\\
&\le \sum_{w \notin \cC(3C_1t)} 8d\mu_A \exp[-K_7\|w\|]
\le K_8 t^{d-1}\exp[-K_9C_1t]
\teq(3.39abc)
\endalign
$$
(see (2.42) in \cite{KSa}). Thus \equ(3.38) and \equ(3.33) hold.

\medskip
\noindent
{\bf Step 2.} We wish to prove the following precise version of \equ(3.4):
for $t \ge s \ge s_0, t\log t \le C_7s^2$ and
for some constant $K_{10}$, independent of $s,t,u$,
$$
\align
P\{A_1(s,t)&\text{ intersects }\La\} \\
&\ge P\big\{\big(\overline{\l^*(s,u)} +A_3(t)\big)
\text{ intersects }\La\big\} -
K_{10}t^{-K-1}.
\teq(3.33abc)
\endalign
$$
To this end we define the following events for
any vector in $v \in \Bbb Z^d$:
$$
\align
&\cI_1(v) :=\Big\{\text{during $[0,s]$ in the process
$\cP^h\big(u,-C_5\ka(s)\big)$
all the  $B$-particles} \\
&\phantom{MMMMi}\text{stay in the set }
\cC(2C_1s) \cap\{x:\langle x,u \rangle
< \langle v,u \rangle + C_5\ka(t)\}\Big\},
\endalign
$$
$$
\align
\cI_2(v)  :=
\Big\{&\text{none of the  particles which were at time 0 in the half-space }\\
& \cS\big(u,\langle v,u \rangle +2C_5\ka(t)\big)
 = \{x:\langle x,u \rangle \ge
\langle v,u \rangle + 2C_5\ka(t)\}\text{ enters }\\
&\text{the set $\cC(2C_1s) \cap\{x:\langle x,u \rangle
< \langle v,u \rangle + C_5\ka(t)\}$ during }[0,s]\Big\}.
\endalign
$$
The following independence property is crucial for our argument:
Let
$\cJ(v)$ be an event which depends only on $v \in \Bbb Z^d$ and
the particles which start
in $\cS\big(u, \langle v,u \rangle + 2C_5\ka(t)\big)$
at time 0, and the {\it paths} of these particles.
Then
$$
\align
&P\{\l^*(s,u) = v,\cI_1(v), \cI_2(v), \cJ(v)\}\\
& =
P\{\l^*(s,u) = v,\cI_1(v), \cI_2(v)\}P\{ \cJ(v)|\cI_2(v)\}.
\teq(3.31)
\endalign
$$
The important feature here is that in the
last conditional probability $v$ is a constant, without relation to
$\l^*(s,u)$.
To see \equ(3.31) we note that
on the event $\cI_1 \cap \cI_2$ none of
the particles which start in $\cS\big(u, \langle v,u \rangle + 2C_5\ka(t)\big)$
coincides with any $B$-particle during $[0,s]$. Therefore, changing
the paths of any of the particles which start in $\cS\big(u, \langle v,u
\rangle + 2C_5\ka(t)\big)$ has no influence on the types of any of the
other particles during $[0,s]$
(and of course no influence on the paths of these
other particles), as long as we stay on $\cI_1 \cap \cI_2$
(compare the argument for (2.36) in \cite {KSb}). In particular,
$$
P\{\l^*(s,u)=v,\cI_1(v)| \cI_2(v), \cJ(v)\}
= P\{\l^*(s,u) = v,\cI_1(v)| \cI_2(v)\}.
$$
This is clearly equivalent to \equ(3.31).

We take
$$
\cJ(v) = \{A_2(s,t,v) \text{ intersects } \La\},
$$
where
$\La$ is some non-random set in $\Bbb Z^d$. By definition, $A_2(s,t,v)$
depends only on $v$ and the particles
which start in the half-space $\cS\big(u,
\langle v,u \rangle +2C_5\ka(t)\big)$. Thus also $\cJ(v)$ depends only
on $v$ and this last collection of particles and their paths.
(This is true despite the fact that we talk about $B$-particles in the
definition \equ(3.29bab). Indeed, these are $B$-particles in $(u,
\langle v,u \rangle +2C_5\ka(t)\big)$ half-space process,
started at $\big(\overline v,s+C_6\ka(t)\big)$, and the types of these
particles are reset at time $s+C_6\ka(t)$ and after that do not depend
on particles which started outside $\cS\big(u,\langle v,u \rangle
+2C_5\ka(t)\big)$.)
With this choice of $\cJ$ we obtain from \equ(3.31) for every fixed $v$
$$
\align
&P\{\l^*(s,u) = v,\cI_1(v), \cI_2(v), A_2(s,t,v) \text{ intersects }
\La\}\\
& \ge
P\{\l^*(s,u) = v,\cI_1(v), \cI_2(v)\}
\big[P\{A_2(s,t,v) \text{ intersects }\La\}-P\{\cI_2^c(v)\}\big]^+.
\teq(3.32)
\endalign
$$

We shall show in step 3 that for suitable choice of constants $0 <
K_i= K_i(K,d) < \infty$, independent of $s,u$ and $v$,
it is the case that for the process $\cP^h\big(u, -C_5\ka(s)\big)$
$$
P\{\l^*(s,u) \in \cC(2C_1s),\; \cI_1^c(\l^*(s,u))\} \le K_{11}s^{-K-1},
\teq(3.45)
$$
$$
P\{\cI_2^c(v)\} \le K_{11}s^{-K-d-1},
\teq(3.46)
$$
and
$$
P\{\text{\equ(3.30) (with $y = \l^*(s,u)$) fails
or \equ(3.4ab) fails}\} \le K_{11}s^{-K-1}.
\teq(3.47)
$$
In the remainder of this step we only show how to complete the proof
of \equ(3.33abc) and the proposition from the estimates \equ(3.45)-\equ(3.47).
 To this end we apply \equ(3.32).
By using \equ(3.32), \equ(3.46) and \equ(3.33) we get
$$
\align
&P\{\l^*(s,u) = v, A_2(s,t,\l^*(s,u))
\text{ intersects }\La\}\\
&\ge  P\{\l^*(s,u) = v, \cI_1(v), \cI_2(v),
A_2(s,t,v) \text{ intersects }\La\}\\
&\ge\Big[ P\{\l^*(s,u) = v\} - P\{\l^*(s,u)= v,  \cI_1^c(v)\} - P\{\l^*(s,u)
= v, \cI_2^c(v)\}\Big]\\
&\phantom{MMMMMMM} \times
\big[P\{ A_2(s,t,v)\text{ intersects }\La\}-P\{\cI_2^c(v)\}\big]^+\\
&\ge P\{\l^*(s,u) = v\}P\{A_2(s,t,v)\text{ intersects }\La\}\\
&\phantom{MMMMMMM} - P\{\l^*(s,u) = v, \cI_1^c(v)\} -2K_{11}s^{-K-d-1}\\
&\ge  P\{\l^*(s,u) = v\}P\{\overline v + A_3(t)
\text{ intersects }\La\}\\
&\phantom{MMMMMMM} - P\{\l^*(s,u) = v, \cI_1^c(v)\}
-(2K_{11}+K_4 + K_5)s^{-K-d-1}.
\teq(3.48)
\endalign
$$
Now recall that
$A_1(s,t) \supset A_2\big(s,t,\l^*(s,u)\big)$
on the intersection of \equ(3.30) and \equ(3.4ab).
Summing \equ(3.48) over
all $v \in \cC(2C_1s)$, and using \equ(3.45) and \equ(3.47),
therefore gives
$$
\align
&P\{A_1(s,t) \text{ intersects }\La\text{ and \equ(3.30),
\equ(3.4ab) occur}\}\\
&\ge P\{A_2(s,t,\l^*(s,u))\text{ intersects }\La\} -
P\{\text{\equ(3.30) or \equ(3.4ab) fails}\}\\
&\ge \sum_{v \in \cC(2C_1s)}P\{\l^*(s,u) = v\}
P\{\overline v + A_3(t) \text{ intersects }\La\}
- K_{12}s^{-K-1}.
\teq(3.49)
\endalign
$$
Finally, since $\l^*(s,u)$ is the location of a
$B$-particle at time $s$ in $\cP^h\big(u, -C_5\ka(s)\big)$,
 we have, essentially as in the estimate for $P\{\cE_{k,4}^c\}$
in \equ(3.18) and the lines following it, or the estimate of the
second term in the right hand side of \equ(3.39)
$$
\align
&P\{\l^*(s,u) \notin \cC(2C_1s)\}
\le P\{\|x_0\| > C_1s \wedge C_5\ka(s)/\sqrt d\}\\
& + \sum_{\|x\|\le C_1s}\mu_AP^{or}\{\text{there
is a $B$-particle outside $\cC(C_1s)$ at time $s$}\}
\le K_{13} s^{-K-1}.
\teq(3.50)
\endalign
$$
Consequently
$$
\align
P\{A_1(s,t) \text{ intersects }\La\} \ge
 \sum_{v \in \Bbb Z^d}P\{\l^*(s,u) = v\}
&P\{\overline v + A_3(t) \text{ intersects }\La\}\\
&- (K_{12}+K_{13})s^{-K-1}.
\teq(3.51)
\endalign
$$
This is the desired \equ(3.33abc).

\equ(3.2) is just the
special case of \equ(3.51) with $\La = \Ga(\al, \be, \ga_s+\ga_t)$. Indeed,
$\{A_1(s,t) \text{ intersects }\La\}$ is the event that there is
a $B$-particle in $\La$ at time $s+t+C_6\ka(t)$ in the process
$\cP^h\big(u, -C_5\ka(s+t+C_6\ka(t))\big)$. For $\La =
\Ga(\al, \be,\ga_s+\ga_t)$ this
event is also denoted by $\cG\big(\al, \be, \ga_s+\ga_t,
\cP^h\big(u, -C_5\ka(s+t+C_6\ka(t))\big), s+t+C_6\ka(t)\big)$.
Thus, the left hand sides of \equ(3.2) and \equ(3.51) are the same for
this choice of $\La$. We leave it to the reader to check that the
right hand side of \equ(3.51) is at least as large as the right hand
side of \equ(3.2), provided we choose $C_8 \ge K_{12} + K_{13}$.
\medskip
\noindent
{\bf Step 3.} Here we prove the relations \equ(3.45)-\equ(3.47). Note
that \equ(3.47)
also supplies the missing estimates for \equ(3.3), to
wit, $P\{\text{\equ(3.30)}$ and \equ(3.4ab)
hold\} $\ge 1-s^{-K-1}$.

Now we start on \equ(3.45). First
$$
P\{\text{in $\cP^h\big(u, -C_5\ka(s)\big)$
some $B$-particle leaves $\cC(2C_1s)$ during} [0,s]\}
= O\big(s^{-K-d-1}\big)
\teq(3.49aa)
$$
(see \equ(3.39aa)). In addition, by definition
of $l^*(s,u),\;\langle l^*(s,u),u \rangle = h\big(s,u,-C_5\ka(s)\big)$. Thus
$$
\align
&P\big\{\text{in $\cP^h\big(u, -C_5\ka(s)\big)$, during $[0,s]$ all
$B$-particles stay in $\cC(2C_1s)$,
but some }\\
&\phantom{MMMMMMM}\text{of them
leave $\{x:\langle x,u \rangle < \langle \l^*(s,u),u
\rangle + C_5\ka(s)\}$ }\big\}\\
&\le P\big\{\text{in $\cP^h\big(u, -C_5\ka(s)\big)$,
at some time $r \le s$ there are }\\
&\phantom{MiM}\text{$B$-particles at some $v \in \cC(2C_1s)$ with}\\
&\phantom{MMi} \langle v,u \rangle \ge
\max\{\langle x,u \rangle: \text{ there is a
$B$-particle at $x$ at time }s\}+C_5\ka(s)\big\}.
\teq(3.49bb)
\endalign
$$
This last event can happen only if some
$B$-particle reaches a vertex $v \in \cC(2C_1s)$
before time $s$ and then this particle moves
to some $x$ at time $s$ with $\langle x,u
\rangle < \langle v,u \rangle -C_5\ka(s)$. The probability
that such particle started outside $\cC(3C_1s)$ is bounded by the
third term in the right hand side of \equ(3.39), with $t$ replaced by $s$.
Therefore, the right hand
side of \equ(3.49bb) is at most
$$
\align
&\text{(third term in right hand side of \equ(3.39) with $t$ replaced
by $s$)}\\
&\qquad +\sum_{w \in
\cC(3C_1s)} E\{N_A(w,0-)P\{\sup_{0 \le r_1,r_2 \le s}\|S_{r_1} -
S_{r_2}\|\ge C_5\ka(s)/\sqrt d\}\\
&\le K_8s^{d-1}\exp[-K_9C_1s] + K_{14} (3C_1s)^d\exp[-K_{15}C_5^2\log s],
\teq(3.49cc)
\endalign
$$
by \equ(3.39abc) and by (2.42) in \cite {KSa}.
Together with \equ(3.49aa) this proves that we can take $C_5$ so large
that \equ(3.45) holds. As observed after \equ(3.36a) we can even
choose  $C_5$ so that \equ(3.36a) is also valid.
Once we have chosen $C_5$ we fix
$$
C_6 = \frac {16 C_5}{C_2}.
\teq(3.49dd)
$$

As for \equ(3.46), we have
$$
\align
&P\{\cI_2^c(v)\}
\le \sum_{w \in \cS(u,\langle v,u \rangle +2C_5\ka(s))}E\{N_A(w,0-\}\\
&\phantom{MMMMMMMM}
\times P\{\sup_{r \le s}\|S_r\| \ge C_5\ka(t) \lor
\big(\|w\|- 2C_1s\big)\}.
\teq(3.49ee)
\endalign
$$
We leave it to the reader to show that this is $O\big(s^{-K-d-1}\big)$
for $t \ge s$ and large enough $C_5$ (again by (2.42) in \cite {KSa}).

Finally, to prove \equ(3.47), we note first that
$P\{\text{\equ(3.4ab) fails}\} = O(s^{-K-1})$, provided $K_3=
K_3(\mu_A,d)$ is taken large enough, just as in \equ(3.13).
Next,
$$
\align
P\{\langle \l^*(s,u),u \rangle < 0\} &\le
P\{h\big(s,u, -C_5\ka(s)\big) \le C_4s\}\\
& \le
[C_5\ka(s)]^{-2K-2}\le s^{-K-1}
\teq(3.52)
\endalign
$$
for large $s$, by virtue of Lemma 2 with $K$ replaced by $2K+2$.
Lastly, we have to show that for the choice of $C_6$ in \equ(3.49dd)
$$
\align
&P\{\text{$z(s,t)$
is not occupied by a $B$-particle in }
\cP^h\big(u,-C_5\ka(s)\big)\text{ at
time }s+ C_6\ka(t)\}\\
&\le P\{\text{$z(s,t)$
is not occupied by a $B$-particle in }\cP^h\big(u,-C_5\ka(s)\big)\text{ at
time}\\
&\phantom{MMMMMMMMMMMM}\text{$s+ C_6\ka(t)$,
but }z(s,t) \in \l^*(s,u) + \cC(C_2C_6\ka(t)/2\}\\
&\phantom{M}+ P\{z(s,t) \notin \l^*(s,u) + \cC(C_2C_6\ka(t)/2\}\\
& = O\big(s^{-K-1}\big).
\teq(3.53)
\endalign
$$
The first inequality here is obvious. The bound $O(s^{-K-1})$ for the middle
member of \equ(3.53) is
formulated as a separate lemma, because the same argument
will be needed once more in the next section. To see that \equ(3.53)
follows from Lemma 4 below, recall that
$z(s,t)$ is occupied at time $s+ C_6\ka(t)$ by some particle
which started at time 0 in
$\cS\big(u,\langle \l^*(s,u),u \rangle + 2C_5\ka(t)\big)$
(by definition of $z(s,t)$). In particular there is some particle at
$z(s,t)$ at time $s+C_6\ka(t)$, so that $z(s,t)$ is occupied in
$\cP^f$ at time $s+C_6\ka(t)$.
Also, $\l^*(s,u)$ is occupied by at least one $B$-particle in
$\cP^h\big(u, -C_5\ka(s)\big)$ at time $s$.
So Lemma 4 with $\wt s = s +C_6\ka(t)$ and $y(s) = \l^*(s,u)$ (and
$C_6$ as in \equ(3.49dd)) shows that
the middle member of \equ(3.53) is at most
$$
\align
&P\{\l^*(s,u)  \notin \cC(2C_1s)\} +P\{\langle \l^*(s,u),u)
\rangle < C_4s/2\} + s^{-K-1}\\
&\phantom{MMMMMMMM}+ P\{z(s,t) \notin \l^*(s,u) +
\cC\big(C_2C_6\ka(t)/2\big)\}.
\teq(3.53efg)
\endalign
$$
Note that we used the
second part of condition \equ(3.2cba) here; we have to choose $C_7$ small
enough to make sure that (3.72)
 holds for $\wt s - s = C_6\ka(t)$.
 The first two terms in \equ(3.53efg)
are $O(s^{-K-1})$, by virtue of \equ(3.49aa)
and \equ(3.52). The fourth term is bounded by
$$
\align
&P\{z(s,t) \notin \l^*(s,u) + \cC\big(C_2C_6\ka(t)/2\big)\}
\le P\{\|z(s,t) - \overline{\l^*(s,u)}\| >  4C_5\ka(t)-1\}\\
&\text{(because $C_2C_6/2 = 8C_5$ and }\|\overline
{\l^*(s,u)} - \l^*(s,u)\| \le 4C_5\ka(t))+1)\\
&\le P\{\l^*(s,u) \notin \cC(2C_1s)\}+ P\{\l^*(s,u)  \in \cC(2C_1s),
\text{ and none of the sites}\\
&\phantom{MMMM}\text{in $\overline{\l^*(s,u) }
+\cC\big(4C_5\ka(t)-1\big)$ are occupied at time $s+ C_6\ka(t)$ by}\\
&\phantom{MMMM}\text{a particle which started in
$\cS(u,\langle \l^*(s,u) ,u \rangle + 2C_5\ka(t))$}\}.
\teq(3.53dcb)
\endalign
$$
We already saw in \equ(3.50) that the first term in the right hand side is
$O\big(s^{-K-1}\big)$. As for the second term in the
right hand side, this is by a decomposition with respect to the
possible values of $\l^*(s,u)$,  analogously to
\equ(3.0c), at most
$$
\align
K_{16} &\sum_{v \in \cC(2C_1s)}
P\{\text{none of the sites in $\overline v+\cC\big(4C_5\ka(t)-1\big)$
is occupied at time}\\
&\text{$s+C_6\ka(t)$ by a particle which started in
$\cS\big(u,\langle v,u\rangle +2C_5\ka(t)\big)$}\}.
\teq(3.53bcd)
\endalign
$$
However, the numbers of particles at sites $\overline v + w$ at time
$s+C_6\ka(t)$ which started in $\cS(u,\langle v,u\rangle +2C_5\ka(t))$
are independent Poisson variables with means $\nu(v,w,s,t)$ given
in \equ(3.34). By the estimate \equ(3.35abc) we have $\nu(v,w,s,t)$
$\ge \mu_A/2$ for $\langle w,u \rangle \ge 0$ and all $v$
(and $t$ large enough). Therefore
\equ(3.53bcd) is at most $K_{17}s^d\exp[-K_{18}\ka^d(t)\mu_A]$.
This
proves the bound $O\big(s^{-K-1}\big)$ in \equ(3.53),
and therefore \equ(3.47) is reduced to the
next lemma.

Roughly speaking, the next lemma guarantees that if a
certain vertex $y(s)$ has a $B$-particle in the
half-space process $\cP^h\big(u,-C_5\ka(s)\big)$ at a time $s$,
then a little later all
occupied sites ``near'' $y(s)$ will actually have a $B$-particle
in $\cP^h\big(u,-C_5\ka(s)\big)$.

\enddemo

\proclaim{Lemma 4}
Let $s, \wt s$ and $t$  be such that
$$
s \le t
\teq(6.0aaaa)
$$
and
$$
\frac{16C_5}{C_2}\ka(t) \le \wt s - s  \le \frac{C_4}{8C_1}s.
\teq(6.0aa)
$$
Let $u \in S^{d-1}$ be fixed and let $y(s) \in \Bbb Z^d$ be a random
point (that is, $y(s)$ may depend on the sample point $\si$).
Define the event $\cK(y)$ by
$$
\align
\cK(y) := \{&\text{there exists a site
$z \in y+\cC\big(C_2(\wt s -s)/2\big)$ such that}\\
&\text{at time $\wt s$, $z$ is occupied in $\cP^f$, but
is not occupied}\\
&\text{by a $B$-particle in $\cP^h\big(u, -C_5 \ka(s)\big)$}\}.
\teq(6.1)
\endalign
$$
Then for each $K > 0$ there exists an $s_1 = s_1(K)$ (independent of
 $u$) such that
$$
\align
P\{&\cB^h\big(y(s),s;u, -C_5\ka(s)\big) \cap \cK(y(s))\}\\
&\le P\{y(s) \notin \cC(2C_1s)\}+P\{\langle y(s), u \rangle < \frac 12C_4s\}
+ s^{-K-1}
\teq(6.2)
\endalign
$$
for $s \ge s_1$ (see \equ(3.0f) for $\cB^h$).
\endproclaim

\demo{Proof} Assume that the space-time point $(y,s)$ is occupied by some
particle in $\cP^h\big(u, -C_5 \ka(s)\big)$. We can then
define the following auxiliary events:
$$
\aligned
\cK_1(y)&:= \{\text{there exists a site
$z \in y+\cC\big(C_2(\wt s -s)/2\big)$
such that $(z, \wt s)$}\\
&\phantom{MM}\text{is occupied in $\cP^f$, but
is not occupied in $\cP^h(u, -C_5 \ka(s))$}\},\\
\cK_2(y)&:= \{\text{there exists a site $z \in y+\cC\big(C_2(\wt s-s)/2\big)$
such that $(z, \wt s)$}\\
&\phantom{MM}\text{is occupied by an $A$-particle
in the full-space process}\\
&\phantom{MM}\text{starting at $(y,s)$}\},\\
\cK_3(y)&:= \{\text{there exists a site $z \in
y+\cC\big(C_2(\wt s -s)/2\big)$
such that $(z, \wt s)$}\\
&\phantom{MM}\text{is occupied by a $B$-particle in the full-space process
starting}\\
&\phantom{MM}\text{at $(y,s)$, but occupied by an $A$-particle in
the $(u, -C_5 \ka(s))$}\\
&\phantom{MM}\text{half-space process starting at $(y,s)$}\},\\
\cK_4(y) &:= \{\text{in the full-space process starting at $(y,s)$ some
$B$-particles}\\
&\phantom{MM}\text{leave $y+\cC\big(2C_1(\wt s - s)\big)$ during $[s,
 \wt s]$}\},\\
\cK_5(y) &:= \{\text{some particles which start outside
$\cS\big(u,-C_5 \ka(s)\big)$
at time 0,}\\
&\phantom{MM}\text{enter $y+\cC\big(2C_1(\wt s -s)\big)$ during $[s,
  \wt s]$}\}.
\endaligned
$$
We shall first show that
$$
\cB^h\big(y,s;u, -C_5\ka(s)\big) \cap\cK(y)
\subset \bigcup_{i=1}^3 \cK_i(y) \text{ and } \cK_3(y) \subset
\cK_4(y) \cup \cK_5(y),
\teq(6.3)
$$
and then estimate $P\{y(s) \in \cC(2C_1s),\langle y(s),u \rangle \ge C_4s/2,
\cK_i(y(s))\}$ for $1 \le i \le 5$.
To prove the first part of \equ(6.3), consider a sample point for
which $\cB^h\big(y,s;u, -C_5\ka(s)\big) \cap \cK(y)$ occurs and let
$z$ be a site in $y+\cC\big(C_2(\wt s-s)/2\big)$ such that $(z,\wt s)$
is occupied in $\cP^f$, but is not occupied by a
$B$-particle in $\cP^h\big(u,-C_5\ka(s)\big)$. Then
it may be that $(z, \wt s)$ is not occupied at all in
$\cP^h\big(u,-C_5\ka(s)\big)$. This would mean that
$\cK_1(y)$ occurs. If this fails, then $(z,\wt s)$ is occupied in
$\cP^h\big(u,-C_5\ka(s)\big)$, necessarily by an
$A$-particle. We claim that
$(z,\wt s)$ is then also occupied by an $A$-particle in
the $\big(u,-C_5\ka(s)\big)$ half-space process starting at $(y,s)$.
This is so, because starting at $(y,s)$ does not remove any particles,
but it may change some types. But on $\cB^h\big(y,s;u,
-C_5\ka(s)\big)$, $y$ has already at least one $B$-particle at time $s$ in
$\cP^h(u, -C_5\ka(s)\big)$. Thus the resetting
at time $s$ only changes some types from $B$ to $A$, and since
$z$ already has type $A$ at time $\wt s$ in $\cP^h\big(u, -C_5\ka(s)\big)$,
it will (by Lemma C) also have type $A$ at time $\wt s$ in the $\big(u,
-C_5\ka(s)\big)$ half-space process started at $(y,s)$, as claimed.
$(z, \wt s)$ is also occupied in the full-space
process starting at $(y,s)$ (since it is occupied in the
full-space process, starting at $(\bold 0, 0)$). The type at $(z,\wt s)$ in
this process may be $A$, in which case $\cK_2(y)$ occurs, or $B$, in
which case $\cK_3(y)$ occurs.
This proves the first inclusion in \equ(6.3).

The second part of \equ(6.3) follows from the argument given for
\equ(3.38). $\cK_3(y)$ requires that there are particles in
$y+\cC\big(C_2(\wt s - s)/2\big)$
which have different types at time $\wt s$ in the full-space
and in the $(u, -C_5\ka(s))$ half-space process, both starting at
$(y,s)$.
This means that in the full-space process starting at $(y,s)$ the
type of some particle which is in $y+\cC\big(C_2(\wt s - s)/2\big)$
at time $\wt s$ is influenced by particles which
started outside $\cS\big(u, -C_5\ka(s)\big)$ at time 0. However, this
can happen only if in the full-space process starting at $(y,s)$,
these particles meet some $B$-particles
during $[s, \wt s]$. In turn, this can happen only if $\cK_4(y)$
or $\cK_5(y)$ occurs. This proves the second inclusion in \equ(6.3).

Our next task is to find bounds for
$$
P\{y(s) \in \cC(2C_1s),
\langle y(s), u \rangle
\ge C_4s/2, \cB^h\big (y(s),s;u, -C_5\ka(s)\big),
\cK_i\big(y(s)\big)\}, i =1,2,4,5.
$$
For $i=1$ we have
$$
\align
&P\{y(s) \in \cC(2C_1s),\langle y(s), u \rangle \ge C_4s/2, \cK_1(y(s))\}\\
& \le \sum_{w \notin \cS\big(u, -C_5 \ka(s)\big)} \; \sum \Sb
\langle z,u \rangle \ge C_4s/2-C_2(\wt s - s)/2  \\
z \in \cC\big((2C_1+C_4) s\big)\endSb
EN_A(w,0-)P\{w+S_{\wt s} = z\} \\
&\le \sum_{w \in \cC\big((4C_1 +2C_4)s\big)}
\mu_A P\{\|S_{\wt s}\| \ge C_4s/(4\sqrt d)\}
+ \sum_{w \notin \cC\big((4C_1 +2C_4)s\big)}
\mu_AP\{\|S_{\wt s}\|\ge \|w\|/2\}\\
&\le s^{-K-1}
\teq(6.4)
\endalign
$$
for all $s \ge$ some $s_1= s_1(K)$.
In the first inequality we used that
$\|z\| \le \|z-y(s)\| + \|y(s)\| \le C_2(\wt s-s)/2 + 2C_1s
\le (2C_1+C_4)s$, by virtue of \equ(6.0aa) and the inequality $C_2 \le C_1$.
In the second inequality we used
that for the summands here we have $\|w-z\| \ge \langle (z-w),u
\rangle/\sqrt d \ge [C_4s/2 -C_2(\wt s - s)/2 + C_5 \ka(t)]/\sqrt d
\ge C_4s/(4\sqrt d)$.
For the third inequality we use $\wt s \le
\big(1+C_4/(8C_1)\big)s$ plus (2.42) in \cite {KSa}; compare
\equ(3.18abc).

Next, we remind the reader that $P^{or}$ is the probability measure
governing the original model, in which one $B$-particle is added at
the origin at time 0. In this notation we have, by \equ(3.0c) and \equ(2.3),
$$
\align
&P\{y(s) \in \cC(2C_1s),\cB^h\big(y(s),s;u, -C_5\ka(s)\big), \cK_2(y(s))\}\\
&\le K_{19} s^d P^{or}\{\text{there exists a $z \in
\cC\big(C_2(\wt s-s)/2\big)$ which is}\\
&\phantom{mMMMMM}\text{occupied by an $A$-particle at time $\wt s-s$}\}\\
&\le 2s^{-K-1}.
\teq(6.5)
\endalign
$$

Again by \equ(3.0c)
$$
\align
&P\{y(s) \in \cC(2C_1s), \cK_4(y(s))\}\\
&\le K_{19} s^d P^{or}\{\text{some
$B$-particles leave $\cC\big(2C_1(\wt s- s)\big)$ during }[0,\wt s-s]\}\\
&\le s^{-K-1} \text{ (by \equ(3.39aa)).}
\teq(6.6)
\endalign
$$

Finally,
$$
\align
&P\{ y(s) \in \cC(2C_1s), \langle y(s), u \rangle
\ge C_4s/2, \cK_5(y(s))\}\\
&\le \sum_{w: \langle w,u \rangle < -C_5\ka(s)} \; \sum \Sb
v \in \cC\big((2C_1+C_4)s\big)\\ \langle v,u \rangle \ge C_4s/4 \endSb
 EN_A(w,0-)P\{w+ S_r = v \text{ for some }r \le \wt s\}\\
&\le \mu_A \sum_{w \in \cC\big((4C_1+2C_4)s\big)}
 P\{\sup_{r \le \big(1+C_4/(8C_1)\big)s}\|S_r\| \ge C_4s/(4\sqrt d)\}\\
&\phantom{MMMMMMMM}+
\mu_A\sum_{w \notin \cC\big((4C_1+2C_4)s\big)}
P\{\sup_{r \le \big(1+C_4/(8C_1)\big)s}\|S_r\| \ge \|w/2\|\}\\
&\le s^{-K-1} \text{ (compare \equ(6.4))}.
\teq(6.7)
\endalign
$$
Together with \equ(6.3) these estimates prove \equ(6.2).
\qedsymbol
\enddemo

\proclaim{Corollary 5}
For every  unit vector $u$ there exists a constant $\la(u) \in [C_4,
2\sqrt d C_1]$ such that
$$
\lim_{t \to \infty} \frac 1t h^*(t,u) = \la(u)
\text{ almost surely and in $L^p$ for all }p >0.
\teq(3.60)
$$
($t$ runs through the reals here).
Moreover, for each  $\eta > 0$ there exist an exponentially increasing
sequence $\{n_1 < n_2 < \dots\}=
\{n_1(\eta) < n_2(\eta) < \dots\}$ (independent of $u$) such that
$$
1 < \frac {n_{j+1}}{n_j} \le 1+\eta,\quad j \ge 1,
$$
and such that every $\ep > 0$,
$$
\sum_{k=0}^\infty P\Big\{\Big|\frac1{n_k}h^*(n_k,u) - \la(u)\Big| > \ep\Big\}
< \infty.
\teq(3.61)
$$
Finally, for given $\eta > 0$ and  $\eta' \le \eta$ and
 $\{n_j(\eta)\}$ corresponding to $\eta$, we can
choose the $n_j(\eta')$ such
that for some $j_0 < \infty, \{n_j(\eta')\} \supset \{n_j(\eta):j \ge j_0\}$.
\endproclaim
\demo{Proof} The basis for this proof is \equ(3.2) with $\be =
\infty$. Since $\Ga(\al, \infty,\ga, u) = \{x \in \Bbb R^d:\langle x,u
\rangle \ge \al\} = \cS(u,\al)$, we have
$$
\cG(\al, \infty, \ga,\cP,t)= \{\text{in $\cP$, at time $t$,
there is a $B$-particle at some $x$ with }\langle x,u\rangle \ge
\al\}.
$$
In particular
$$
\align
\cG\big(\al, \infty, \ga, \cP^h(&u, -C_5\ka(s+t+C_6\ka(t)),s+t+
C_6\ka(t)\big)\\
& = \{h\big(s+t+C_6\ka(t),u, -C_5\ka(s+t+C_6\ka(t))\big) \ge \al\}\\
& =\{h^*(s+t+C_6 \ka(t),u) \ge \al\}.
\endalign
$$
Similarly,
$$
\cG\big(\al, \infty, \ga, \cP^h\big (u,-C_5\ka(t)\big),t\big) =
\{h^*(t,u) \ge \al\}.
$$
Thus, \equ(3.2) with $\be = \infty$ says that, under \equ(3.2cba),
$$
P\{h^*(s+t+C_6\ka(t),u) \ge \al\}
\ge P\{h_1^*(s,u) +h^*_2(t,u) \ge \al\} -C_8s^{-K-1},
\teq(3.66)
$$
where $h_1^*(s,u)$ and $h^*_2(t,u)$, are independent copies of
$h^*(s,u)$ and $h^*(t,u)$, respectively.

The Corollary will be derived from this relation by more or less
standard subadditivity  techniques. To apply these techniques we
first derive some simple properties
of $h^*(s,u)$. The first is the following tail estimate:
$$
P\{h^*(s,u) \ge \al\} +P\{\|m^*(s,u)\| \ge \al\}
\le \exp[-K_1\al] \text{ for }\al \ge 2\sqrt d
C_1s,\quad s \ge s_3,
\teq(3.62)
$$
where $s_3= s_3(K)$ is some constant independent of $\al, u$.
The second and third property are semi-continuity properties in $s$, namely
$$
\align
&P\{\inf_{r \le t} h^*(s+r,u) - h^*(s,u) \le -\al\}\\
&\le K_3s^{-K}+ P\{\|\sup_{r \le t} S_r\| \ge \al\} \le
K_3s^{-K}+ 8d \exp\Big[ -\frac {K_2\al^2}{t+\al}\Big],
\;\al \ge 0,
\teq(3.67)
\endalign
$$
and
$$
\align
&P\{\sup_{r \le t} h^*(s+r,u) - h^*(s+t,u) \ge \al\}\\
&\le K_3s^{-K}+
K_4(s+t)^d \exp\Big[ -\frac {K_2\al^2}{t+\al}\Big],
\;\al \ge 0.
\teq(3.67bb)
\endalign
$$

To prove \equ(3.62) take  $\al \ge 2 \sqrt d C_1s$. Since
$\langle x,u \rangle \le \|x\|_2 \le \sqrt d \|x\|$, as well as
$\|x^\perp\| \le  \|x\|_2 \le \sqrt d \|x\|$, the left hand
side of \equ(3.62) is then bounded by
$$
\align
&2P\{\text{in $\cP^h\big (u, -C_5\ka(s) \big)$ there is a
$B$-particle outside $\cC(\al/\sqrt d)$}\\
&\phantom{MMMMMMMMMMMMM}\text{at some time during } [0,s]
\subset \big[0, [2\sqrt d C_1]^{-1}\al\big] \}\\
&\le 2P\{\text{nearest site $x$ to $\bold 0$ in $\cS\big(u,
-C_5\ka(s)\big)$ with $N_A(x,0-) > 0$}\\
&\phantom{MMMMMMMMMMMMM}\text{lies outside }\cC\big(\al/(2\sqrt d)\big)\}\\
&+ \sum_{x \in \cC\big(\al/(2\sqrt d)\big)} \mu_A
4E^{or}\{\text{number of  $B$-particles outside }\cC\big(\al/(2\sqrt d)\big)\\
&\phantom{MMMMMMMMMMMMM}\text{at time } [2\sqrt d C_1]^{-1}\al\big]\},
\endalign
$$
by an application of \equ(3.0b) and the argument following \equ(3.18),
very much as in \equ(3.39aa). The
first term in the right hand side here is at most $2\exp[-K_5\al^d]$,
and the second term is at most $K_6 \al^d\exp\big[-[2\sqrt d
C_1]^{-1}\al\big]$, by virtue of Theorem 1 in \cite {KSb}. Thus
\equ(3.62) holds.

The argument for \equ(3.67) is basically already given in
\equ(3.12a). Moreover, it is similar to, but simpler than the proof of
\equ(3.67bb) so we only prove the latter. If $h^*(s+t,u)= h$, then
all $B$-particles in $\cP^h\big(u,-C_5\ka(s+t)\big)$ are located in
$\{x:\langle x,u \rangle
\le h\}$ at time $s+t$. If for some $0 \le r \le t, h^*(s+r,u) \ge h +
\al$, then there is some $B$-particle $\rho$ in
$\cP^h\big(u,-C_5\ka(s+r)\big)$ in
$\{x:\langle x,u \rangle \ge h+\al\}$ at time $s+r$. This $\rho$ is
also a particle present in $\cP^h\big(u,-C_5\ka(s+t)\big)$
and even of type $B$ in
$\cP^h\big(u,-C_5\ka(s+t)\big)$ at time $s+t$, provided $\|x_0\| \le C_5\ka(s)/\sqrt
d$ (see \equ(3.0cde)). Thus in this
case $\rho$ moved over a
distance at least $\al/\sqrt d$ during $[s+r,s+t]$.
Therefore, the left hand side of
\equ(3.67bb) is at most
$$
\align
&P\{\|x_0\| > C_5\ka(s)/\sqrt d\}\\
&+P\{\text{some particle which starts outside $\cC\big(3C_1(s+t)\big)$
 becomes}\\
&\phantom{MMMMMMMMMMMMMMM}\text{a
$B$-particle in $\cP^f$ before time $s+t$}\}\\
&+\sum_{x \in \cC\big(3C_1(s+t)\big)} \mu_A P\{\sup_{r \le t}
\|S_r - S_t\| \ge
\al/\sqrt d\}\\
&\le K_7s^{-K} + K_8(s+t)^d\exp\Big[ -\frac {K_9\al^2}{d(t+\al)}\Big]
\endalign
$$
(see \equ(3.39)-\equ(3.39abc), as well as (2.42) in \cite{KSa}).

We can now proceed with subadditivity arguments.

\noindent
We introduce the random variables
$$
X(s) = [2 \sqrt d C_1 s - h^*(s,u)]^+
$$
and the deterministic quantities  $Y(t) = 2\sqrt d C_1C_6\ka(t)$,
and let $X'(t)$ be a copy of $X(t)$ which is independent of $X(s)$.
Then \equ(3.66) shows that, under \equ(3.2cba), these random variables satisfy
$$
P\{X\big(s+t+C_6\ka(t)\big) \le \be\} \ge P\{X(s) + X'(t) + Y(t) \le \be\}
- C_8s^{-K-1}
\teq(3.68)
$$
for $\be \ge 0$. Of course this inequality also holds trivially for
$\be < 0$. This is very close to the principal hypothesis of
the lemma on p. 674 of \cite{Ha} but we have to do some
extra work because of the $C_6\ka(t)$ which appears in the argument on
the left hand side of \equ(3.68).
From now on we take $K=4$. Combining \equ(3.68) with
$$
\align
&EX^p\big(s+t+C_6\ka(t) \big) \\ &= p \int_0
^{2\sqrt d C_1\big(s+t+C_6\ka(t)\big)} \al^{p-1}P\{
X\big(s+t+C_6\ka(t) \big)\ge \al\} d\al\\
&\le p \int_0^{2\sqrt d C_1\big(s+t+C_6\ka(t)\big)}
\al^{p-1}P\{X(s)+X'(t) + Y(t) \ge \al\}d\al\\
&\phantom{MMMMMMMMMMMMMM}+C_8[2\sqrt d C_1]^p\big(s+t+C_6\ka(t)\big)^ps^{-K-1},
\endalign
$$
we obtain
$$
\align
&EX\big(s+t+C_6\ka(t) \big)\\
&\le EX(s) + EX(t) + 2 \sqrt d C_1C_6\ka(t)
+K_{10}[s+t+ C_6\ka(t)]s^{-K-1}
\teq(3.68aa)
\endalign
$$
and
$$
\align
EX^2\big(2s+C_6\ka(s)\big) &\le E\big[X(s)+X'(s) + Y(s)\big]^2
+ K_{10}[s+ C_6\ka(s)]^2s^{-K-1}\\
&\le E\big[X(s)+X'(s) + Y(s)\big]^2 +4K_{10}s^{-K+1}.
\teq(3.69)
\endalign
$$
\equ(3.69) holds for any $s \ge s_0$, but so far \equ(3.68aa) has only
been proven under \equ(3.2cba). But there is a simple replacement
for this inequality that holds as soon
as $s_0\le s \le t$. Indeed, assume that $s_0 \le s \le t$, but
$t \log t  > C_7s^2$. We first observe that then
$$
X\big(s+t+C_6\ka(t)\big) - X(t) \le 2\sqrt dC_1[s+t+C_6\ka(t)] \le
5\sqrt d C_1 t,
$$
provided $s_0$ is large enough. Further, it follows from the simple inequality
$$
[a+b-c]^+ -[a-d]^+ \le |b| + [a-c]^+-[a-d]^+ \le |b| +[c-d]^-
\teq(3.69aa)
$$
that
$$
X\big(s+t+C_6\ka(t)\big) - X(t) \le 2\sqrt d C_1[s+C_6\ka(t)] +
[h^*(s+t+C_6\ka(t),u) - h^*(t,u)]^-.
$$
It follows that
$$
\align
EX\big(&s+t+ C_6\ka(t)\big)-  EX(t)\\
& \le  2\sqrt d C_1[s+C_6\ka(t)]+
E\Big[[h^*(s+t+C_6\ka(t),u) - h^*(t,u)]^-
\land 5\sqrt d C_1 t\Big]\\
&\le 2\sqrt d C_1[s+C_6\ka(t)]+\int_0^{5\sqrt d C_1 t}
 P\{h^*\big(s+t+C_6\ka(t),u\big)
\le h^*(t,u) -\al\} d\al\\
& \le 2\sqrt d C_1[s+C_6\ka(t)]+ 5K_3 \sqrt d C_1 t t^{-3}
+ 8d \int_0^\infty
\exp\Big[ -\frac {K_2\al^2}{s+C_6\ka(t) +\al}\Big]
d\al\\
& \text{(by \equ(3.67) with $K$ taken as 3)}\\
&\le 2\sqrt d C_1[s+C_6\ka(t)]+K_{11}[s+C_6\ka(t)]^{1/2}  \le K_{12}\ka(t).
\teq(3.69bb)
\endalign
$$
Together with \equ(3.68aa) this shows
$$
EX\big(s+t+C_6\ka(t) \big) \le EX(s) + EX(t) + K_{13}\ka(s+t)
\teq(3.70)
$$
for all $s_0 \le s \le t$.
\comment
For $s_0 \le t \le s$ we  have similarly
$$
\align
EX\big(s+t +C_6\ka(t)\big) &\le EX\big(s+t+C_6\ka(s)\big)
+2 \sqrt dC_1C_6|\ka(t)-\ka(s)|\\
&+ \int_0^\infty P\{h^*\big(s+t+C_6\ka(t)\big) - h^*\big(s+t+
C_6\ka(s)\big) \le -\al\}d\al\\
&\le EX(s) + EX(t) + K_{12}\ka(s+t),
\endalign
$$
so that \equ(3.70) holds for all $s \land t \ge s_0$ (after an
adjustment of $K_{12}$).
\endcomment

We shall next use a small variation on the argument of \cite{Ha}
to show that \equ(3.70) implies
$$
\la(u) := \lim_{t \to \infty} \frac 1t Eh^*(t,u)) \text{ exists}.
\teq(3.71)
$$
It suffices for \equ(3.71) to show that
$$
\lim_{t \to \infty} \frac 1t EX(t) = 2 \sqrt d C_1 - \la(u),
\teq(3.72)
$$
because
$$
\limt \frac 1t E\{h^*(t,u); h^*(t,u) \ge 2 \sqrt dC_1 t\} =0,
$$
by virtue of \equ(3.62). Now define for any $M \ge e$,
$$
t_0(M) = M, \; t_{k+1}(M) = 2t_k(M) + C_6\ka(t_k(M)).
$$
Note that $t_{k+1}/t_k \ge 2$, and hence $t_k(M) \ge 2^k M$, and for
large $k$
$$
1 < \frac {t_{k+1}(M)}{2t_k(M)} \le 1 +
\frac {C_6}2\Big( \frac{k\log 3}{2^kM}
\Big)^{1/2},
$$
and for some $K_{14}$, independent of $k \ge 0$,
$$
1 \le \prod_{j=0}^{k-1} \frac {t_{j+1}(M)}{2t_j(M)} = \frac {t_k(M)}{M2^k}
\le 1+ \frac{K_{14}}{M^{1/2}}.
\teq(3.73)
$$
Also, by \equ(3.70), for all $M \ge s_0+e$,
$$
EX(t_k(M)) \le 2EX(t_{k-1}(M)) +
K_{13}\ka\big(t_k(M)\big),
\; k \ge 1.
$$
Consequently,
$$
\align
&\frac {EX(t_k(M))}{t_k(M)}\\
 &\le \frac {EX(M)}M \prod_{j=1}^k
\frac {2t_{j-1}(M)}
{t_j(M)} + K_{13} \sum_{\l=0}^{k-1}  \frac
{\ka\big(t_{k-\l}(M)\big)}{t_{k-\l}(M)}
\prod_{j=k-\l+1 }^k\frac {2t_{j-1}(M)}{t_j(M)}\\
&\le \frac{EX(M)}M + K_{15} \frac{[\log M]^{1/2}}{M^{1/2}}, \quad k
\ge 0.
\endalign
$$
In particular, for given $\ep > 0$ we can choose $M \ge s_0+e$ so large, that
$$
K_{15}[\log M]^{1/2}M^{-1/2} < \ep
\text{ and }EX(M)/M \le \liminf_{s \to\infty} EX(s)/s + \ep.
$$
Then
$$
\frac {EX(t_k(M))}{t_k(M)} \le \liminf_{s \to \infty} \frac {EX(s)}s +
2\ep, \quad k \ge 0.
\teq(3.73abc)
$$

Now let $q_0 \ge s_0+M$ be large. We shall expand $q_0$ as a sum of
the form $\sum t_{k(i)}$ plus some error terms (see (3.97)) and
obtain a corresponding bound for $EX(q_0)$ in (3.96).
We define $k(1)$ as the unique integer $k$
for which $t_k \le q_0 < t_{k+1}$. We distinguish two cases.
We are in the first
case if $q_0 \ge t_{k(1)} + C_6\ka\big(t_{k(1)}\big)+ s_0 +M$.
In this case we set
$q_1 = q_0 -t_{k(1)} - C_6\ka\big(t_{k(1)}\big) < t_{k(1)+1}
-t_{k(1)} - C_6\ka\big(t_{k(1)}\big) = t_{k(1)}$. Then $s_0 +M \le q_1
< t_{k(1)}$ and
$$
EX(q_0) \le EX\big(t_{k(1)}\big) +EX(q_1) + K_{13}\ka(q_0),
$$
by virtue of \equ(3.70).
If $t_{k(1)} \le q_0 < t_{k(1)} +   C_6\ka\big(t_{k(1)}\big) + s_0+M$,
then, as in  \equ(3.69aa),\equ(3.69bb),
$$
\align
&EX(q_0) \\
&\le EX\big(t_{k(1)}\big) + 2\sqrt d C_1[q_0-t_{k(1)}]
+\int_0^{2\sqrt dC_1q_0}
P\{h^*(q_0,u) -h^*\big(t_{k(1)},u\big) \le -\al\}d\al\\
&\le EX\big(t_{k(1)}\big) + K_{16}\ka(q_0) \text{ (by \equ(3.67))}
\endalign
$$
for a suitable large constant $K_{16}$. If we are in the
first case, we repeat
the above  procedure with $q_0$ replaced by $q_1$. That is, we find
$k(2)$ such that $t_{k(2)} \le q_1 < t_{k(2)+1}$ etc. We continue to
determine $k(i)$ and $q_i$ until for the first time $q_i$ is in
the second case, i.e., $t_{k(i+1)} \le q_i < t_{k(i+1)} + C_6
\ka\big(t_{k(i+1)}\big) +s_0+M$. Suppose this first happens at the index
$i_0$. We then have
$$
\align
EX(q_0) &\le EX\big(t_{k(1)}\big) +EX(q_1) + K_{13}\ka(q_0)\le \cdots\\
&\le \sum_{i=1}^{i_0+1} EX\big(t_{k(i)}\big)
+ (K_{13} + K_{16})\sum_{i=0}^{i_0} \ka(q_i).
\teq(3.74aa)
\endalign
$$
Note that by construction, $q_i < t_{k(i)}$ for $1 \le i \le i_0$,
and consequently,
$k(i+1) < k(i)$ for $i < i_0$. Therefore the above procedure ends at a
finite $i_0$, and
$$
(K_{13} + K_{16})\sum_{i=0} ^{i_0}\ka(q_i) \le
(K_{13} + K_{16})\Big[\sum_{k:t_k \le q_0} \ka(t_k) +\ka(q_0)\Big]
\le K_{17} \big[q_0\log q_0\big]^{1/2}.
$$
In addition we have either $i_0 =0$ and $q_0 \ge t_{k(1)}$ or $i_0 \ge 1$
and
$$
q_0 = t_{k(1)} + C_6\ka\big(t_{k(1)}\big) + q_1=\cdots =
\sum_{i=1} ^{i_0}\big[ t_{k(i)}
+ C_6\ka\big(t_{k(i)}\big)\big] +q_{i_0} \ge \sum_{i=1} ^{i_0+1} t_{k(i)}.
\teq(3.74xy)
$$
Finally, we note that by definition of $i_0,\;
q_{i-1} \ge s_0+M$, and therefore $t_{k(i)} \ge M$, for $i \le i_0$.
\equ(3.74aa) and \equ(3.73abc) now show that
$$
\align
\frac {EX(q_0)}{q_0} \le \frac 1{q_0}
\sum_{i=1}^{i_0+1} t_{k(i)} &\Big[\liminf_{s \to \infty} \frac {EX(s)}s
+2\ep \Big]\\
&+K_{17} \Big[\frac{\log q_0}{q_0}\Big]^{1/2}
+I[t_{k(i_0+1)} < M] \frac {\max_{j < M} EX(j)}{q_0},
\endalign
$$
whence
$$
\limsup_{q \to \infty}\frac{EX(q)}q \le
\liminf_{s \to \infty} \frac {EX(s)}s
+3\ep.
$$
Thus the limit in \equ(3.72) exists and we can use \equ(3.72) to
define $\la(u)$.

We next turn our attention to the second moments. By definition $0 \le
X(s) \le 2\sqrt d C_1s$. With this inequality as a replacement of
(11.10) and (11.12) in \cite {Ha}, we can start from \equ(3.69)
and imitate without essential changes the
computations on p. 676 of \cite {Ha} or pp. 21, 22 of \cite{SW}
to obtain for any $M \ge$ some $s_4$
$$
\sum_{k=0}^\infty \frac {\si^2\big(X(t_k(M)\big)}{(M2^k)^2} < \infty.
$$
Since $t_k(M)/(M2^k) \ge 1$ (see \equ(3.73)) we even
have
$$
\sum_{k=0}^\infty \frac {\si^2\big(X(t_k(M)\big)}{[t_k(M)]^2} < \infty,
\teq(3.74cba)
$$
and hence for any $\ep >0$,
$$
\sum_{k=0}^\infty P\Big\{\frac 1{t_k(M)}\big|X(t_k(M))-2 \sqrt d C_1
+\la(u)\big| \ge \ep\Big\} <\infty.
$$
By \equ(3.62)  also
$$
\sum_{k=0}^\infty P\{ X\big(t_k(M)\big) \ne 2\sqrt d C_1 t_k(M) -
h^*(t_k(M),u)\} < \infty,
$$
so that for each fixed $M \ge s_4$ and $ u \in S^{d-1}$,
$$
\sum_{k=0}^\infty P\Big\{ \frac 1{t_k(M)}\big|h^*(t_k(M),u) -
\la(u)\big| \ge \ep\Big\} < \infty.
\teq(3.75)
$$

Of course \equ(3.75) implies $h^*(t_k(M),u)/t_k(M) \to \la(u)$, almost
surely.
Since $X(s) \ge 0$ by definition, $2 \sqrt dC_1 - \la(u) \ge 0$ in
\equ(3.72), and
hence $\la(u) \le 2\sqrt d C_1$, as claimed. Finally,
$\la(u) \ge C_4$ follows from Lemma 2 and the almost sure convergence of
$h^*(t_k(M),u)/t_k(M)$ to $\la(u)$.
In fact, \equ(3.10) shows
that almost surely $h^*(t_k(M),u) = h\big(t_k(M),u,-C_5\ka(t_k(M)\big)
\ge C_4t_k(M)$
for all large $k$.

Now choose a large $M_0 \ge s_4$ and for some large integer $r$ take
$M_i =  M_02^{i/r}, i=0, 1, \dots, r-1$. Further take
$M_r = t_1(M_0)    $
and note that $M_{i+1}/M_i \to 2^{1/r}$ as $M_0 \to \infty$ for fixed $r$
and $0 \le i \le r-1$. For given $\eta > 0$ we can
therefore first choose $r$ large, such that
$1 < 2^{3/r} < 1+\eta$, and then $M_0$ so large that
$$
2^{1/(2r)} \le
\frac {M_{i+1}}{M_i} \le 2^{2/r},  \; 0 \le i \le r-1.
$$
By \equ(3.73) we may further take $M_0$ so
large that
$$
2^{-1/(4r)} \frac {M}{M'} \le
\frac {t_k(M)}{t_k(M')} \le
2^{1/r} \frac {M}{M'}, \text{ for } M \ge M' \ge M_0,\, k \ge 0.
$$
Once these choices have been made we take for
$\{n_j\}_{j \ge 0}$ the collection of all distinct
$t_k(M_i), 0 \le i \le r-1, k
\ge 0$, arranged in increasing order.
Note that $i$ only runs to $r-1$ here. We claim that the collection $\{n_j\}$
in increasing order is $\{M_0, M_1, \dots, M_{r-1}, t_1(M_0), \dots,
t_1(M_{r-1}), t_2(M_0), \dots \}$. To verify this we merely need to
check that $t_k(M_0) > t_{k-1}(M_{r-1})$, since the other orderings
are obvious from the monotonicity of $t_j(\cdot)$. However, $t_k(M_0)
= t_{k-1}(t_1(M_0)) > t_{k-1}(M_{r-1})$ is also easy from $t_1(M_0)
\ge 2M_0 > M_{r-1}$. This proves our claim.

By construction we now have for all $j \ge 0$,
$$
\align
2^{1/(4r)} &\le
2^{-1/(4r)}\inf\Big\{\frac {t_k(M_{i+1})}{t_k(M_i)}:
k\ge 0,\; 0 \le i \le r-1\Big\}
\le \frac{n_{j+1}}{n_j}\\
& \le
2^{1/r}\sup\Big\{\frac {t_k(M_{i+1})}{t_k(M_i)}:
k\ge 0,\; 0 \le i \le r-1\Big\} \le 2^{3/r} \le 1+\eta.
\teq(3.72cba)
\endalign
$$
\comment
$$
\align
\frac{n_{j+1}}{n_j} &\le
2^{1/r}\sup\Big\{\frac {t_k(M_{i+1})}{t_k(M_i)}:
k\ge 0,\; 0 \le i \le r-1\Big\}\\
&\le
2^{2/r}\sup\Big\{\frac {M_{i+1}}{M_i}:
k\ge 0,\; 0 \le i \le r-1\Big\} \le 2^{3/r} \le 1+\eta.
\teq(3.72cba)
\endalign
$$

It is obvious that the preceding construction allows us to choose
$\{n_j(\eta')\}$ such that it contains $\{n_j(\eta)\}$, if $\eta' <
\eta$, as claimed at the end of the lemma.
\endcomment
The leftmost inequality here shows that $n_j$ increases exponentially
with $j$.

It is simple to see that one can choose
$\{n_j(\eta')\}$ such that it contains the $\{n_j(\eta)\}$ from some
index on, if $\eta' <
\eta$, as claimed at the end of the Corollary. In fact
if the $n_j(\eta)$ are constructed by the above method for some
$M_0,r$, then one can use the same construction for the $n_j(\eta')$
based on $M_0', r'$ with $M_0' = t_k(M_0)$ for some  $k$ and $r'$
some integer multiple of $r$.

Next, \equ(3.61) holds, because by \equ(3.75)
$$
\align
&\sum_{k=0}^\infty P\Big \{\Big|\frac 1{n_k}h^*(n_k,u) -\la(u)\Big| >
\ep\Big\}\\
&= \sum_{i=0}^{r-1}\sum_{k=0}^\infty
P\Big\{\Big|\frac1{t_k(M_i)}h^*(t_k(M_i),u)
- \la(u)\Big| > \ep\Big\} < \infty.
\teq(3.72abc)
\endalign
$$
Thus also
$$
\lim_{k \to \infty} \frac 1{n_k}h^*(n_k,u) = \la(u) \text{ a.s.}
\teq(3.72bcd)
$$

Now let $0 < \ep \le C_4/2 \le \la(u)$ and $\eta 4 \sqrt d C_1 < \ep/ 2$.
If
$$
\frac 1{n_k} h^*(n_k,u),
\; \frac 1{n_{k+1}} h^*(n_{k+1},u) \le 2\la(u) \le 4\sqrt d C_1
\teq(3.76)
$$
and
$$
h^*(n_k,u) - K_{18}\ka(n_k) \le h^*(t,u) \le h^*(n_{k+1},u)
+ K_{18}\ka (n_{k+1})
\teq(3.77)
$$
for all $n_k \le t \le n_{k+1}$,
then, for $k$ large enough and all $n_k \le t \le n_{k+1}$
$$
\align
\frac 1{n_k} h^*(n_k,u) -\ep &\le
\frac 1{n_{k+1}} h^*(n_k,u) - \frac {K_{18} \ka(n_k)}{n_{k+1}}
\le \frac 1t h^*(t,u) \\
&\le \frac 1{n_k} h^*(n_{k+1},u) +
\frac {K_{18} \ka(n_{k+1})}{n_k} \le \frac 1{n_{k+1}} h^*(n_{k+1},u) + \ep.
\teq(3.78)
\endalign
$$
By choosing $K_{18}$ large enough, and applying \equ(3.67),
\equ(3.67bb), we can make
$$
\align
&\sum_{k=0}^\infty P\{\equ(3.78) \text{ fails for some }t \in
[n_k, n_{k+1}]\} \\
&\le
\sum_{k=0}^\infty P\{\equ(3.76) \text{ fails}\}
+\sum_{k=0}^\infty P\{\equ(3.77) \text{ fails}\}\\
&\le 2 \sum_{k=0}^\infty P\Big \{\Big|\frac 1{n_k}h^*(n_k,u) -\la(u)\Big| >
\ep\Big\} + 2\sum_{k=0}^\infty K_3[n_k]^{-4}\\
&\phantom{MMMM}+ \sum_{k=0}^\infty
\Big[8d +K_4\big(n_{k+1}\big)^d\Big] \exp[-K_{19} K_{18}^2
\eta^{-1} \log n_k]\\
& < \infty.
\teq(3.79)
\endalign
$$
Since $\ep >0$ can be taken arbitrarily small,
this, together with \equ(3.72bcd),
proves the almost sure convergence in \equ(3.60).
The $L^p$ convergence along all reals
in \equ(3.60) follows from the almost sure convergence and the tail
estimate \equ(3.62).
\qedsymbol
\enddemo

\subhead
4. From half-space to full-space processes
\endsubhead
\numsec=4
\numfor=1

The goal for this section is to prove that the $B$-particles in the
full-space process do not spread faster than in the appropriate
half-space process (see Corollary 8 for a precise statement). The first
lemma establishes that for every $u \in S^{d-1}$ there are
deterministic vectors $V_k$ such that for all
$\eta > 0$ there is with a probability close to 1 a $B$-particle in
$\cP^h\big(u, -C_5\ka((1+\eta)n_k)\big)$ ``near'' $V_k$ at time $n_k$, for all
large $k$. Here $n_k$ is the $n_k(\eta)$ of Corollary 5 and $\langle
V_k,u \rangle$ has to grow essentially like $h^*(n_k,u) \sim n_k
\la(u)$ (see (4.1)).
Apart from this growth condition the behavior of $V_k$ as a function
of $k,u$ is unimportant
for us. The only important aspect is that it is non-random, so that we
can find, with high probability, a $B$-particle in a non-random
location at which $h^*(n_k,u)$ is (almost) achieved. This will be used
in the second lemma to concatenate $\cP^h\big(u, -C_5\ka(n_k)\big)$
with another process which runs from time $(1+\eta)n_k$ to
$(1+\eta)n_k+r_k$ with $r_k$
also of order $n_k$.
By starting the second process at $\big(V_k,(1+\eta)n_k\big)$
we will be able to assure
that a $B$-particle at time $(1+\eta)n_k+r_k$ in
the second process is also a $B$-particle in $\cP^h\big(u,
-C_5\ka((1+\eta)n_k+r_k)\big)$.
\proclaim{Lemma 6} Let $u \in S^{d-1}$ be fixed, and let $n_k=
n_k(\eta)$ be as in Corollary 5. Then, for all $0 < \eta < C_4/(8C_1)$
there exists a deterministic sequence of vectors $\{V_k\} =
 \{V_k(\eta,u)\}$ such that
$$
\langle V_k(\eta, u), u \rangle = n_k(\eta)\la(u),
\teq(9.2)
$$
and such that
$$
\align
&\sum_{k=0}^\infty P\{\text{in $\cP^h\big(u, -C_5\ka(n_k)\big)$
there is at time $(1+\eta)n_k$ either no particle}\\
&\phantom{MM}\text{at all in
$V_k+\cC(C_2 n_k \eta/4)$ or there is an $A$-particle in
$V_k+ \cC(C_2n_k\eta/4)$}\}\\
& <\infty.
\teq(9.13)
\endalign
$$
\endproclaim
\demo{Proof} Fix $u \in S^{d-1}$ and $\ep > 0$.
Let $\si$ be a time which is so large that $\si \ge s_0$ (with $s_0$
is as in Proposition 3) and such that
$$
\Big|\frac 1\si Eh^*(\si,u)- \la(u)\Big| \le \frac 14 \ep
\teq(9.3)
$$
(see \equ(3.13aa) for $h^*$). Define the further times
$$
\si_1 = \si, \si_{j+1} = \si + \si_j +C_6\ka(\si_j), \; j\ge 1.
$$
Now apply \equ(3.2) with the following choices:
$s=\si, t = \si_j, \ga = \ga_s = E m^*(\si,u)$ (see \equ(3.13cc) for $m^*$)
and $\ga_t = j Em^*(\si,u)$. This yields
$$
\align
&P\{\cG\big(\al, \be, (j+1) E m^*(\si,u),
\cP^h\big(u, - C_5\ka(\si_{j+1})\big), \si_{j+1}\big)\\
&\ge \int_{0 \le h < \infty} \int_{m \in \Bbb R^d}
P\{h^*(\si,u) \in dh, m^*(\si,u) \in \ga +dm\}\\
&\phantom{MMMMMM}P\{\cG\big(\al-h, \be-d, j\ga -m, \cP^h\big(u,
-C_5\ka(\si_j) \big), \si_j\big)\} -C_8 \si^{-K-1},
\endalign
$$
provided \equ(3.2cba) holds, that is, provided
$(\si_j+1)\log (\si_j +1) \le C_7 \si^2$. We start with $j=
r-1$, then use the case $j=r-2$ with $\al, \be$ replaced by $\al-h$
and $\be -d$, respectively,  etc., all the way down to $j=1$.
With $(h^*_j, m^*_j), \; j \ge 1$, i.i.d. copies of $\big(h^*(\si,u),
m^*(\si,u)\big)$ we obtain
$$
\align
&P\{\cG\big(\al, \be, r E m^*(\si,u),
\cP^h\big(u, - C_5\ka(\si_r)\big), \si_r\big)\}\\
&\ge \int\limits \Sb h_j \ge 0, \\1 \le j \le r-1 \;\endSb \int \limits
\Sb m_j \in
\Bbb R^d\\ 1 \le j \le r-1 \endSb
\prod_{j=1}^{r-1} P\{h^*(\si,u) \in dh_j, m^*(\si,u) \in \ga +
dm_j\}\\
&\phantom{MMMM} \times
P\Big\{\cG\big(\al - \sum_{j=1}^{r-1}h_j, \be -(r-1)d,
\ga - \sum_{j=1}^{r-1}m_j, P^h\big(u, -C_5\ka(\si)\big), \si\big)\Big\}\\
&\phantom{MMMMMMMMMMMMMMMMMMMMMMM} - (r-1)C_8\si^{-K-1}\\
&=P\Big\{\text{in $\cP^h\big(u,-C_5\ka(\si)\big)$ there is at time
$\si$ a
$B$-particle at some $x$ with }\\
&\phantom{MMMM}\langle x,u \rangle +\sum_{j=1}^{r-1}
h^*_j \ge \al
\text{ and }\big\|x^\perp + \sum_{j=1}^{r-1}m^*_j- \ga\big\| \le \be
-(r-1)d\Big\}\\
&\phantom{MMMMMMMMMMMMMMMMMMMMMMM} -(r-1)C_8\si^{-K-1}\\
&\ge P\Big\{ \sum_{j=1}^r h^*_j \ge \al, \big\|\sum_{j=1}^r (m^*_j -
\ga)\big\| \le \be -(r-1)d\Big\} - (r-1)C_8\si^{-K-1},
\teq(9.6)
\endalign
$$
provided
$$
(\si_r+1) \log (\si_r+1) \le C_7\si^2.
\teq(9.7)
$$

It is easy to see by induction that each $\si_j$ is a continuous,
increasing function of $\si$ on $[0,\infty)$. We further see by induction that
$\si_k \ge k\si$ and $\si_j$ increases with $j$. Finally, we can for
any fixed $\si \ge 1$ find a $K_1= K_1(\si)$ such that
$$
\si K_1 2k(\log k +1) \ge C_6\ka\big(\si K_1k^2(\log k+1)\big),\quad k \ge 1,
$$
and $\si_1 \le \si K_1 \log 2$. One more induction argument then shows
that for all $k \ge 1$, $\si_k \le \si K_1 k^2 (\log k+1)$.
Now let $s\ge s_0$ be large and take $r = \lfloor s^{1/3}\rfloor$.
The preceding argument shows that we can fix $\si$ such that $\si_r
= s$.  Thus for $j-1 \le r$ we have $\si_{j-1} \le \si_r = s$
and $\si_j \le \si + \si_{j-1} + C_6\ka(s)$. Consequently, $s = \si_r
\le r\si +rC_6 \ka(s)= r\si + \lfloor s^{1/3}\rfloor C_6\ka(s)= r\si +
o(s)$, and necessarily $\si \sim s/r \sim
s^{2/3}$ for large $s$. \equ(9.7) is therefore automatically
satisfied. If we further take
$$
\al = r \si [\la(u) - \frac 12 \ep],
$$
then, by \equ(9.3) and the fact that Variance $(h^*_j) \le K_2 \si^2$ (by
\equ(3.62)),
$$
P\big\{\sum_{j=1}^r h^*_j \le \al \big\}
\le P\big \{\sum_{j=1}^r\big(h^*_j - Eh^*_j\big) \le -r\si \ep/4\big\}
\le \frac{K_3}{r \ep^2}.
\teq(9.8)
$$
Further, fix $s_5$ so large that $2s\ep \ge r\si \ep \ge
(1/2)s \ep \ge 2 rd \sim
2s^{1/3}d$ for $s \ge s_5$. Then we have similarly to \equ(9.8), for
$s \ge s_5$ and $\be = s \ep$
$$
P\big\{\big \|\sum_{j=1}^r \big(m^*_j - \ga\big)\big\| > \be -(r-1)d\big\}
\le P\big\{\big \|\sum_{j=1}^r \big(m^*_j - \ga\big)\big\| > r \si \ep/4
\big\}\le \frac {K_4}{r\ep^2}.
\teq(9.9)
$$
The last two inequalities provide us with a lower bound for the
right hand side of \equ(9.6). We conclude that for $s \ge s_5$
$$
\align
P\{\cG\big(\al, \be,\; &r E m^*(\si,u), \cP^h\big(u, -
C_5\ka(\si_r)\big), \si_r\big)\}\\
&\ge 1- \frac{(K_3+K_4)}{r\ep^2} - (r-1)C_8\si^{-K-1}
\ge 1- \frac {K_5}{s^{1/3}\ep^2}
\teq(9.10)
\endalign
$$
(use any $K \ge 1$). Let $n_j(\eta)$ be as in Corollary 5, and
take $s = n_k = n_k(\eta)$. In agreement with our previous choice for
$r, \si$ we then take $r = \lfloor n^{1/3}_k(\eta) \rfloor$ and $ \si$ such
that $\si_r = n_k(\eta)$. Then, by going over to the complementary
event in \equ(9.10), we find for any $\eta > 0$, that
$$
\align
&\sum_{k=0}^\infty
P\big\{\text{in $\cP^h\big(u, -C_5\ka(n_k)\big)$
there is at time $n_k$ no $B$-particle }\\
&\phantom{MMMMMMMMMMMMMM}
\text{in } \Ga\big(n_k[\la(u) - \frac 14 \ep],n_k \ep, rE
m^*(\si,u)\big)
\big\}\\
&\le \sum_{k=0}^\infty \frac {K_5}{n_k^{1/3}\ep^2}< \infty
\teq(9.11)
\endalign
$$
(recall that the $n_j$ grow exponentially).
But \equ(3.61) says in particular that
$$
\align
&\sum_{k=0}^\infty P\big\{\text{in $\cP^h\big(u, -C_5\ka(n_k)\big)$
there is at time $n_k$ a $B$-particle  }\\
&\phantom{MMMMMMMMMMMMMM}
\text{in }\Ga\big(n_k[\la(u) + \frac 14 \ep],n_k \ep, rE
m^*(\si,u)\big)
\big\}\\
&< \infty.
\teq(9.12)
\endalign
$$
We now take
$$
V_k = V_k(\eta, u) = n_k(\eta)\la(u) u + rEm^*(\si,u).
\teq(9.12a)
$$
Since $m^*$ is orthogonal to $u$ (by definition), this
choice of $V_k$ satisfies \equ(9.2). Moreover,
\equ(9.11) and \equ(9.12) together give
$$
\align
&\sum_{k=0}^\infty
P\{\text{in $\cP^h\big(u, -C_5\ka(n_k)\big)$
there is at time $n_k$
no $B$-particle}\\
&\phantom{iMMMMMMMMMMMMMMMMM} \text{at any site }
x \in V_k + \cC(2 \sqrt d n_k \ep)\}\\
&\le \sum_{k=0}^\infty
P\{\text{in $\cP^h\big(u, -C_5\ka(n_k)\big)$
there is at time $n_k$ no $B$-particle at any site $x$ }\\
&\phantom{iMM}\text{with }\langle x, u \rangle \in
\Big[n_k[\la(u) - \frac{\ep}4],
n_k[\la(u) +\frac \ep 4]\Big], \|x^\perp - rEm^*(\si,u)\|\le n_k \ep\}\\
& < \infty.
\teq(9.14)
\endalign
$$

The convergence of the sums in \equ(9.14)
shows that almost surely, for all large
$n_k(\eta)$, there is in $\cP^h\big(u,-C_5\ka(n_k)\big)$
a $B$-particle in  $V_k + \cC(2 \sqrt d  n_k \ep)$
at time $n_k(\eta)$. We claim that this implies that if we take $\ep
= C_2\eta/(16 d)$, then, in $\cP^h\big(u,-C_5\ka(n_k)\big)$
at time $(1+\eta)n_k$, all occupied sites
in $V_k + \cC(C_2 n_k \eta/4)$ are occupied by $B$-particles (and there are
such occupied sites). More precisely, we claim that \equ(9.13) holds.
To see this we shall apply Lemma 4 with the following choices: $s =
n_k, \wt s = (1+\eta)n_k, t = (1+\eta)n_k$
and finally $y(n_k)$ is the location of any $B$-particle in
$\cP^h\big(u, -C_5\ka(n_k)\big)$ at time $n_k$ in the set $V_k+ \cC\big(C_2
n_k \eta/(8\sqrt d)\big)$, if such a $B$-particle exists.
If several such $B$-particles
exist we pick the location of one of them according to some
deterministic rule chosen in advance. On the event that no such
$B$-particle exists we cannot apply Lemma 4, but this does not
cause any problems, because \equ(9.14) already tells us that
$$
\sum_{k=0}^\infty P\{\text{no choice for $y(n_k)$ exists}\} < \infty.
\teq(9.15)
$$
If $y(n_k)$ exists, then there is automatically a particle in
$\cP^h\big(u, -C_5\ka(n_k)\big)$
at time $n_k$ at $y(n_k) \in V_k+\cC\big(C_2 n_k\eta/(8\sqrt d)\big)$.
If this particle does not move a distance $> C_2n_k \eta/8$
during $[n_k, (1+\eta)n_k]$, then it is in $y(n_k) +\cC(C_2n_k\eta/8)
\subset V_k + \cC(C_2n_k \eta/4)$ at time $(1+\eta)n_k$.
We recall further that
all particles in $\cP^h\big(u, -C_5 \ka(n_k)\big)$
are also particles in $\cP^f$.
It follows that the $k$-th summand
in \equ(9.13) is bounded by the $k$-th summand in \equ(9.15) plus
$$
P\{\|S_{n_k \eta}\| > C_2n_k \eta/8\}
+ P\{\cB^h\big(y(n_k) ,n_k;u,
-C_5\ka(n_k)\big) \cap \cK\big(y(n_k)\big)\}
\teq(9.16)
$$
(see \equ(6.1) for $\cK(y)$).
The first probability in \equ(9.16) is at most $K_6\exp[-K_7n_k \eta]$
by (2.42) in \cite {KSa}.
The last probability in
\equ(9.16) is by Lemma 4 at most
$$
P\{y(n_k) \notin \cC(2C_1n_k)\} +P\{\langle y(n_k),u \rangle
< \frac 12 C_4n_k\} + n_k^{-K-1}.
\teq(9.17)
$$
The first probability in \equ(9.17) is $O\big(n_k^{-K-1}\big )$ by the
estimates used for \equ(3.50). The second probability in \equ(9.17) is
zero, because, by construction, $y(n_k) \in V_k +\cC\big(C_2n_k
\eta/(8\sqrt d)\big)$, so
that
$$
\align
\langle y(n_k),u \rangle &\ge \langle V_k, u \rangle -C_2n_k \eta /8 \\
&= n_k\la(u) - C_2n_k \eta/ 8 \ge n_k(C_4 - C_2 \eta/8)
\text{ (see Corollary 5) }
\ge \frac 12 C_4n_k.
\endalign
$$
It follows that the sum of \equ(9.17) over $k$ is also finite, and
this proves \equ(9.13).
\qedsymbol
\enddemo

We can now show how to concatenate two processes as outlined before
the last lemma.

\proclaim{Lemma 7}
Define
$$
\align
H(t,u)&= h(t,u, - \infty)\\
&= \max \{\langle x,u \rangle:
x \text{ is occupied by a $B$-particle in }\cP^f \text{ at time $t$}\}.
\teq(9.19)
\endalign
$$
Assume that for some fixed $u \in S^{d-1}$ and $\mu \ge 0$
$$
P\{\limsup_{t \to \infty} \frac 1t H(t,u) \ge \mu\} > 0.
\teq(9.20)
$$
Then
$$
\la(u) \ge \mu.
\teq(4.53aa)
$$
\endproclaim
\demo{Proof} We divide the proof into 4 steps. Without loss of
generality we assume $\mu > 0$.
\newline
{\bf Step 1.} For each small $\eta >0$ we choose
$$
K_1 > 2 \sqrt d C_1 \ge  \la(u) \text{ and } K_2 = \frac 1{4C_1 \sqrt d K_1}.
\teq(4.53bb)
$$
We then define
$$
m_k = m_k(\eta)= K_2n_k(\eta),
\teq(9.20abc)
$$
where $n_k = n_k(\eta)$ is as in Corollary 5. We take $\eta_0 =
\eta_0(\ep) > 0$ so small that
$$
1+\eta_0 \le \frac{\mu-\ep/2}{\mu - 3\ep/4}.
$$
Note that these definitions imply that for $\eta \le \eta_0$,
$$
\frac {m_{k+1}} {m_k} = \frac {n_{k+1}} {n_k} \le 1+\eta
\le \frac{\mu-\ep/2}{\mu - 3\ep/4}.
\teq(9.20aa)
$$
Further, for small $\ep >0$, define the events
$$
\align
&\cL_k(\eta, \mu-\ep)\\
&= \big\{\text{in $\cP^f$ there
is a $B$-particle in the half-space }
\text{$\cS\big(u,m_k(\mu - \ep)\big)$ at time }m_k\big\}\\
&=\{H(m_k,u) \ge m_k(\mu - \ep)\}.
\endalign
$$
In this step we shall show that for fixed $\ep >0$ and all
$0 < \eta \le \eta_0(\ep)$,
$$
\sum_{k=0}^\infty P\{\cL_k(\eta, \mu -\ep)\} = \infty.
\teq(9.21)
$$
To prove this we shall show that
$$
P\{\cL_k(\eta, \mu-\ep) \text{ occurs for infinitely many } k\} > 0.
\teq(9.22)
$$
\equ(9.21) then follows from the Borel-Cantelli lemma.
Now, \equ(9.20) says that for every $\ep > 0$
$$
P\{\text{for infinitely many $k,  H(t,u)> (\mu-\ep/2)t$ for some }t
\in [m_k,m_{k+1}]\} > 0.
\teq(9.23)
$$
However, by \equ(3.67bb) with $h^*$ replaced by $H$ (this amounts to
taking $C_5 = \infty$, which does not influence the estimate \equ(3.67bb))
and with $\al = (\ep /4)m_{k+1} \le (\mu - \ep/2)m_k -(\mu-\ep)m_{k+1}$
(see \equ(9.20aa))
$$
\align
&P\{H(t,u)> (\mu -\ep/2)t \text{ for some $t
\in [m_k,m_{k+1}]$ but } H(m_{k+1},u) \le (\mu - \ep)m_{k+1}\}\\
&\le P\{\sup_{r \in [m_k,m_{k+1}]}\big[ H(r,u) - H(m_{k+1},u)\big] \ge
(\mu-\ep/2)m_k -(\mu -\ep)m_{k+1}\}\\
&\le K_3(\ep, \eta)[m_{k+1}]^{-K}.
\endalign
$$
In particular, by Borel-Cantelli, the event in the left hand side here
occurs almost surely only finitely often. Together with \equ(9.23)
this shows that
$$
P\{\text{for infinitely many }k,  H(m_{k+1},u) \ge  (\mu-\ep)m_{k+1}\} > 0.
$$
This is the required \equ(9.22).
\medskip \noindent
{\bf Step 2.}
The remaining steps are based on \equ(9.21) only; \equ(9.20) itself is not
needed.
With $V_k= V_k(\eta,u)$ as in \equ(9.12a) we
define an auxiliary process $\cQ_k = \cQ_k(\eta,u)$
which is more or less the full-space
process started at the deterministic space time point $(V_k,(1+\eta)n_k)$. The
only difference is that $\cQ_k$ only uses the particles which are {\it
at time} 0 in
the ``slab''
$$
\{x:-n_k/K_1 \le \langle x,u \rangle - n_k\la(u) < K_1n_k\},
\teq(9.18)
$$
with $K_1$ given by \equ(4.53bb).
Thus $\cQ_k$ is defined only from time $(1+\eta)n_k$ on. At time
$(1+\eta)n_k$ it has at any $x$ only the particles which started at
time 0 in the set \equ(9.18). If no such particles exist, then
there never are any particles in the process $\cQ_k$. Otherwise,
let $z_k$ be the nearest site
to $V_k$ which is occupied at time $(1+\eta)n_k$ by some particle,
which at time 0 was in \equ(9.18). The types of all particles in
$\cQ_k$ at time $(1+\eta)n_k$ are
reset to type $A$, except for the particles
at $z_k$, which are reset to type $B$. From time $(1+\eta)n_k$ the
process then develops by our standard rules.
Even though the process $\cQ_k$
is defined for all times in $[(1+\eta)n_k,\infty)$ we are only interested in
what happens during $[(1+\eta)n_k, (1+\eta)n_k+m_k]$.
Specifically, we define the events
$$
\align
&\cM_k =\cM_k(\eta, \mu-\ep) = \{\text{in $\cQ_k$ there is a $B$-particle
in the half-space}\\
&\phantom{MMM}\cS\big(u,n_k\la(u)  + m_k(\mu-\ep) \big)
\text{ at time } (1+\eta)n_k +m_k\}.
\teq(4.57)
\endalign
$$

In this step we shall prove that
$$
\sum_{k = 0}^\infty P\{\cM_k\} = \infty.
\teq(4.57aa)
$$
To this end let us shift the event $\cL_k$ by $(1+\eta)n_k$ in time and by
$V_k$ in space. Then $\cL_k$ goes over into the  event
$$
\align
\cL'_k := \{&\text{in the full-space process started at $(V_k,
(1+\eta)n_k)$ there}\\
&\text{is a $B$-particle in the half-space
$\cS\big(u,  n_k\la(u)+m_k(\mu-\ep)\big)$}\\
&\text{at time $(1+\eta)n_k+m_k$}\}.
\endalign
$$
$\cL'_k \setminus \cM_k$ can occur only if
one of the following two events occurs:
$$
\align
&\{\text{at time $(1+\eta)n_k$, some particle at
the nearest occupied site to }\\
&\text{$V_k$ in the
full-space process  started at time 0 outside
 the set \equ(9.18)}\},
\teq(9.24)
\endalign
$$
or
$$
\align
\{&\text{in the full-space process started at
$\big(V_k,(1+\eta)n_k\big)$
there is a }\\
&\text{particle which starts at time 0
outside the set \equ(9.18) and}\\
&\text{which coincides with a
$B$-particle during $[(1+\eta)n_k,(1+\eta)n_k+ m_k]$}\}
\teq(9.25)
\endalign
$$
(compare the argument for \equ(3.38)).
It follows that
$$
P\{\cL'_k\setminus M_k\}
\le P\{\text{\equ(9.24) or \equ(9.25) occurs}\}.
$$
But
$$
\align
&P\{\text{\equ(9.24) occurs}\} \\
&\le P\{\text{nearest occupied site to $V_k$ in
$\cP^f$ at time $(1+\eta)n_k$ has}\\
&\phantom{MMMMMMM}\text{distance more than $K_4 \log k$ from $V_k$}\}\\
&+P\{\text{some particle which starts at time 0 outside the set
\equ(9.18)}\\
&\phantom{MMMMMMM}\text{is in $V_k + \cC(K_4\log k)$ at
time $(1+\eta)n_k$}\}.
\teq(9.25a)
\endalign
$$
Also,
$$
\align
&P\{\text{\equ(9.25) occurs}\}\\
&\le
P\{\text{in the full-space process started at $(V_k,(1+\eta)n_k)$
 there are $B$-particles}\\
&\phantom{MMM}\text{outside
$V_k+\cC(2C_1m_k)$ at some time during }[(1+\eta)n_k,(1+\eta)n_k+m_k]\}\\
&+ P\{\text{some
  particle which starts at time 0 outside the set  \equ(9.18)}\\
&\phantom{MMM}\text{visits $V_k+ \cC(2C_1m_k)$ during
 }[0,(1+\eta)n_k +m_k]\}.
\teq(4.58)
\endalign
$$

The first probability in the right hand side of
\equ(9.25a) ) can be made  $O\big(k^{-K}\big)$ for any given $K$, by
choosing $K_4$ large (compare \equ(3.13)).
The second probability in the right hand side of \equ(9.25a) is
for large $k$ no
more than the second probability in the right hand
side of \equ(4.58). To estimate the latter,
we merely point out that a particle which starts at some $z$ outside the set
\equ(9.18) and visits $V_k+ \cC(2C_1m_k)$ during $[0,(1+\eta)n_k+ m_k]$
has to move over a distance of at least
$$
\align
\|z-&V_k\| -2C_1m_k \ge d^{-1/2}|\langle (z-V_k),u \rangle | - 2C_1m_k\\
&=d^{-1/2}|\langle z, u \rangle - n_k\la(u)| - 2C_1m_k
\ge  n_k/(\sqrt d K_1) - 2C_1m_k \ge n_k/(2\sqrt d K_1),
\endalign
$$
by virtue of our choice of $m_k$.
We leave it to the reader to use this to check that the last
probability in \equ(4.58) is $O\big([n_k]^{-K}\big)$ (see also the
estimates in \equ(3.18abc) and \equ(3.39abc) or \equ(6.4)). Finally, the first
probability in the right hand side of \equ(4.58) equals
$$
P\{\text{in $\cP^f$ there are $B$-particles outside $\cC(2C_1m_k)$
during $[0,m_k]$}\},
$$
and this is $O\big([m_k]^{-K-2}\big)$, as in \equ(3.18) and the lines
following it.
It follows from these estimates that  $\sum_kP\{\cL'_k \setminus
\cM_k\} < \infty$.
In view of \equ(9.21) and the fact that $P\{\cL'_k\} = P\{\cL_k\}$,
this implies \equ(4.57aa).
\medskip
\noindent
{\bf Step 3.} In this step we show that
$$
P\{\cM_k \text{ occurs for infinitely many }k \} = 1.
\teq(9.28)
$$
This is an easy application of Borel-Cantelli, because $\cM_k$ and
$\cM_\l$ depend on particles which start at disjoint sets of sites
(and are therefore independent) as soon as
the set \equ(9.18) and the corresponding set with $k$ replaced by $\l$
are disjoint. If $\l > k$, this is the case if $n_k(\la(u)+K_1) <
n_\l(\la(u) -1/K_1)$ and similarly if $k > \l$. In particular, there
is some integer $K_5 = K_5(\eta)$ such that $\cM_k$ and $\cM_\l$ are
independent as soon as $|k-\l| \ge K_5$. Moreover, by \equ(4.57aa),
there is some integer $j \in [0, K_5-1]$ such that
$$
\sum_{k \equiv j \pmod {K_5}} P\{\cM_k\} = \infty.
$$
Thus \equ(9.28) is true.

\medskip
\noindent
{\bf Step 4.} We now complete the proof of the lemma by showing that,
almost surely, for all large $k$ for which $\cM_k$ occurs, also
$$
\align
&\big\{\text{in $\cP^h\big(u, -C_5((1+\eta)n_k+m_k)\big)$
there
is a $B$-particle in the half-space}\\
&\phantom{MMMMMMMMM}
\text{$\cS\big(u,n_k\la(u)+ m_k(\mu - \ep)\big)$
at time }(1+\eta)n_k+m_k\big\}\\
&=\{h^*\big((1+\eta)n_k+m_k,u\big) \ge n_k\la(u)+ m_k(\mu - \ep)\}
\teq(9.29)
\endalign
$$
occurs.
This will indeed complete the proof, since we already know from
Corollary 5 that $\big((1+\eta)n_k+m_k\big)^{-1}
h^*\big((1+\eta)n_k+m_k,u\big)\to
\la(u)$.
Thus \equ(9.28) and \equ(9.29)
will imply, for all $\ep > 0, 0 < \eta < \eta_0(\ep)$,
$$
\align
\la(u) &\ge \liminf_{k \to \infty} \Big[\frac {n_k}{(1+\eta)n_k+m_k} \la(u)
+ \frac {m_k}{(1+\eta)n_k+m_k} (\mu-\ep)\Big] \\
&= \frac 1{1+\eta +K_2}\la(u) + \frac{K_2}{1+\eta+K_2}(\mu-\ep),
\endalign
$$
and hence
$$
\la(u) \ge \frac {K_2}{\eta+K_2} (\mu -\ep).
\teq(9.29cba)
$$

Now to prove \equ(9.29), we write, as in the lines following
\equ(9.18), $z_k$ for the nearest site to $V_k$ at time $(1+\eta)n_k$
which is occupied by a particle which started at time 0 in \equ(9.18).
We already proved that, almost surely, \equ(9.24) occurs only finitely often.
Thus, except for finitely many $k,\; z_k$ actually equals the nearest
occpied site to $V_k$ at time $(1+\eta)n_k$ in $\cP^f$. Since the set
\equ(9.18) is contained in $\cS(u, 0) \subset \cS\big(u, -C_5\ka(n_k)\big)$,
$z_k$ is also the nearest
occupied site to $V_k$ at time $(1+\eta)n_k$
in $\cP^h\big(u, -C_5\ka(n_k)\big)$. By virtue of Lemma 6,
we further know that, a.s. for all large $k$, $z_k$ is occupied by
$B$-particles at time $(1+\eta)n_k$ in $\cP^h\big(u, -C_5(n_k)\big)$ for
all large $k$. By using the monotonicity property of Lemma C
we conclude that, almost surely, for all large $k$ all the $B$-particles
in $\cQ_k$ at time $(1+\eta)n_k+m_k$ are also $B$-particles in
$\cP^h\big(u, -C_5\ka((1+\eta)n_k+m_k)\big)$. In particular,
$$
h^*\big((1+\eta)n_k+m_k,u) \ge n_k\la(u)+m_k(\mu-\ep)
$$
for all large $k$ for which $\cM_k$ occurs. This is the required
\equ(9.29).
\qedsymbol
\enddemo

\proclaim{Corollary 8} For every  unit vector $u$
$$
\lim_{t \to \infty} \frac 1t H(t,u) = \la(u)
\text{ almost surely and in $L^p$ for all }p >0.
\teq(9.30)
$$
{\rm (}$t$ runs through the reals here{\rm)}.
Moreover, for $n_k= n_k(\eta)$ as
in Corollary 5, we have for any $\de > 0$ and $\eta >0$,
$$
\sum_{k=0}^\infty
P\big\{\big|\frac 1{n_k}H(n_k,u) - \la(u)\big| > \de \big\} < \infty.
\teq(9.31)
$$
\endproclaim
\demo{Proof} By the monotonicity property of Lemma C
$$
H(t,u) \ge h^*(t,u) \text{ on the event } \{\|x_0\| \le
C_5\ka(t)/\sqrt d\}
\teq(9.32)
$$
(see also the lines after \equ(3.39)). Thus, by
the estimate \equ(3.13)
$$
\liminf_{t \to \infty} \frac 1tH(t,u)
\ge \limt \frac 1t h^*(t,u) = \la(u)
$$
(see Corollary 5). In the other direction, we have from Lemma 7 that
$$
P\{\limsup_{t \to \infty} \frac 1t H(t,u) \ge \mu\} = 0 \text{ for all
} \mu > \la(u).
$$
This proves the almost sure convergence in \equ(9.30). The $L^p$
convergence follows from the almost sure convergence and the tail
estimate
$$
P\{H(s,u) \ge \al\} \le \exp[-K_1\al] \text{ for }\al \ge 2\sqrt
dC_1s, \quad s \ge s_3,
\teq(9.32ab)
$$
which can be proven in the same way as \equ(3.62) (or we can take $C_5
= \infty$ in \equ(3.62)).

As for \equ(9.31), we have by \equ(9.32), \equ(3.61) and an estimate
like \equ(3.13) that
$$
\sum_{k=0}^\infty
P\big\{\frac 1{n_k}H(n_k,u) < \la(u)-\de \big\} < \infty.
\teq(9.33)
$$
For the other direction, we begin with an indirect argument. Assume, to derive
a contradiction, that for some $\de > 0$ and $0 < \eta \le C_4/(8C_1)$
$$
\sum_{k=0}^\infty
P\big\{\frac 1{m_k}H(m_k,u) > \la(u)+\de/2 \big\} = \infty,
$$
with $m_k=m_k(\eta)$ as in \equ(9.20abc).
This is just \equ(9.21) with $\mu-\ep$ replaced by $\la(u)+\de/2$. By
steps 2-4 of the proof of Lemma 7 we then have that \equ(9.29cba), again
with $\mu-\ep$ replaced by $\la(u)+\de/2$, holds.
 This is impossible
for $\eta < K_2\de/(2\la(u))$. Thus for all $\de > 0,
0 < \eta < C_4/(8C_1) \wedge K_2\de/(2\la(u)),$
it is the case that
$$
\sum_{k=0}^\infty
P\big\{\frac 1{m_k}H(m_k,u) > \la(u)+\de/2 \big\} < \infty.
\teq(9.34)
$$
Finally, for given $k$, let $\l = \l(k)$ be determined by $m_\l < n_k \le
m_{\l+1}$. We now use that
$$
\align
P\{H(n_k,u) &> n_k\big(\la(u) + \de\big)\} \le
P\{H(m_{\l+1},u) > m_{\l+1}\big(\la(u)
+ \de/2\big)\} \\
&+ P\{H(m_{\l+1},u) - H(n_k,u) \le m_{\l+1}\big(\la(u)+\de/2\big) -
n_k\big(\la(u)+\de\big)\}.
\teq(9.35)
\endalign
$$
But, by \equ(3.67) (with $C_5$ taken to be infinity) we have
$$
P\{\inf_{r \le t} H(s+r,u) - H(s,u) \le -\al\}
\le K_3s^{-K}+
8d \exp\Big[ -\frac {K_2\al^2}{t+\al}\Big],
\; \al \ge 0.
\teq(9.35ab)
$$
Moreover, $m_{\l+1} \le (1+\eta) m_\l \le (1+\eta)n_k$
(see \equ(9.20aa)). Therefore the second term in the right hand side
of \equ(9.35) is at most
$$
\align
&P\{H(m_{\l+1},u) - H(n_k,u) \le n_k [(1+\eta)\big(\la(u)+\de/2\big)
 - \big(\la(u)+\de\big)] \le
-n_k\de/4\}\\
&\le K_3n_k^{-K}+K_6\exp[- K_7n_k \de^2/(\eta + \de)],
\endalign
$$
provided
$$
\eta < \min\Big\{\frac{C_4}{8C_1},
\frac{K_2\de}{2\la(u)},\frac \de{4(\la(u) +\de/2)}\Big\}.
\teq(9.35bc)
$$
It follows that under this last condition
$$
\align
\sum_{k=0}^\infty
&P\{H(n_k,u) > n_k(\la(u) + \de)\} \\
&\le \sum_{k=0}^\infty
P\{H(m_{\l(k)+1},u) > m_{\l(k)+1}(\la(u)
+ \de/2)\} + O(1).
\endalign
$$
The right hand side here is finite by virtue of \equ(9.34),
because $m_{\l(k)} = K_2n_{\l(k)}
< n_k \le K_2n_{\l(k)+1}$ forces $|\l(k) -k| \le K_8$ for some
$K_8$ which is independent of $k$ (see \equ(3.72cba)).
Finally, we may drop the condition \equ(9.35bc), because if $\eta$
does not satisfy this condition, but $\eta'$ does satisfy this
condition, then we may choose $\{n_k(\eta')\}$ so that it contains
the tail of $\{n_k(\eta)$\}, by Corollary 5. By this inclusion and
by what we just proved
$$
\align
&\sum_{k=0}^\infty
P\{H(n_k(\eta),u) > n_k(\eta) (\la(u) + \de)\}\\
&\quad \le \sum_{k=0}^\infty P\{H(n_k(\eta'),u) > n_k(\eta')(\la(u) + \de)\}
+O(1)< \infty.
\tag "$\blacksquare$"
\endalign
$$

\enddemo

\subhead
{\bf 5. Proof of the shape theorem}
\endsubhead
\numsec=5
\numfor=1
Now that we have shown that the spread of the $B$-particles in the full
space process has a definite speed in each direction, the half-space
processes are no longer of importance. In fact {\it Corollary 8 contains
Theorem 1 in the one-dimensional case} (with $B_0 = [-\la(e_1), \la(e_1)]$).
For the higher dimensional case, we shall show in this section how
to go from the existence of $\limt (1/t)H(t,u)$ for all $u \in S^{d-1}$ to the
full shape theorem. This should work for a fairly general class of
processes. The idea to derive the shape theorem via results on the
propagation of half-spaces we learned from \cite {GG}. However, the
details in our case differ from those in \cite {GG}.

The remaining problem in dimension $d > 1$
is that even if we know that $H(t,u)$
grows at rate $\la(u)$, it only tells us that there
exist $B$-particles at time $t$ at some random site $x_t$
for which $\langle x_t, u \rangle \sim t\la(u)$. It
does not tell us where the points $x_t$ near the hyperplane $\{x:
\langle x,u \rangle = t\la(u)\}$ are. In particular, it does not
guarantee that we can find $x_t$ which converge in direction to a
prescribed unit vector, i.e., for given $v \in S^{d-1}$ we do not know
whether we can choose $x_t$ such that $x_t/\|x_t\|_2 \to v$.

To attack this problem we first write down the conjectured limiting
shape $B_0$ in terms of the function $\la(\cdot)$ on $S^{d-1}$. This
conjectured $B_0$ is convex (for trivial reasons). We
then show that we can guarantee $x_t/\|x_t\|_2 \to v$ if $v$
corresponds to a so-called exposed point of the convex set
$B_0$. Using some further properties of convex sets, as well as
approximate convexity properties of the set of points which can be
reached by the
$B$-particles in a large time, we  can then show
that the limiting shape result \equ(1.3) holds.

The convergence result \equ(9.30) suggests that
the limit set $B_0$ in \equ(1.3) should be given by
$$
B_0 = \{z \in \Bbb R^d:\langle z, u \rangle \le \la(u) \text { for all
} u \in S^{d-1}\}.
\teq(4.4)
$$
Clearly this set $B_0$ is a closed convex set. In fact it is also
bounded and hence compact, because $\la(u) \le 2\sqrt d C_1$ for all $u$.
The origin is an interior point of $B_0$ because $\la(u) \ge C_4$.
We call a point $w \in \partial B_0$ an {\it exposed point} of $B_0$
if there exists a supporting hyperplane
$\{z \in \Bbb R^d: \langle a,z \rangle =b\}$ of $B_0$ which contains
$w$, but no other point of $B_0$. Thus
$$
\langle a,w\rangle = b \text{ but } \langle a,z \rangle  < b \text { for
all } z \in B_0 \setminus \{w\}.
\teq(4.5)
$$
Note that this forces $a \ne \bold 0$.
We now show that $\cP^f$
indeed grows in the direction of an exposed point at the
rate which is necessary for \equ(1.3).

\proclaim{Lemma 9}
Let $w$ be an exposed point of $B_0$ and let $(a,b) \in \Bbb
R^d \times \Bbb R$ satisfy \equ(4.5). Let $u = a/\|a\|_2$.
Then, there exists a sequence $\ep_n \downarrow  0$ such that
$$
P\{\cN_n(w, \ep_n) \text{ occurs for all large integers }n\} = 1,
\teq(4.5aa)
$$
where
$$
\align
\cN_n(w,\ep) := \{&\text{in $\cP^f$ there are at
time $\big(1 + 8\ep/C_2\big)n$ occupied }\\
&\text{sites in $nw + \cC(2\ep n)$
and all these sites are in fact}\\
&\text{occupied by $B$-particles
at time }\big(1 + 8\ep/C_2\big)n\}.
\teq(4.6)
\endalign
$$
Also, define
$$
\cO_n(w,\de) = \big\{\text{in $\cP^f$
there is at time $n$ a $B$-particle in $n w +\cC(\de n)$}\big\}.
$$
Finally, let $n_k = n_k(\eta)$ be as in Corollary 5. Then for all $\de
,\eta > 0$
$$
\sum_{k=0}^\infty \big[1- P\{\cO_{n_k(\eta)}(w, \de) \}\big] < \infty.
\teq(4.6aba)
$$
\endproclaim
\demo{Proof} Order the vertices of $\Bbb
Z^d$ in some deterministic way, for instance in the lexicographic way.
Let $x_t$ be the first vertex $x$ in this order which is occupied by a
$B$-particle in $\cP^f$ at
time $t$ and with $\langle x, u \rangle = H(t,u)$.
By \equ(9.30), almost surely,
$$
\frac 1t \langle x_t,u \rangle  \to \la(u) = \limt \frac 1t H(t,u)
\teq(4.6a)
$$
as $t \to \infty$.
Moreover, by \equ(9.31), for each $\de > 0, \eta >0$
$$
\sum_{k=0}^\infty P\big\{\big|\frac 1{n_k}\langle  x_{n_k},u \rangle -
\la(u)\big| > \de \big\} < \infty.
\teq(4.6bcd)
$$

We want to show that for each $\de > 0$
$$
P\Big\{\big\|\frac 1n x_n - w\big\| \le \de
\text{ for all large integers }n\Big\} = 1.
\teq(4.7aa)
$$
Note that $w \in B_0$ implies
$$
\langle w,u \rangle \le \la(u).
\teq(4.7)
$$
Recall next that
$P\{x_n \notin \cC(2C_1n)\}\le K_6n^{-K-d-1}$,
by virtue of \equ(3.39aa) or the estimates for \equ(3.18). So,
$$
P\{x_n \in \cC(2C_1n)\text{ for all large }n\}=1.
\teq(4.7bb)
$$
Also
$$
\sum_{k=0}^\infty P\{x_{n_k} \notin \cC(2C_1n_k)\} < \infty.
\teq(4.7bbb)
$$
So, we can ignore the events $\{x_n \notin \cC(2C_1n)\}$.

Next, let $v \in S^{d-1}$ be a
unit vector which is not a multiple of $w$.
We claim that there exists some $\de = \de(v) > 0$ such that
$$
P\Big\{\Big\|\frac{x_n}{\|x_n\|_2} -v\Big\|< \de \text{ i.o.}\Big\} = 0
\teq(4.9)
$$
and
$$
\sum_{k=0}^\infty P\Big\{\Big\|\frac{x_{n_k}}{\|x_{n_k}\|_2} -v\Big\|
< \de\Big\} < \infty
\teq(4.9a)
$$
(i.o. stands for infinitely often).
To prove this, note first that \equ(4.9) holds if
$\langle v,u \rangle = 0$, because
$$
\liminf_{n \to \infty} \langle \frac{x_n}{\|x_n\|_2} ,u\rangle \ge
\frac {\la(u)}{\limsup_{n \to \infty} \|x_n\|_2/n} \ge \frac {C_4}
  {2\sqrt dC_1}
$$
by virtue of \equ(4.6a), \equ(4.7bb) and the fact the $\la(u) \in
[C_4, 2 \sqrt d C_1]$.
Similarly, \equ(4.9a) holds if $\langle v,u\rangle = 0$, by virtue
of \equ(4.6bcd) and \equ(4.7bbb).
To take care of other vectors
$v$, define for any $y \in \Bbb R^d \setminus \{\bold 0\}$,
with $\langle y,u \rangle \ne 0$,
$$
\wt y = \text{  the unique multiple of $v$ which
satisfies }\langle \wt y,u \rangle = b/\|a\|_2.
$$
In particular, $\wt y$ lies
in the in the  supporting hyperplane $\{z: \langle a,z\rangle
=b\}$. Now, by assumption $\wt v \ne w$, so that $\wt v \notin B_0$. By
definition of $B_0$ this means that there exists some $u' \in S^{d-1}$
such that $\langle \wt v, u'\rangle > \la(u')$. We can then find
$\de>0$ and $\eta >0$
such that $\langle \wt z ,u'\rangle > (1+\eta)\la(u')$ for all $z \in
S^{d-1}$ with $\|z-v\| < \de$.
Thus, if
$$
\Big\|\frac{x_n}{\|x_n\|_2} -v\Big\|< \de,
$$
then
$$
\langle \wt x_n, u' \rangle = \langle \wt {\big(x_n/\|x_n\|_2\big)}, u'
\rangle > (1+\eta) \la(u').
\teq(4.9aa)
$$
In addition, by \equ(4.6a) and \equ(4.7),
$$
\limn \frac 1n \langle x_n,u \rangle = \la(u) \ge \langle w,u \rangle =
\frac b{\|a\|_2} \text{ (see \equ(4.5))},
$$
while, by definition of $\wt y$,
$$
\langle \wt x_n,u \rangle = \frac b{\|a\|_2}.
$$
Moreover, we must have
$$
\|a\|_2\langle w,u \rangle =b > 0
\teq(4.9xy)
$$
by \equ(4.5) and the fact that $\bold
0 \in B_0$.
Consequently, $x_n = \ga_n \wt x_n$ for some reals $\ga_n$ which satisfy
$\ga_n/n \to 1$. Together with \equ(4.9aa)
this would imply
$$
\langle x_n ,u'\rangle > n(1+\eta/2)\la(u')
$$
for large $n$. But, $P\{\langle x_n ,u'\rangle > n(1+\eta/2)\la(u')
\text{ i.o.}\} =0$,
by virtue of \equ(9.30) with $u$ replaced by $u'$ and
the fact that $H(n,u') \ge \langle x_n,u' \rangle$ (by definition
of $H$). Thus \equ(4.9) holds for the chosen $\de$. Similarly, \equ(4.9a)
follows by means of \equ(4.6bcd) with $u'$ instead of $u$.

Now, for any
$\ep > 0$
the compact set
$$
\align
W(\ep) :&= \{z \in S^{d-1}: z = \frac x{\|x\|_2} \text{ for some $x \in
  \cC(2C_1n)$ with }\\
&\phantom{MMMMMMMMMM}
\langle x,u \rangle \ge n\la(u)/2, \|z-\frac w{\|w\|_2}\| \ge \ep \} \\
&=\{z \in S^{d-1}: z = \frac x{\|x\|_2} \text{ for some $x \in
  \cC(2C_1)$ with }\\
&\phantom{MMMMMMMMMM} \langle x,u \rangle \ge \la(u)/2,
\|z-\frac w{\|w\|_2}\| \ge \ep \}
\endalign
$$
is independent of $n$ and is covered by finitely many neighborhoods
$U_1, \dots, U_N$ of the form $U_i = \{z \in S^{d-1}:\|z-v_i\|< \de
(v_i)\}$ with $v_i \in S^{d-1}$.
Thus, by \equ(4.9), $P\{x_n/\|x_n\|_2 \in W(\ep)\text{ i.o.}\} = 0$. This
holds for all $\ep > 0$. In view of
\equ(4.6a) and \equ(4.7bb), this implies
$$
P\Big\{\frac{x_n}{\|x_n\|_2 }\to \frac w{\|w\|_2}\Big\} = 1.
\teq(4.9cba)
$$
In turn, this together with \equ(4.6a) implies
$$
\limn\frac {n\la(u)}{\|x_n\|_2}
= \limn \frac{\langle x_n,u \rangle}{\|x_n\|_2} = \frac {\langle
  w,u\rangle}{\|w\|_2} \text{ a.s.}
$$
Since $\langle w,u \rangle \ne 0$ (see \equ(4.9xy)), $\|x_n\|_2 \sim n\|w\|_2
\la(u)/\langle w,u \rangle$
and
$$
\limn \frac 1n x_n = \frac {\la(u)}{\langle w,u \rangle}w \text{ a.s.}
\teq(4.9ab)
$$

To complete the proof of \equ(4.7aa) we show that
$$
\la(u) = \langle w,u \rangle.
\teq(4.10)
$$
Indeed, we already saw that $\langle \wt x_n,u \rangle = b/\|a\|_2 =
\langle w,u \rangle$. We also saw that $x_n =
\ga_n \wt x_n$ with $\ga_n \sim n$. Therefore  $\langle x_n/n,u
\rangle \sim \langle \wt x_n,u \rangle = \langle  w,u \rangle$. On the
other hand, \equ(4.9ab) implies that $\limn \langle x_n/n,u \rangle
= \la(u)$. Thus \equ(4.10) and \equ(4.7aa) hold.

We now also obtain \equ(4.6aba). Indeed, essentially the same argument
as for \equ(4.9cba), but now using \equ(4.9a) instead of \equ(4.9)
gives
$$
\sum_{k=0}^\infty P\Big\{\Big\|\frac{x_{n_k}}{\|x_{n_k}\|_2} - \frac
w{\|w\|_2}\Big\| > \de \Big\} < \infty.
\teq(4.10aa)
$$
Consequently also
$$
\sum_{k=0}^\infty P\big\{\big|\frac 1{n_k}
\langle x_{n_k},u \rangle - \frac{\|x_{n_k}\|_2}
{n_k \|w\|_2}\langle w,u\rangle \big| > \frac {\de \sqrt d} {n_k}\|x_{n_k}\|_2
\big\} < \infty.
$$
Together with \equ(4.6bcd), \equ(4.10) and \equ(4.7bbb) this last
relation  yields
$$
\sum_{k=0}^\infty P\big\{\big|\la(u)- \frac{\|x_{n_k}\|_2}
{n_k \|w\|_2} \la(u)\big|> \de 3C_1 d \big\} < \infty.
$$
Thus, for a suitable constant $K_7$
$$
\sum_{k=0}^\infty P\Big\{\Big|\frac{\|x_{n_k}\|}{n_k}- \|w\|_2\Big| > K_7
 \de \Big\} < \infty.
$$
Together with \equ(4.10aa) this finally gives for some other constant $K_8$
$$
\align
&\sum_{k=0}^\infty \big[1-P\{\cO_{n_k(\eta)}(w, K_8\de)\}\big]\\
&\le
\sum_{k=0}^\infty P\big\{\big\|\frac{x_{n_k}}{n_k}- w\big\| > K_8\de
 \big\} < \infty.
\endalign
$$
Since this holds for any $\de >0$, this is equivalent to \equ(4.6aba).

The preceding (see \equ(4.7aa))
shows that there exists a sequence $\ep_n \to 0$, and
 random vertices $x_n$ such that with probability 1, for all large $n$,
$$
x_n \in nw + \cC(\ep_nn) \text{ and $\cB^f(x_n,n)$ occurs},
\teq(4.11)
$$
where
$$
\cB^f(x,s) :=  \{\text{there is $B$-particle at $x$ at time $s$ in
 $\cP^f$}\}.
$$
Now take
$$
\wt n := n\big(1 +  \frac{8\ep_n}{C_2}\big)
$$
and define the event
$$
\align
\cR(x,n) = \{&\text{there is some particle in $\cP^f$
which lies in }x + \cC\big(C_2(\wt n
  -n)/2\big)\\
& = x + \cC(4\ep_nn) \text{ but is of type $A$
at time $\wt n$}\}.
\endalign
$$

We shall complete the proof by proving that the event
$$
\align
\{\text{for infinitely many $n$ there }&\text{exists an $x_n$ for which}\\
&\cB^f(x_n,n) \cap \cR(x_n,n) \text{ occurs}\}
\teq(4.12)
\endalign
$$
has probability 0. First we show that this will indeed prove the lemma.
The probability that
any particle which is in $nw + \cC(\ep_nn)$ at time $n$
is outside $nw +\cC(2\ep_n n)$ at time
$\wt n$ is bounded by
$$
K_9[\ep_n n]^d
P\{\sup_{r \le \wt n - n} \|S_r\| \ge \ep_n n\}.
\teq(4.11aa)
$$
Without loss of generality we can let $\ep_n$ go to
 0 so slowly that for large $n$
this expression is no more than $n^{-K-1}$
(by (2.42) in \cite {KSa}) and such that
$$
\ep_n \ge n^{-1/2}.
\teq(4.11bb)
$$
From this  and the
fact that $\cB^f(x_n,n)$ occurs for all large $n$, we conclude via the
Borel-Cantelli lemma that almost surely, for all large $n$ there
are particles in $\cP^f$ in the set $nw + \cC(2\ep_n n)$ at time $\wt n$.
The fact that \equ(4.12) has probability 0
 will then imply that $\cR(x_n,n)$ must
fail for all large $n$.  
But this implies that a.s. there are particles in $\cP^f$ which
lie in $nw + \cC(2 \ep_n n) \subset x_n +\cC(4\ep_n n)$ at time $\wt
n$,
and all of these particles must have type $B$.
This is the desired result \equ(4.5aa).

It remains to prove \equ(4.12). But this is almost immediate from
Proposition 1. Indeed,
$$
\align
&P\{\cB^f(x_n,n)\cap \cR(x_n,n)\}\\
& \le P\{x_n \notin \cC(2C_1n)\}
+\sum_{x \in \cC(2C_1n)}P\{\cB^f(x,n) \text{ but there is a particle
in $\cP^f$
of}\\
&\phantom{MMMMMMMMMMMM}
\text{ type $A$ at some $z \in x + \cC(4\ep_n n)$ at time $\wt n$}\}\\
&\le K_6n^{-K-d-1} + \sum_{x \in \cC(2C_1n)}
P\{\text{$x$ is occupied at time $n$ in $\cP^f$ and
in the full-space}\\
&\phantom{iMMMMMMMMMMMM}\text{process started at $(x,n)$
there is an $A$-particle}\\
&\phantom{iMMMMMMMMMMMM}
\text{at some $z \in x + \cC(4\ep_n n)$ at time $\wt n$}\},
\endalign
$$
where we used Lemma C for the last inequality.
As in the estimate for $\cK_2$ in \equ(6.5), by
\equ(3.0c) and \equ(2.3) with $K$ replaced by $2K+2d$,
the last sum here is at most
$$
 K_{10}n^d (\text{left hand side of \equ(2.3) with $t =
\wt n-n = 8\ep_n n/C_2 $}) \le K_{11} n^{-K}
$$
(see \equ (4.11bb)).
\hfill $\blacksquare$
\enddemo

The preceding lemma shows that the set $B(t)$ grows in the direction
of the exposed points of $B_0$ in $\partial B_0$ at the ``right''
speed. More specifically, if $w$ is such a point, then almost surely,
for all large $t$, there exist points $w(t) \in (1/t)B(t)$ such that
$w(t) \to w$. We merely have to choose  $n$ in Lemma 9 such that
$n(1+8\ep_n/C_2) \le t$ but $n/t \to 1$, and then $w(t)$ a point in
$\wt B(n) \cap [nw +\cC(2\ep_nn)]$. Lemma 9 guarantees that this last
intersection is nonempty for large $n$. The next two lemmas will show
that the same is true for any point $w \in \partial B_0$. This is
basically done by  concatenating a number of paths which
produce $B$-particles at $\al_i n w_{n,i}$ for exposed points
$w_{n,i}$ and $\sum_{i=1}^k\al_i w_{n,i} \to w, \al_i \ge 0,
\sum_{i=1}^k \al_i = 1$. Lemma 10 contains the basic technical
step. It explains how the concatenation works; this is basically the
same construction as in the proof of Lemma 7.

\proclaim{Lemma 10} Let $w_1, w_2 \in \partial B_0$.
Assume that there exist $\ep_n > 0$ such that $\ep_n \to 0$ and
such that \equ(4.5aa)  holds with $w$ replaced by $w_1$, that is,
$$
P\{\cN_n(w_1,\ep_n) \text{ occurs for all large integers }n\} = 1.
\teq(4.53cc)
$$
{\rm (}We are not
assuming that $w_1$ is an exposed point of $B_0$.{\rm)}
In addition, assume that for all $\de, \eta >0$
$$
\sum_{k=0}^\infty \big[1- P\{\cO_{n_k(\eta)}(w_2, \de) \}\big] < \infty
\teq(4.53)
$$
(see Corollary 5 for $n_k = n_k(\eta)$). Let $0 < \al < 1$ and $\eta >
0$. Then there exist $\de_n > 0$ such that $\de_n \to 0$ and such that
$$
P\{\cN_n(\al w_1 +(1-\al)w_2,\de_n) \text{ occurs for all large }n\} = 1.
\teq(4.54)
$$
\endproclaim
\demo{Proof} Fix $0 < \al < 1$. Also fix
$$
\de > 0 \text{ and } 0 < \eta < \de/2
$$
for the time being. Take
$$
p_k =p_k(\eta) = \Big \lfloor \frac \al{1-\al} n_k(\eta)\Big\rfloor
$$
and
$$
q_k = q_k(\eta) = \big(1+ 8 \ep_{p_k(\eta)}/C_2\big)p_k.
$$
Define
$\cO'_{n_k}(w_2, \de)$ as the translate by $\big (p_k(\eta)w_1,
q_k(\eta)\big) $ (in space-time) of
$\cO_{n_k}(w_2,\de)$.
Explicitly,
$$
\align
\cO'_{n_k}(w_2,\de) =\{&\text{in the full-space process started at
$\big (p_k(\eta)w_1,q_k(\eta)\big)$ there}\\
&\text{is at time $
q_k+n_k$  a $B$-particle in }p_kw_1+n_kw_2
+\cC(\de n_k)\}.
\endalign
$$
(We suppress the dependence on $w_1$ and $\eta$ in this notation).
Also let
$$
z_k= \text{nearest occupied site to $p_kw_1$ in $\cP^f$ at time }q_k.
$$
Since $P\{\cO'_{n_k}(w_2, \de)\} = P\{\cO_{n_k}(w_2, \de)\}$,
assumption \equ(4.53) implies that almost surely,
$$
\cO'_{n_k}(w_2,\de) \text{ occurs for all large }k.
\teq(4.55)
$$
Also, by assumption \equ(4.53cc), almost surely,
$$
\cN_{p_k}(w_1, \ep_{p_k})  \text{ occurs for all large }k.
\teq(4.56)
$$
Now consider a $k$ for which $\cN_{p_k}(w_1,\ep_{p_k})\cap
\cO'_{n_k}(w_2,\de)$
occurs. By the definition of $\cN_{p_k}$ this implies that $z_k$ lies
in $p_k w_1 +\cC(2\ep_{p_k}p_k)$ and that the particles at $z_k$
at time $q_k$ have type $B$ in $\cP^f$. Therefore the resetting of the types to
start the full-space process at $(p_kw_1, q_k)$ does not change the
type at $z_k$. By the monotonicity  property of Lemma C, $\cP^f$
therefore has at least as many $B$-particles
 at any space-time point $(x,t)$ with $t \ge q_k$ as
the full state process started at $(p_kw_1,q_k)$. Since $\cO'_{n_k}(w_2,
\de)$ occurs this implies that in $\cP^f$ there
is a $B$-particle in $p_kw_1+n_kw_2 +\cC(\de n_k)$ at time $q_k+n_k$.

Let the nearest $B$-particle to $p_kw_1+n_kw_2$ in $\cP^f$ at time
$q_k+n_k$ be at the position $y_k$, so that $\cB^f(y_k,q_k+n_k)$ occurs.
The last paragraph gives us that $\|y_k-p_kw_1 -n_kw_2\| \le \de n_k$.
These are only statements for the times $q_k+n_k$. Since
\equ(4.53cc) requires that certain events happen for all large $n$ we
now first show how to go from the $q_k +n_k$ to general integers
$n$. For any large $n$ let $k(n)$ be such that
$q_k + n_k \le n < q_{k+1}+ n_{k+1}$.
Then for large $n$
$$
q_k + n_k \le n \le (q_k+n_k)(1+2\eta)  \le (q_k+n_k)(1+\de),
$$
since $n_{k+1}/n_k \le 1+\eta$. Also by our choice of $p_k,q_k$
$$
\|p_kw_1 + n_kw_2 - n[\al  w_1 + (1-\al) w_2]\| \le K_{12}\de  n.
$$
Thus, on $\cB^f(y_k, q_k+n_k)$, there is a $B$-particle at $y_k \in n[\al w_1 +
(1-\al)w_2] + \cC\big((K_{12} +1)\de n\big)$ at time $q_k+n_k$.
Moreover, as in \equ(4.11aa) we have
$$
\align
&P\{\text{in $\cP^f$ there is a $B$-particle in
$n[\al w_1 +
(1-\al)w_2] + \cC\big((K_{12} +1)\de n\big)$}\\
&\phantom{Mi}\text{at time $q_k+n_k$ which is no
longer in $n[\al w_1 +
(1-\al)w_2] + \cC\big((K_{12} +2)\de n\big)$}\\
&\phantom{Mi}\text{at time $n$}\}\\
&= O\big(n^{-K}\big).
\endalign
$$
Thus, almost surely, there is in $\cP^f$ for all large $n$ a $B$-particle in
$n[\al w_1 +
(1-\al)w_2] + \cC\big((K_{12} +2) \de n\big)$ at time $n$.
We can now proceed as in Lemma 9.
Essentially as in \equ(4.12) and in the lines following it we now have
that almost surely
$$
\align
\{&\text{there is some $y \in
n[\al w_1+ (1-\al)w_2] +\cC\big((K_{12}+2)\de n\big)$ for which}\\
&\text{$\cB^f(y,n)$ occurs, but
in $\cP^f$ there are either no particles or an $A$-particle in}\\
&\text{$n[\al w_1+ (1-\al)w_2] +\cC\big(2(K_{12}+2 )\de n \big)$
at time } \big(1+ 8(K_{12}+2) \de/C_2\big)n\big\}
\teq(4.61)
\endalign
$$
occurs only for finitely many $n$. This shows that
$$
P\{\cN_n\big(\al w_1+(1-\al)w_2,(K_{12} +2)\de\big)
\text{ occurs for all large }n\} = 1.
\teq(4.61abc)
$$
This holds for all $\de >0$ and $\eta < \de/2$. However, \equ(4.61abc)
is already independent of $\eta$, so that it holds for all $\de > 0$.
There then also exists a sequence $\de_n \to 0$ such that
$\cN_n\big(\al w_1 + (1-\al)w_2, \de_n\big)$
occurs almost surely for all large $n$.
\qedsymbol
\enddemo

{\bf Proof of Theorem 1.} We shall prove \equ(1.3) with the $B_0$
defined in \equ(4.4). For the right hand inclusion in \equ(1.3) we
note that for any $\ep > 0$ there exists finitely many halfspaces
$ \{z \in \Bbb R^d: \langle z,u_i\rangle
\le \la(u_i)\}, \;1 \le i \le
N$, with $u_i \in S^{d-1}$ such that
$$
\bigcap_{i=1}^N \{z \in \Bbb R^d: \langle z,u_i\rangle
\le \la(u_i)\}
\subset (1+\ep/3)B_0.
\teq(4.62)
$$
Indeed, $B_0$ is contained in the cube $\wt \cC
:=\cap_{i=1}^N \{z \in \Bbb R^d:-\la(e_i) \le
 \langle z,e_i \rangle \le \la(e_i)\}$ (with $e_i = i$-th coordinate vector),
and by compactness,
$\wt \cC \setminus$(interior of $(1+\ep/3)B_0$ is covered by finitely
many relatively open subsets of $\wt\cC$ of the form $\wt \cC \cap
\{z \in \Bbb R^d:
\langle z,u \rangle > \la (u)\}$. In addition to \equ(4.62)
we know from \equ(9.30) that, almost surely, $H(t,u_i)
< t(1+\ep/3) \la(u_i)$ for all large $t$ and $i=1, \dots,
N$. Consequently, almost surely
$$
\wt B(t) \subset t (1+\ep/3) \bigcap_{i=1}^N \{z \in \Bbb R^d:
\langle z,u_i\rangle  \le \la(u_i)\}
\subset (1+\ep/3)^2 t B_0
$$
for all large $t$. Thus the right hand inclusion
in \equ(1.3) holds.

For the left hand inclusion in \equ(1.3) we first observe that by
Lemma 9, the hypotheses \equ(4.53cc)
and \equ(4.53) of Lemma 10 hold for
all exposed points $w_1,w_2 \in \partial B_0$. It then follows from Lemma 10
that \equ(4.54) holds. In turn, \equ(4.54) states that the hypothesis
\equ(4.53cc) with $w_1$ replaced by $\al w_1 +(1-\al)w_2$ is
satisfied. Therefore,
if $w_3 \in \partial B_0$ is also an exposed point of $B_0$ and $0 <
\be < 1$, then we get from Lemma 10 that there exist $\de'_n \to 0$
such that
$$
P\{\cN_n(\be \al w_1 +\be(1-\al)w_2+(1-\be)w_3, \de'_n)
\text{ occurs for all large }n\} = 1.
$$
But as $\al$ and $\be$
vary over $(0,1),\; \be \al w_1 +\be(1-\al)w_2+(1-\be)w_3$ varies over
the convex combinations $\al_1 w_1 + \al_2w_2 + \al_3 w_3$ with $\al_i
> 0, \sum_{i=1}^3 \al_i = 1$. We can repeat this procedure to obtain
that for each convex combination $\sum_{i=1}^k \al_i w_i$ with $\al_i
\ge 0, \sum_{i=1}^k \al_i = 1$ and $w_i \in \partial B_0$ exposed
points of $B_0$, there exist $\de_n \to 0$ such that
$$
P\{\cN_n\big(\sum_{i=1}^k \al_iw_i , \de_n\big)
\text{ for all large integers }n\}= 1.
$$
In particular (see \equ(4.6)), for each such $\sum_{i=1}^k \al_i w_i$
and each fixed $\eta > 0$
$$
\align
P\{&\text{in $\cP^f$ there are at time $(1+8\de_n/C_2)n\; B$-particles in}\\
&\text{$n \sum_{i=1}^k \al_iw_i + \cC(2\eta n)$ for all large integers }n\}=1.
\endalign
$$
In turn, this means that if for a given vector $v$ and $\eta > 0$ we
can find $\al_i, w_i$ as above such that $\|v - \sum_{i=1}^k \al_iw_i\|
\le \eta$, then also
$$
\align
P\{&\text{in $\cP^f$ there are at time $(1+8\de_n/C_2)n \;B$-particles in}\\
&\text{$n v + \cC(3\eta n)$ for all large integers }n\} = 1.
\teq(4.63)
\endalign
$$
If $v$ is such that there exist $k^{(r)} < \infty,
\al_i^{r} \ge 0$ and $ w_i^{(r)} \in
\partial B_0$ exposed points of $B_0$ such that $\sum_{i =1}^{k^{(r)}}
\al_i = 1$ and $\|v-
\sum_{i=1}^{k^{(r)}} \al_i^{(r)} w_i^{(r)}\| \to 0$ (as $r \to \infty$),
then \equ(4.63) holds for each $\eta > 0$. For such $v$ there then
exist $\eta_n \to 0$ such
that almost surely, for all large $n$ there exist $B$-particles within
distance $4\eta_n n$ of $nv$ at time $(1+8\de_n/C_2)n$, for some $\de_n
\to 0$ ($\de_n$ and $\eta_n$ may depend on $v$).

The last statement applies to each $v \in B_0$, because each such $v$
is a convex combination of at most $(d+1)$ extreme points of $B_0$
(see \cite{Ru}, Theorem 3.22 and Lemma following Theorem 3.25)
and the exposed points of $B_0$ are dense
in the extreme points (Strascewicz' theorem; see Theorem 18.6 in \cite{Ro}).
Thus, by applying the last result to a fixed
$v \in B_0$ with $n = \lfloor (1-\ep)t \rfloor$
and $0 < \ep < 1$, we find that
almost surely for all large $t$,
$$
\align
&\text{at time $(1+8\de_n/C_2)n$ there exists a site $v_n$ with}\\
&\text{$\|v_n-nv\| \le 4 \eta_n n$,
which is occupied in $\cP^f$ by $B$-particles}.
\teq(4.64)
\endalign
$$
We claim that
$$
\align
&P\{\text{\equ(4.64) holds, but not all sites in $(1-\ep)tv + \cC(C_2 \ep t/4)$
belong to $\wt B(t)$}\}\\
& \le K_{13}t^{-K}.
\teq(4.65)
\endalign
$$
This is an easy consequence of \equ(3.0c) and  Theorem A. Indeed,
from \equ(3.0c) with $(X,s)$ taken to be $\big(v_n,(1+8\de_n/C_2)n\big)$
and
$$
\align
\cA &= \{\text{not all vertices in $\cC(C_2\ep n/2)$ have been visited
by a $B$-particle by time $\ep n/2$}\}\\
& =\{\cC(C_2 \ep n/2) \not
\subset B(\ep n/2)\},
\endalign
$$
we see that the probability in \equ(4.65) is for large $t$ at most
$$
K_{14} n^d P^{or} \{\cC(C_2 \ep n/2) \not
\subset B(\ep n/2)\} \le K_{15} n^{-K} \le K_{13} t^{-K}
$$
(for the first inequality here we used Theorem A with $K+d$
in the place of $K$).
This establishes the claim \equ(4.65).

To obtain Theorem 1 we now choose for a given $\ep$ a finite number of
vectors $v^{(1)}, \dots, v^{(N)}$ in $B_0$ such that each $v \in B_0$
satisfies $\|v-v^{(r)}\| < C_2 \ep/4$ for at least one $r$.
This means that
$$
B_0 \subset \bigcup_{1 \le r \le N} \big[v^{(r)} + \cC(C_2 \ep/4)\big].
$$
Moreover, by \equ(4.64)
and \equ(4.65) it holds almost surely for all large $t$ that
$$
\bigcup_{1 \le r \le N} \big[(1-\ep)t v^{(r)} + \cC(C_2\ep t/4)\big]
\subset \wt B(t).
$$
Together, these last two inclusions imply that
almost surely the left hand inclusion in \equ(1.3) holds
for all large $t$.

\enddocument